\documentclass{article}
\usepackage[utf8]{inputenc}
\usepackage{url}
\usepackage{hyperref}
\usepackage{graphicx}
\usepackage{amsmath,amssymb}
\usepackage{xcolor}
\allowdisplaybreaks
\usepackage{float}
\usepackage{graphicx}
\usepackage{appendix}
\usepackage{listings}
\usepackage{subfigure}
\usepackage{cite}
\usepackage{bm}
\usepackage{makecell}
\usepackage{multirow}
\usepackage{algorithm}
\usepackage{algorithmicx}
\usepackage{algpseudocode}
\usepackage[T1]{fontenc}
\usepackage[utf8]{inputenc}
\usepackage{authblk}

\usepackage{multicol}
\newcommand\keywords[1]{\textbf{Keywords}:#1}
\newtheorem{theorem}{Theorem}[section]

\newtheorem{remark}{Remark}[section]
\newtheorem{corollary}{Corollary}[section]

\makeatletter
\newenvironment{breakablealgorithm}
  {\begin{center}
     \refstepcounter{algorithm}
     \hrule height.8pt depth0pt \kern2pt
     \renewcommand{\caption}[2][\relax]{
       {\raggedright\textbf{\ALG@name~\thealgorithm} ##2\par}
       \ifx\relax##1\relax
         \addcontentsline{loa}{algorithm}{\protect\numberline{\thealgorithm}##2}
       \else
         \addcontentsline{loa}{algorithm}{\protect\numberline{\thealgorithm}##1}
       \fi
       \kern2pt\hrule\kern2pt
     }
  }{\kern2pt\hrule\relax
   \end{center}
  }
\makeatother

\newcommand{\bd}{\bm d}
\newcommand{\bx}{\bm x}

\newcommand{\bX}{\bm X}

\newcommand{\bu}{\bm u}

\newcommand{\bv}{\bm v}

\newcommand{\bV}{\bm V}

\newcommand{\bpsi}{\bm \psi}

\newcommand{\bxi}{\bm \xi}

\title{A time-segmented non-intrusive reduced-order modeling for
fluid-structure interactions using greedy algorithms}
\author[c,a]{Qijia Zhai \thanks{zhaiqijia@163.com}}
\author[a]{Shiquan Zhang \thanks{shiquanzhang@scu.edu.cn}}
\author[b]{Pengtao Sun \thanks{pengtao.sun@unlv.edu}}
\author[a]{Xiaoping Xie \thanks{Corresponding author: xpxie@scu.edu.cn}}
\affil[a]{School of Mathematics, Sichuan University, Chengdu 610064, China}
\affil[b]{Department of Mathematical Sciences, University of Nevada Las Vegas, Las Vegas, Nevada 89154, USA}
\affil[c]{Shenzhen Research Institute of Big Data, Shenzhen, 518116, China}

\date{}
\begin{document}

\maketitle
\begin{abstract}
In this paper, a novel non-intrusive reduced order modeling (ROM)
approach is developed for solving fluid-structure interaction (FSI)
problems by combining with the greedy algorithms of greedy parameter
selection and orthogonal greedy training technique. This method
employs a time-segmentation technique to handle long-term
simulations, applying differentiated processing to characteristics
within distinct time windows. During each time segment, a greedy
parameter selection strategy is developed to iteratively identify
the most informative ``physical parameter--time'' pairs by
maximizing projection errors onto the current ROM bases. The
selected snapshots are processed using Gram-Schmidt
orthogonalization to construct orthonormal ROM bases for each
variable separately. Furthermore, a shallow neural network trained
using the orthogonal greedy algorithm is employed to establish the
relationship between physical parameters and basis coefficients of
ROM, which not only provides theoretical convergence guarantees but
also avoids common issues such as overfitting and entrapment in
local minima. The proposed method is validated through comprehensive
numerical experiments using a benchmark FSI problem that involves a
vibrating elastic beam. Numerical results illustrate that the
developed ROM approach significantly enhances computational speed
while maintaining accuracy comparable to the full-order model,
effectively handles physical parameters' perturbations during the
online phase and exhibits dimensionality independence, facilitating
its extension to high-dimensional FSI problems. This study
demonstrates that the integration of greedy parameter selection with
convergent machine learning methods marks a significant advancement
in non-intrusive surrogate modeling for complex multiphysics
systems.
\end{abstract}

\section{Introduction}
Fluid-structure interaction (FSI) problems concern the coupled
behavior of fluid flow and solid motion or deformation: fluid forces
influence the structural response, which in turn modifies the flow.
This interplay is central to cardiovascular dynamics, where blood
flow interacts with deformable tissues. In engineering, it shapes
the response of bridges and tall buildings to wind, the motion of
parachutes, and the behavior of turbine blades and fluid-driven
rotating machinery \cite{DowellMODELING}. Despite their distinct
physical settings, these applications share a common modeling and
simulation challenge: capturing the mutual influence of fluid
dynamics and structural mechanics. Accurate simulation of FSI
problems enables comprehensive analysis of system performance and
durability under complex operating conditions. However,
high-fidelity FSI modeling presents significant computational
challenges due to the need to couple distinct fluid and structural
solvers.

Numerical methods for solving FSI problems can be categorized into
two main approaches based on mesh conformity and coupling
strategies. In terms of mesh conformity, the arbitrary
Lagrangian-Eulerian (ALE) method
\cite{Hughes1981,Huerta1988,Nitikitpaiboon1993,Souli2010} represents
the most common body-fitted approach, where the fluid mesh adapts to
structural deformations, naturally satisfying interface conditions.
Alternatively, the fictitious domain/immersed boundary method
(FD/IBM) \cite{Peskin1972,Peskin2002,LeVeque1994,ZLi2001} employs a
fixed background Eulerian mesh with a moving Lagrangian structural
mesh, enforcing interface conditions throughout the structural
domain. Regarding coupling strategies, the monolithic approach
\cite{WANG;Yan;ZHANG2012,Ryzhakov2010a,langer2016robust,Barker.A;Cai.X2009a}
solves fluid and structural equations simultaneously, which offers
unconditional stability and immunity to systematic errors while
implementing interface conditions. The partitioned approach
\cite{grandmont_2001_numerical,causin2005added,tuprints254,Hou2012}
solves fluid and structural equations separately through
Dirichlet/Robin-Neumann alternating iterations and may fail to
converge when fluid and structural densities are comparable due to
the added-mass effect \cite{Idelsohn_Pin2009,causin2005added}.

In this paper, we adopt the ALE-finite element method (ALE-FEM)
within the monolithic framework as our full-order model (FOM),
following the methodology in
\cite{Wang;Sun2016,HaoSun2021,SunLan2020,SunLan2019}. 
While the FOM
provides high accuracy, it requires large-scale computations over
extended time periods, resulting in time-consuming simulations with
substantial computational costs in practical applications.

To address the computational challenges of FSI simulations, the
reduced-order modeling (ROM) have been developed to construct
low-dimensional surrogate models that can significantly reduce
computational complexity and cost while maintaining accuracy
comparable to that of FOM 
\cite{yang2016active,li2019numerical,kalashnikova_2013_a,ballarin_2016_podgalerkin,ballarin2017reduced20,nonino2019overcoming18}.

Traditional ROM for fluid dynamics employ Galerkin projection of the
Navier-Stokes equations onto spaces spanned by time-invariant proper
orthogonal decomposition (POD) bases, resulting in reduced ODE
systems with respect to spatial variables
\cite{alfioquarteroni_2014_reduced}. For FSI problems, ROM methods
have been developed using both partitioned and monolithic
approaches, where separated ROMs are typically constructed for fluid
and structural domains and then coupled together
\cite{beran_2004_reducedorder,vierendeels_2007_implicit,kalashnikova_2013_a,alfioquarteroni_2014_reduced}.
When using ALE-based body-fitted mesh methods, the moving fluid mesh
requires multi-POD approaches that account for grid shifts
\cite{anttonen_2003_podbased}. A notable advancement in intrusive
ROM for FSI problems is the projection-based time-segmented reduced
order model \cite{zhai2025projection}, where a novel ROM is
developed for dynamic FSI problems using the ALE-FEM in a monolithic
frame. In this approach, variables are spatially separated by
attributions (fluid/structure), components (velocity/pressure) and
categories (2D/3D), while POD bases are temporally constructed in
deliberately partitioned time segments, which enables prolonged
simulations under prescribed accuracy thresholds and achieves
comparable accuracy to the FOM with much higher efficiency.

However, intrusive ROM approaches face several limitations,
including the need to access and modify the source code, advance
the reduced equations sequentially in time, maintain the stability
of the reduced dynamics, and repeatedly evaluate nonlinear terms
\cite{alfioquarteroni_2014_reduced}. These
challenges have motivated the development of non-intrusive reduced
order modeling (NIROM) methods
\cite{audouze2009non,guo2018non,guo2019data}, which are independent
from source code modifications and physical system implementations.


The construction of NIROM typically follows a two-step procedure.
The first step involves the identification of appropriate ROM bases
to approximate the solution, the numerical methods for which include
proper orthogonal decomposition (POD)
\cite{sirovich_1987_turbulence,lumley1967structure},
tensor train decomposition (TTD) \cite{swischuk2019projection},
dynamic mode decomposition (DMD) \cite{chen2020machine}, and polynomial chaos expansion (PCE) \cite{peherstorfer2018survey}. The second step focuses on establishing the relationship between input parameters and ROM basis coefficients through regression techniques such as Gaussian Process Regression \cite{raissi2019physics} or machine learning approaches \cite{chen2018neural,lu2019deepxde}. 

Several NIROM strategies have been developed for FSI problems by
combining snapshot-based dimensionality reduction with data-based
approximation of reduced quantities.
For instance, in \cite{xiao2016non} 
the reduced bases are directly constructed using proper orthogonal
decomposition (POD) of high-fidelity FSI simulation snapshots, and
the radial basis function (RBF) interpolation is employed to
approximate the evolution of the reduced coefficients.
The POD-RBF framework has been extended to compressible fluid and
fractured solid coupling applications \cite{xiao2017non},
demonstrating its versatility across diverse FSI scenarios. More
recently, a different non-intrusive approach combining tensor-train
decomposition (TTD) with polynomial chaos expansion (PCE) is
proposed for uncertainty quantification in parameterized FSI
problems \cite{xu2025reduced}, emphasizing the method's capability
to handle parametric variations without source code modifications.
In related work on parametric NIROM, Gaussian process regression is
introduced to approximate the discrepancy between full-order and
reduced-order solutions \cite{xiao2019error}, where machine learning
serves as an error-estimation and correction mechanism supplementing
an interpolation-based reduced model. Taken together, these
developments highlight the potential of data-driven NIROM for
efficient FSI simulation and parametric analysis, particularly when
modifying the source code is impractical or repeated evaluations
over a range of parameter values are required.

Despite these advances, current NIROM approaches exhibit several
limitations. The development of fast and accurate ROM necessitates
careful consideration of error estimation and accuracy assessment
strategies. The construction of ROM usually requires numerous
simulations of the high-fidelity full model across different
parameter configurations, and the existing methods frequently employ
uniform parameter sampling strategies, which can not capture the
most informative regions of the parameter space and result in
suboptimal ROM. Once constructed, accurate assessment of ROM
performance through various error estimation methods becomes crucial
\cite{quarteroni2015reduced,hesthaven2018non}. Another significant
gap remains for non-intrusive surrogate modeling is that there is no
existing method combines greedy algorithms with convergent machine
learning methods. While greedy algorithms have been proven effective
for physical parameters' selection in ROM construction
\cite{quarteroni2015reduced}, and machine learning methods have
shown promise for approximating complex parameter-to-solution
mappings \cite{hesthaven2018non}, their integration remains largely
unexplored. Thus, combining advanced machine learning techniques
with greedy algorithms holds the promise of achieving more efficient
and accurate non-intrusive ROM approach, thereby fully leveraging
the strengths of both methodologies.

In this paper, we propose a non-intrusive surrogate model that
combines greedy parameter selection with orthogonal greedy training
for FSI problems, while considering distinct physical features of
FSI during different time segments. In particular, for a long-term
FSI simulation, the entire time interval is divided into multiple
time segments according to significantly different physical
characteristics, then the ROM basis is constructed separately in
each time segment, where a greedy selection strategy is developed as
well to iteratively identify the most informative ``physical
parameter--time'' pairs by maximizing projection errors onto the
current reduced bases. The selected snapshots are then processed
using Gram-Schmidt orthogonalization to construct orthonormal ROM
bases. To establish the relationship between parameters and ROM
basis coefficients, we employ shallow neural networks trained using
the orthogonal greedy algorithm (OGA), which provides theoretical
convergence guarantees and avoids the common issues of overfitting
and local minima encountered in traditional neural network training.
The proposed approach is applied to each variable (fluid velocity
and pressure, structural displacement, ALE displacement),
separately, enabling the variable-specific optimization and accuracy
improvement as well. Through comprehensive numerical experiments on
a benchmark FSI problem which involves a vibrating elastic beam, we
demonstrate that the developed NIROM approach achieves comparable
accuracy to the FOM while providing significant computational
speedup, effectively handles physical parameters' perturbations in
the online phase and exhibits dimensionality independence as well,
facilitating its extension to high-dimensional FSI problems.

The structure of this paper is organized as follows. In Section
\ref{S2}, a generic FSI model and its fully discrete ALE-FEM are
introduced. The novel non-intrusive surrogate modeling approach
using greedy algorithms is proposed in Section \ref{S3}. Numerical
experiments are carried out in Section \ref{S4} to validate the
accuracy and efficiency of the proposed method. Finally, concluding
remarks are given in Section \ref{s5}.

\section{Mathematical model and ALE-FEM/full-order modeling of FSI}\label{S2}

\begin{figure}[H]
    \centering
    \includegraphics[width = 7cm, height = 3cm]{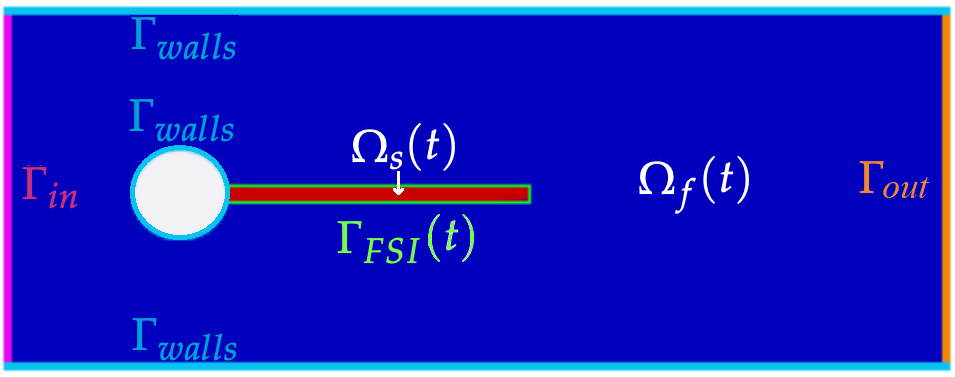}
    \caption{A schematic domain of benchmark FSI problem \cite{turek2006proposal, turek2010numerical}
    with the fluid channel $\Omega_f(t)$ in blue and an elastic beam
    $\Omega_s(t)$ (the structural domain) behind a rigid cylinder in red.
    From the left to right, lines in purple, cyan, green and orange
    represent the inlet $\Gamma_{in}$, the fluid channel wall $\Gamma_{walls}$,
    the fluid-structure interface $\Gamma_{FSI}(t)$, and the outlet
    $\Gamma_{out}$, respectively.}    \label{domainsetting}
\end{figure}

Consider an open bounded domain $\Omega \subset \mathbb{R}^d$
($d=2,3$) characterized by a polygonal boundary $\partial\Omega$.
The domain is partitioned by a time-dependent interface
$\Gamma_{FSI}(t)$ into two distinct subdomains: $\Omega=\Omega_f(t)
\cup \Omega_s(t)$ and $\Gamma_{FSI}(t)=\bar{\Omega}_f(t) \cap
\bar{\Omega}_s(t)$ for $t\in [0,T]$, where $\Omega_f(t)$ is the
fluid domain and $\Omega_s(t)$ the structural domain, and
$\Omega_f(t) \cap \Omega_s(t)=\varnothing$. Figure
\ref{domainsetting} illustrates a domain of a benchmark FSI problem
\cite{turek2006proposal, turek2010numerical}, where
$\partial\Omega=\Gamma_{in}\cup\Gamma_{walls}\cup\Gamma_{out}$ with
$\Gamma_{in}$, $\Gamma_{walls}$ and $\Gamma_{out}$ representing the
inlet, walls and outlet of the fluid channel, respectively.

\subsection{Model description of FSI}
In this paper, we model the fluid motion within FSI problems using
the Newtonian-type incompressible Navier-Stokes equations in terms
of the fluid velocity $\boldsymbol{u}_f$ and fluid pressure $p_f$,
defined as: find $(\boldsymbol{u}_f, p_f) \in
H^1(0,T;H^2(\Omega_f(t))^d) \times L^2(0,T;H^1(\Omega_f(t)))$ such
that
\begin{equation}\label{nseq}
\begin{cases}
\rho_f\left(\frac{\partial \boldsymbol{u}_f}{\partial t} + (\boldsymbol{u}_f \cdot \nabla) \boldsymbol{u}_f\right) - \operatorname{div} \boldsymbol{\sigma}_f(\boldsymbol{u}_f, p_f) = \boldsymbol{b}_f, & \text{in } \Omega_f(t) \times (0,T], \\
\operatorname{div} \boldsymbol{u}_f = 0, & \text{in } \Omega_f(t)
\times (0,T],
\end{cases}
\end{equation}
which is subject to the following boundary conditions and initial
conditions:
\begin{equation}\label{fluid-BICs}
\begin{array}{rclll}
\boldsymbol{u}_f &=& \boldsymbol{u}_{in}, & \text{on } \Gamma_{in}, & \text{(the inlet flow condition)} \\
\boldsymbol{u}_f &=& \boldsymbol{0}, & \text{on } \Gamma_{walls}, & \text{(the no-slip condition)} \\
\boldsymbol{\sigma}_f \boldsymbol{n}_f &=& \boldsymbol{0}, & \text{on } \Gamma_{out}, & \text{(the free outflow condition)} \\
\boldsymbol{u}_f(\boldsymbol{x},0) &=& \boldsymbol{u}_{f}^0, &
\text{in } \hat{\Omega}_f, & \text{(the initial condition)}
\end{array}
\end{equation}
and the fluid stress tensor is defined as
\begin{equation}
\boldsymbol{\sigma}_f(\boldsymbol{u}_f, p_f) = \rho_f \nu_f \left(\nabla \boldsymbol{u}_f + (\nabla \boldsymbol{u}_f)^T\right) - p_f \boldsymbol{I},
\end{equation}
where $\rho_f$ is the fluid density, $\nu_f$ is the kinematic
viscosity, $\boldsymbol{b}_f$ is the body force per unit volume, and
$\boldsymbol{I}$ is the identity tensor.

In general, the fluid equation (\ref{nseq}) is modeled in Eulerian
description, i.e., the fluid domain $\Omega_f(t)$ is associated with
time-dependent spatial coordinates $\boldsymbol{x}_f$, which can be
mapped from the reference domain $\hat\Omega_f$ that is associated
with the spatial coordinates $\hat\bx_f$ through a bijective mapping
$\hat\bX_f(\hat\bx_f,t)$, i.e.,
$\bx_f=\hat{\boldsymbol{X}}_f(\hat{\boldsymbol{x}}_f, t)$. 
It is well known that the mapping $\hat\bX_f(\hat\bx_f,t)$ may not
be unnecessarily defined for fluid equations unless the fluid is
modeled in Lagrangian description itself, or the fluid couples a
structural motion that is modeled in Lagrangian description, in
which a fluid domain displacement, $\hat\bd_f$, may be introduced to
define the coordinate transformation as follows
\begin{equation}
\boldsymbol{x}_f = \hat{\boldsymbol{X}}_f(\hat{\boldsymbol{x}}_f, t)
= \hat{\boldsymbol{x}}_f + \hat{\boldsymbol{d}}_f,
\end{equation}
where $\hat{\boldsymbol{d}}_f$ can be treated as either the material
displacement or an artificial domain displacement defined in the
reference domain $\hat{\Omega}_f$. Throughout this paper, the
notation $\hat{\cdot}$ indicates quantities associated with the
reference configuration $\hat{\Omega}_f$, or $\hat{\Omega}_s$ as
used below for the structural model description.

In this paper, the linear elasticity constitution is chosen to model
the structural motion in Lagrangian description, which can be
extended to more complex nonlinear constitutive relations.
Particularly, the structural motion is governed by the following
linear elasticity equation with respect to the structural
displacement $\hat{\boldsymbol{d}}_s$: find
$\hat{\boldsymbol{d}}_s(t) \in H^2(0,T;H^2(\hat{\Omega}_s)^d)$ such
that
\begin{equation}\label{soliddisplacement}
\hat{\rho}_s \frac{\partial^2 \hat{\boldsymbol{d}}_s}{\partial t^2}
- \widehat{\operatorname{div}}
\hat{\boldsymbol{P}}(\hat{\boldsymbol{d}}_s) =
\hat{\boldsymbol{b}}_s, \quad \text{in } \hat{\Omega}_s \times
(0,T],
\end{equation}
with the first Piola stress tensor,
$\hat{\boldsymbol{P}}(\hat{\boldsymbol{d}}_s)$, being defined as
\begin{equation}
\hat{\boldsymbol{P}}(\hat{\boldsymbol{d}}_s) = 2\mu_s \hat{\boldsymbol{\varepsilon}}(\hat{\boldsymbol{d}}_s) + \lambda_s \operatorname{tr}(\hat{\boldsymbol{\varepsilon}}(\hat{\boldsymbol{d}}_s)) \boldsymbol{I},
\end{equation}
where $\hat{\rho}_s$ is the structural density,
$\hat{\boldsymbol{b}}_s$ is the body force per unit volume, $\mu_s$
and $\lambda_s$ are the Lam\'e constants, the linearized strain
tensor, $\hat{\boldsymbol{\varepsilon}}(\hat{\boldsymbol{d}}_s)$, is
given by
\begin{equation}
\hat{\boldsymbol{\varepsilon}}(\hat{\boldsymbol{d}}_s) =
\frac{1}{2}\left(\hat{\nabla} \hat{\boldsymbol{d}}_s + (\hat{\nabla}
\hat{\boldsymbol{d}}_s)^T\right),
\end{equation}
with the differential operators $\hat{\nabla}$ and
$\widehat{\operatorname{div}}$ being associated with the Lagrangian
coordinates $\hat{\boldsymbol{x}}_s$ in the reference configuration.

The structural equation (\ref{soliddisplacement}) is subject to the
following boundary condition and initial condition:
\begin{equation}\label{solid-BICs}
\begin{array}{rclll}
\hat{\boldsymbol{d}}_s &=& \boldsymbol{0}, & \text{on } \hat{\Gamma}_{walls} \cap \hat{\Gamma}_{FSI}, & \text{(the fixed support condition)} \\
\hat{\boldsymbol{d}}_s(\hat{\boldsymbol{x}},0) &=&
\bm{0}, 
& \text{in } \hat{\Omega}_s, & \text{(the initial condition)}
\end{array}
\end{equation}
which results that the reference/Lagrangian domain of structure is
just its initial domain, i.e., $\hat\Omega_s=\hat\Omega_s(0)$, and
the current structural domain $\Omega_s(t)$ can be updated through
its material/Lagrangian mapping $\hat{\boldsymbol{X}}_s$ such that
    \begin{equation*}
        \begin{aligned}
            \hat\bX_s:~&\hat\Omega_s\rightarrow \Omega_s(t)
            \\&\hat{\bm{x}}_s\rightarrow
            \boldsymbol{x}_s = \hat{\boldsymbol{X}}_s(\hat{\boldsymbol{x}}_s, t)
= \hat{\boldsymbol{x}}_s + \hat{\boldsymbol{d}}_s.
        \end{aligned}
    \end{equation*}

Therefore, in the scenario of FSI with zero initial structural
displacement, the reference domain of fluid is also identical with
its initial domain, i.e., $\hat\Omega_f=\hat\Omega_f(0)$, which is
separated from $\hat\Omega_s$ by the initial interface
$\hat\Gamma_{FSI}$. Now we impose the no-slip type of interface
conditions across $\hat\Gamma_{FSI}$ below to ensure continuity of
velocity and balance of forces between the fluid and structure
through the interface $\hat\Gamma_{FSI}$, defined by
\begin{equation}\label{interface}
\begin{array}{rclll}
\boldsymbol{u}_f(\boldsymbol{x}(\hat{\boldsymbol{x}},t),t) &=& \frac{\partial\hat{\boldsymbol{d}}_s}{\partial t},
& \text{on } \hat{\Gamma}_{FSI} \times [0,T], \text{ (the kinematic condition)} \\
\hat{J}\boldsymbol{\sigma}_f(\boldsymbol{x}(\hat{\boldsymbol{x}},t),t)
\hat{\boldsymbol{F}}^{-T}\hat{\boldsymbol{n}}_s &=&
\hat{\boldsymbol{P}}\hat{\boldsymbol{n}}_s, & \text{on }
\hat{\Gamma}_{FSI} \times [0,T], \text{ (the dynamic condition)}
\end{array}
\end{equation}
where 
$\hat{\boldsymbol{n}}_s$ is the outward unit normal vector on
$\hat{\Gamma}_{FSI}$ pointing from the structural domain into the
fluid domain, the deformation gradient tensor,
$\hat{\boldsymbol{F}}$, is defined as
\begin{align}
\hat{\boldsymbol{F}} &
= \hat{\nabla}
\hat{\boldsymbol{X}}_s = \boldsymbol{I} + \hat{\nabla}
\hat{\boldsymbol{d}}_s,
\end{align}
with the Jacobian $\hat{J} =
\operatorname{det}(\hat{\boldsymbol{F}})$.

\subsection{Full-order modeling of FSI via monolithic ALE-FEM}
In this subsection, we introduce a monolithic, fully discrete ALE-FEM
for the coupled FSI model, which serves as the full-order model
(FOM) for the reduced-order modeling developed in Section \ref{S3}
\cite{HaoSun2021,Wang;Sun2016}. We first describe the ALE mapping
and the corresponding fluid equations, and then present the time
and finite element discretizations of the coupled system.

To accommodate the moving fluid domain, we introduce an artificial
fluid domain displacement $\hat\bd_f:\hat\Omega_f\rightarrow
\mathbb{R}^d$. This displacement follows the structural motion on
$\hat\Gamma_{FSI}$ and extends it into the fluid domain, thereby
avoiding the mesh distortions that may arise from tracking fluid
material points. Here we adopt a harmonic extension to define
$\hat\bd_f\in H^2(\hat\Omega_f)^d$ such that
\begin{equation}\label{alemap}
\left\{\begin{array}{rcll}
-\widehat{\operatorname{div}}\left(\hat{\nabla} \hat{\boldsymbol{d}}_f\right)&=&\bm 0,& \text { in } \hat{\Omega}_f, \\
\hat{\boldsymbol{d}}_f&=&\hat{\boldsymbol{d}}_s, & \text { on }
\hat{\Gamma}_{FSI},\\
\hat{\boldsymbol{d}}_f&=&\bm 0, & \text { on }
\partial\hat{\Omega}_f\backslash\hat{\Gamma}_{FSI},
\end{array}\right.
\end{equation}
The resulting ALE mapping is
${\hat{\bm X}}_f(\hat\bx_f,t)=\bm x_f=\hat{\bm x}_f+\hat{\bm d}_f$
for $\hat\bx_f\in\hat\Omega_f$. It updates the fluid domain
$\Omega_f(t)$ according to the structural displacement at the
shared interface. The associated domain velocity $\bm v_f$ is
defined by
\begin{equation}\label{mesh_velocity}
\hat\bv_f(\hat\bx_f,t)={\bm v}_{f}({\bm x_f},t)\circ{{\hat{\bm
X}}}_f=\frac{\partial{\hat{\bm X}}_f}{\partial t}
=\frac{\partial{\bm x}_f}{\partial t}=\frac{\partial\hat{\bm
d}_f}{\partial t}(\hat{\bm x}_f,t).
\end{equation}
The material derivative on the moving fluid domain can then be
written as
\begin{equation}
D_t {\bm u}_f = \partial_{t}^{{\hat{\bm X}}_f} {\bm u}_f +
\big(({\bm u}_f - {\bm v}_f) \cdot\nabla\big){\bm u}_f,
\end{equation}
where $\partial_t^{{\hat{\bm X}}_f} {\bm u}_f=\frac{\partial {\bm
u}_f}{\partial t}+({\bm v}_f\cdot\nabla) {\bm u}_f$ denotes the ALE
temporal derivative. Hence, the fluid model (\ref{nseq}) can be
reformulated in the ALE frame as
\begin{equation}\label{ALE-momentum}
\left\{
\begin{aligned}
\rho_f \big(\partial_{t}^{{\hat{\bm X}}_f} {\bm u}_f + (({\bm u}_f -
{\bm v}_f) \cdot\nabla){\bm u}_f\big) - \operatorname{div} {\bm
\sigma}_f &= {\bm b}_f,& \hbox{~~~in~} \Omega_f(t)\times (0,T],\\
\operatorname{div}{\bm u}_f&=0,&\hbox{~in~}\Omega_f(t)\times (0,T].
\end{aligned}
\right.
\end{equation}

We now discretize the ALE fluid equations (\ref{ALE-momentum}),
the structural equation (\ref{soliddisplacement}), and the harmonic
extension (\ref{alemap}), subject to the boundary and interface
conditions (\ref{fluid-BICs}), (\ref{solid-BICs}), and
(\ref{interface}). For the time discretization, we introduce a
uniform partition
$0=t_0<t_1<\cdots<t_{N_{T}}=T$ with the time-step size $\Delta
t=T/{N_{T}}$. Set $t_n=n\Delta t$,
$\varphi^n=\varphi(\bx(\hat\bx,t_n),t_n)$ for $n=0,\cdots,N_{T}$,
$d_t\varphi^{n}=\frac{\varphi^{n}-\varphi^{n-1}}{\Delta t}$ and
$d_{tt}\varphi^{n}=\frac{\varphi^{n}-2\varphi^{n-1}+\varphi^{n-2}}{\Delta
t^2}$ for $n=1,\cdots,N_{T}$, where $\varphi^{-1}=\varphi^{0}$.

For the spatial discretization, we triangulate $\hat\Omega_f$ and $\hat\Omega_s$ to obtain two
quasi-uniform meshes, $\hat{\mathcal{T}}_{f,h}$ in $\hat\Omega_f$
and $\hat{\mathcal{T}}_{s,h}$ in $\hat\Omega_s\ (0<h<1)$, which are
conforming across $\hat\Gamma_{FSI}$. Then for any $t\in(0,T]$, we
numerically solve the ALE mapping (\ref{alemap}) for $\hat{\bm
d}_{f,h}$ in the finite element space
$\hat\bV^m_{h}:=\big\{\hat\bxi_f\in H^1(\hat\Omega_f)^d : \hat\bxi_f\big|_K\in
P_1(K)^d,\forall K\in\hat{\mathcal{T}}_{f,h}\big\}$, where $P_k$
denotes the $k$-th degree piecewise polynomial space, to attain the
discrete ALE mapping ${\hat\bX}_{f,h}$ that is smooth and
invertible,
representing the moving fluid mesh, i.e., 
\begin{equation}\label{discretemapping_fluid}
{\hat\bX}_{f,h}(\hat{\bm x}_{f,h},t)=\hat{\bm x}_{f}+\hat{\bm
d}_{f,h},
\end{equation}
where $\hat{\bm {d}}_{f,h}\in \hat\bV^m_{D,h}:=\big\{\hat\bxi_f\in
\hat\bV^m_h:\hat\bxi_f=\bm 0 \text{ on
 }\partial\hat\Omega_f\backslash\hat\Gamma_{FSI}
\big\}$. Then, $\mathcal{T}_{f,h}(t)$ becomes the image of
$\hat{\mathcal{T}}_{f,h}$ through the discrete ALE mapping
${\hat\bX}_{f,h}$, i.e.,
$$
\mathcal{T}_{f,h}(t)={\hat\bX}_{f,h}\left(\hat{\mathcal{T}}_{f,h}\right)
=\hat{\mathcal{T}}_{f,h}+\hat{\bm d}_{f,h}.
$$

Accordingly, the semi-discrete ALE temporal derivative is defined
as:
$$\partial_t^{{\hat\bX}_{f,h}} {\bpsi}_{f,h}
=\frac{\partial {\bpsi}_{f,h}}{\partial
t}+\left(\frac{\partial\hat{\bm d}_{f,h}}{\partial
t}\circ{\hat\bX}_{f,h}^{-1}\right)\cdot\nabla {\bpsi}_{f,h},$$ which
is approximated by the following fully discrete form at $t=t_{n}$:
\begin{equation}\label{BDF2ALE}
\partial_t^{{\hat\bX}_{f,h}} {\bpsi}_{f,h}\approx d_{t}^{{\hat\bX}_{f,h}}{\bpsi}_{f,h}^{n}=\frac{{\bpsi}_{f,h}^{n}
-{\bpsi}_{f,h}^{n-1}\circ{\hat\bX}^{n-1}_{f,h}\circ
({\hat\bX}^{n}_{f,h})^{-1}}{\Delta t}.
\end{equation}

The remaining finite element spaces are defined for
$n=1,\cdots,N_{T}$ as follows:
\begin{equation*}
    \begin{array}{rcl}
        \bV^{f,n}_{h}&:=&\big\{\bpsi_f\in H^1(\Omega_f(t_n))^d : \bpsi_f\big|_K\in
P_1(K)^d,\forall K\in{\mathcal{T}}_{f,h}(t_n)\big\},\\
Q^{f,n}_{h}&:=&\big\{q_f\in L^2(\Omega_f(t_n)) : q_f\big|_K\in P_1(K),\forall
K\in{\mathcal{T}}_{f,h}(t_n)\big\},\\
\bV^{f,n}_{0,h}&:=&\{\bpsi_f\in \bV^{f,n}_h:\bpsi_f=\bm 0 \text{ on
}\Gamma_{in}\cup\Gamma_{walls}\},\\
\bV^{f,n}_{D,h}&:=& \{\bpsi_f\in \bV^{f,n}_h:\bpsi_f=\bu_{in} \text{
on }\Gamma_{in},\bpsi_f=\bm 0 \text{ on
}\Gamma_{walls}\},\\
\hat\bV^s_{h}&:=&\big\{\hat\bpsi_s\in H^1(\hat\Omega_s)^d :
\hat\bpsi_s\big|_K\in P_1(K)^d,\forall
K\in\hat{\mathcal{T}}_{s,h}\big\},\\
\hat\bV^s_{0,h}&:=&\{\hat \bpsi_s \in \hat\bV^s_h:\hat\bpsi_s=\bm 0
\text{ on }\hat\Gamma_{walls}\cap
\hat\Gamma_{FSI}\},\\
\hat\bV^m_{0,h}&:=&\{\hat\bxi_f\in \hat\bV^m_h:\hat\bxi_f=\bm 0
\text{ on }\partial\hat\Omega_f\},
\end{array}
\end{equation*}
Thus, we employ the lowest equal-order mixed finite element,
$P_1/P_1$, for the fluid velocity and pressure, together with a
pressure stabilization term
\cite{Hughes.T;Franca.L;Balestra.M1986,Tezduyar.T1992}, and
piecewise linear elements for the structural and fluid domain
displacements. Combining these approximations with the interface
conditions yields the following fully discrete, monolithic
ALE-FEM as our FOM: find
$({\bm{u}^{n}_{f,h}},p^{n}_{f,h}$,
$\hat{\bm{d}}^{n}_{s,h},\hat{\bm{d}}^{n}_{f,h}) \in
\bm{V}^{f,{n}}_{D,h} \times Q^{f,{n}}_h \times \hat\bV^{s}_{0,h}
\times{\hat\bV^{m}_{D,h}}$, which satisfy
$d_t\hat{\bm{d}}^{n}_{s,h}={\bm{u}^{n}_{f,h}}\circ{\hat\bX}^n_{f,h}$
on $\hat\Gamma_{FSI}$ and
$\hat{\bm{d}}^{n}_{f,h}=\hat{\bm{d}}^{n}_{s,h}$ on
$\hat\Gamma_{FSI}$, such that
\begin{eqnarray}
&&(\hat{\rho}_{s} d_{tt} \hat{\bm d}^{n}_{s,h}, \hat\bpsi_{s}
)_{\hat{\Omega}_{s}} + ({\bm P}(\hat{\bm d}^{n}_{s,h}), \nabla
\hat\bpsi_{s})_{\hat{\Omega}_{s}}+ ({\rho}_{f}
d_{t}^{{\hat\bX}_{f,h}} {\bm u}^{n}_{f,h}, {\bpsi_{f}}
)_{{\Omega}_{f}^n}\notag\\
&& + \left(\big({\bm u}^{n}_{f,h} - d_t\hat{\bm
d}^{n}_{f,h}\circ({\hat\bX}^n_{f,h})^{-1}\big)\cdot\nabla {\bm
u}_{f,h}^{n}, {\bpsi_{f}} \right)_{{\Omega}_{f}^n} +({\bm
\sigma}_f({\bm u}^{n}_{f,h},p^{n}_{f,h}),
\nabla{\bpsi_{f}})_{{\Omega}_{f}^n}
\notag\\
&&+ (\nabla \cdot {\bm u}^{n}_{f,h}, q_{f})_{{\Omega}_{f}^n}
+(\nabla\hat{\bm d}^n_{f,h},\nabla\hat\bxi_f)_{\hat{\Omega}_{f}} +
\delta\frac{h^2}{\rho_f\nu_f}(\nabla p^{n}_{f,h},\nabla
q_{f})_{{\Omega}_{f}^n}\notag\\
&&
= (\hat{\bm b}_s^{n},\hat\bpsi_{s} )_{\hat{\Omega}_{s}}+(\bm
b_f^{n},\bpsi_{f} )_{{\Omega}_{f}^n},\ \forall
\left(\bpsi_{f},q_{f},\hat{\bpsi}_{s},\hat\bxi_f \right) \in
\bV^{f,n}_{0,h}\times Q^{f,n}_{h}\times
\hat\bV^s_{0,h}\times\hat\bV^m_{0,h}\notag\\
&&\text{ satisfying }\hat{\bpsi}_{s}=\bpsi_{f}\circ{\hat\bX}^n_{f,h}
\text{ on } \hat\Gamma_{FSI},\quad
n=1,\cdots,N_{T},\label{disALE-FEM}
\end{eqnarray}
where $\Omega_f^n=\Omega_f(t_n)$ and the last term on the
left-hand side of (\ref{disALE-FEM}) is the pressure stabilization
term with a suitably chosen parameter $\delta$. The system is
nonlinear because of fluid convection and the fluid mesh update
through (\ref{discretemapping_fluid}). To solve it, we adopt a
fixed-point iteration to update the fluid mesh through the discrete
ALE mapping ${\hat\bX}_{f,h}$. At each fixed-point iteration,
Newton's method is used to linearize the fluid convection term,
and the resulting linear systems are solved iteratively until
convergence. Further details of the solution algorithm can be found
in \cite{Wang;Sun2016,HaoSun2021}.

\section{A time-segmented non-intrusive reduced order modeling}\label{S3}
In this section, based on the FOM (\ref{disALE-FEM}), we attempt to
develop a proper non-intrusive ROM to obtain fast FSI simulation
during the online phase. In the offline phase, we first adopt a
greedy parameter-time selection strategy and Gram-Schmidt
orthogonalization to construct the time segmented reduced order
subspaces, then we employ shallow neural networks trained using the
orthogonal greedy algorithm to establish the relationship between
physical parameters and ROM basis coefficients. During the online
phase, we evaluate the trained coefficient maps and reconstruct the
solution for each variable at the requested times, without modifying
the FOM source code. Each prediction is computed directly from its
parameter--time input, so predictions for different variables and
time points can be obtained independently and in parallel. No
sequential time integration of reduced equations or iterative
solution of nonlinear FSI systems is required online. We discuss
both phases and give the theoretical analysis below.

\subsection{Offline Phase}\label{sec:offline}
In order to suppress the growth of approximation errors and achieve
the goal of long-term FSI simulation, we employ an advanced
non-intrusive surrogate modeling approach that combines time
segmentation, greedy parameter selection, and orthogonal greedy
training. The offline phase consists of three main components: time
segmentation, greedy parameter selection and Gram-Schmidt basis
construction, and neural network training.

For the benchmark problem in this paper, the physical parameter
vector $\Theta \in \bm{\Theta}$ consists of the maximum incoming
fluid velocity $\hat{U}$ and the structural shear modulus $\mu_s$,
i.e., $\Theta = (\hat{U}, \mu_s)$.

We use the index $\xi\in\{f,p,s\}$ to distinguish the three
physical quantities: fluid velocity ($f$), fluid pressure ($p$),
and structural displacement ($s$). Let $n_f$, $n_p$, and $n_s$
denote the numbers of spatial degrees of freedom in their respective
finite element representations in the FOM. The combined dimension
of these three variable blocks is $n_f+n_p+n_s$.

For each $\xi$, let
$\bm U_\xi(\Theta,t)\in\mathbb R^{n_\xi}$ denote the corresponding
FOM snapshot vector at parameter $\Theta$ and discrete time $t$,
and let $u_\xi(\cdot;\Theta,t)$ denote its associated finite element
function. Thus, uppercase $\bm U_\xi$ represents a vector of degrees
of freedom, while lowercase $u_\xi$ represents a function of the
spatial coordinate. The offline algorithms operate on these vectors;
the online phase reconstructs a predicted vector and then recovers
the corresponding function through the finite element basis, as
described in Section~\ref{sec:online}.

For greedy selection, Gram--Schmidt orthogonalization, and
coefficient extraction, we use the standard Euclidean inner product
$\langle\cdot,\cdot\rangle_{\ell^2_\xi}$ and its induced
$\ell^2$ norm $\|\cdot\|_{\ell^2_\xi}$ on $\mathbb R^{n_\xi}$.
The subscript $\xi$ identifies the variable block to which the
vectors belong. The convergence analysis in
Section~\ref{sec:convergence} uses spatial $L^2$ norms of the
associated finite element functions.

For the first step of offline, we divide the total time interval
$[0,T]$ into $G$ time segments $[T_{g},T_{g+1}]$, $g =
0,\ldots,G-1$, where $T_0=0$, $T_{G}=T$,
$[0,T]=\bigcup\limits_{g=0}^{G-1}[T_{g},T_{g+1}]$, and $T_g\ (g =
0,\ldots,G)$ coincide with the discrete time points $t_n$ at some
time steps $n\ (0\leq n \leq N_T)$. All discrete time points $t_0,
t_1, \ldots, t_{N_T}$ are reallocated into each $[T_{g},T_{g+1}]\ (g
= 0,\ldots,G-1)$. Each time segment $[T_{g},T_{g+1}]$ is divided
into $N_{T_{g}}-1$ subintervals with the time step size $\Delta t$,
where $N_{T_{g}}=\frac{T_{g+1}-T_g}{\Delta t}+1$, $\sum\limits_{g =
0}^{G-1}(N_{T_{g}}-1) = N_{T}$, and $t_n=T_g+(n-1)\Delta t\in
[T_{g},T_{g+1}]$, $n=1,\cdots,N_{T_g}$.

In the offline phase, as the construction of NIROM for each variable
and each time segment are implemented independently and in parallel,
we just use $\xi$ and $g$ to denote variable and time segment in the
following description.

\subsubsection{Greedy Parameter Selection}\label{sec:greedypra}
To construct reduced spaces from representative FOM snapshots, we
apply greedy parameter--time selection independently for each time
segment $g$ and each variable $\xi$. We first sample
$N_{\Theta,t}$ candidate parameter--time pairs
$(\Theta_i,t_i)\in\bm\Theta\times[T_g,T_{g+1}]$,
$i=1,\ldots,N_{\Theta,t}$, and obtain the corresponding FOM
snapshot vectors $\bm U_\xi(\Theta_i,t_i)$.

\begin{remark}[Parameter--time sample counts]\label{rem:sample_counts}
\normalfont Here $N_{\Theta,t}$ denotes the number of candidate
parameter--time samples for one variable in one fixed time segment;
its dependence on the segment is suppressed in the local algorithm.
For a Cartesian sampling design with $N_\Theta$ physical parameter
combinations and $N_{t,g}$ sampled times in segment $g$, we have
$N_{\Theta,t}=N_\Theta N_{t,g}$. This sample count is distinct from
the physical parameter dimension, the spatial degrees of freedom
$n_\xi$, and the reduced basis dimension $N_{\xi,g}$. In the
numerical experiments, the two velocity components, the two
displacement components, and pressure are processed separately: each
scalar component uses the same candidate input pairs, but its
snapshot values and greedily selected pairs may differ.
\end{remark}

For each fixed variable $\xi$ and time segment $g$, let
$\mathcal S_{g,\xi}$ denote the set of selected parameter--time
pairs, and let $\mathcal R_g$ denote the set of candidate pairs
that have not yet been selected. We initialize
$\mathcal S_{g,\xi}=\emptyset$ and
$\mathcal R_g=\{(\Theta_i,t_i)\}_{i=1}^{N_{\Theta,t}}$.
As pairs are selected, they are added to $\mathcal S_{g,\xi}$ and
removed from $\mathcal R_g$. The remaining set is maintained
separately for each variable; its dependence on $\xi$ is suppressed
in the notation $\mathcal R_g$.

We then choose an initial pair $(\Theta_1^*,t_1^*)\in\mathcal R_g$
whose snapshot vector is nonzero, and update $\mathcal S_{g,\xi}=\{(\Theta_1^*,t_1^*)\}$ and
$\mathcal R_g=\mathcal R_g\setminus\{(\Theta_1^*,t_1^*)\}$.
The corresponding snapshot vector is
$\bm U_{\xi,g,1}=\bm U_\xi(\Theta_1^*,t_1^*)$, and the initial
set of selected snapshot vectors is
$S_{\xi,g,1}=\{\bm U_{\xi,g,1}\}$.
Normalizing this vector gives
$\phi_{\xi,g,1}=\bm U_{\xi,g,1}/\|\bm U_{\xi,g,1}\|_{\ell^2_\xi}$
and the initial orthonormal basis
$B_{\xi,g,1}=\{\phi_{\xi,g,1}\}$.

At iteration $j=2,\ldots,N_{\Theta,t}$, the snapshot set is
$S_{\xi,g,j-1}=\{\bm U_{\xi,g,k}\}_{k=1}^{j-1}$, where
$\bm U_{\xi,g,k}=\bm U_\xi(\Theta_k^*,t_k^*)$. Gram--Schmidt
orthogonalization of these vectors yields the orthonormal basis
$B_{\xi,g,j-1}=\{\phi_{\xi,g,k}\}_{k=1}^{j-1}$.
For each remaining pair $(\Theta,t)\in\mathcal R_g$, we compute
the projection error:

\begin{equation}
e_{\xi,g}(\Theta, t) = \|\bm{U}_\xi(\Theta, t) - \bm{U}_{\xi,\text{proj}}(\Theta, t)\|_{\ell^2_\xi},
\end{equation}
where the projected solution is given by:
\begin{equation}
\bm{U}_{\xi,\text{proj}}(\Theta, t) = \sum_{k=1}^{j-1} \langle \bm{U}_\xi(\Theta, t), \phi_{\xi,g,k} \rangle_{\ell^2_\xi} \phi_{\xi,g,k}.
\end{equation}

We select the parameter-time pair with the largest projection error: $(\Theta_j^*, t_j^*) = \arg\max_{(\Theta, t) \in \mathcal{R}_g} e_{\xi,g}(\Theta, t)$.
If $e_{\xi,g}(\Theta_j^*, t_j^*)$ is small enough, the greedy selection stops here and the final orthonormal reduced-order bases are just  $B_{\xi,g,j-1}$.
Otherwise the parameter-time sets are updated as $\mathcal{S}_{g,\xi} = \mathcal{S}_{g,\xi} \cup \{(\Theta_j^*, t_j^*)\}$ and $\mathcal{R}_g = \mathcal{R}_g \setminus \{(\Theta_j^*, t_j^*)\}$,
and the selected snapshot set is updated as
$S_{\xi,g,j}=S_{\xi,g,j-1}\cup\{\bm U_{\xi,g,j}\}$, where
$\bm U_{\xi,g,j}=\bm U_\xi(\Theta_j^*,t_j^*)$.
Gram--Schmidt orthogonalization of this new vector against the
previous basis yields $\phi_{\xi,g,j}$, and the orthonormal basis
is augmented as $B_{\xi,g,j}=B_{\xi,g,j-1}\cup\{\phi_{\xi,g,j}\}$.

The complete greedy parameter-time selection algorithm is given in Algorithm 1.
\begin{breakablealgorithm}
\caption{Greedy Parameter Selection for variable $\xi$ on $[T_g,T_{g+1}]$}
\begin{algorithmic}
\State \textbf{Input:} $N_{\Theta,t}$ parameter-time pairs $(\Theta_i, t_i)$ in $\bm{\Theta} \times [T_g, T_{g+1}]$, corresponding full model solutions  $\bm{U}_\xi(\Theta_i, t_i), i=1,2,\dots, N_{\Theta,t}$, tolerence $\varepsilon_{\xi,g}$
\State \textbf{Output:} Orthonormal basis $B_{\xi,g}=\{\phi_{\xi,g,k}\}_{k=1}^{N_{\xi,g}}$

 \State Initialize $\mathcal{S}_{g,\xi} = \emptyset$ and remaining candidate set $\mathcal{R}_g = \{(\Theta_i, t_i)\}_{i=1}^{N_{\Theta,t}}$
        \State Choose an initial pair $(\Theta_1^*,t_1^*)\in\mathcal R_g$ with $\bm U_\xi(\Theta_1^*,t_1^*)\ne\bm 0$
        \State Set $\mathcal{S}_{g,\xi} = \{(\Theta_1^*, t_1^*)\}$, $\mathcal{R}_g = \mathcal{R}_g \setminus \{(\Theta_1^*, t_1^*)\}$
        \State Obtain first snapshot: $\bm{U}_{\xi,g,1} = \bm{U}_\xi(\Theta_1^*, t_1^*)$, and set $S_{\xi,g,1}=\{\bm U_{\xi,g,1}\}$
        \State Normalize $\phi_{\xi,g,1}=\bm{U}_{\xi,g,1}/\|\bm{U}_{\xi,g,1}\|_{\ell^2_\xi}$, and set $B_{\xi,g,1} = \{\phi_{\xi,g,1}\}$ and $N_{\xi,g}=1$

    \For{$j = 2$ to $N_{\Theta,t}$}
            \For{each $(\Theta, t) \in \mathcal{R}_g$}
                \State Project onto current bases:
                    \[ \bm{U}_{\xi,\text{proj}}(\Theta, t) = \sum_{k=1}^{j-1} \langle \bm{U}_\xi(\Theta, t), \phi_{\xi,g,k} \rangle_{\ell^2_\xi} \phi_{\xi,g,k} \]
                \State Compute projection error:
                    \[ e_{\xi,g}(\Theta, t) = \|\bm{U}_\xi(\Theta, t) - \bm{U}_{\xi,\text{proj}}(\Theta, t)\|_{\ell^2_\xi} \]
            \EndFor
            \State Select $(\Theta_j^*, t_j^*) = \arg\max_{(\Theta, t) \in \mathcal{R}_g} e_{\xi,g}(\Theta, t)$
            \If {$e_{\xi,g}(\Theta_j^*, t_j^*) < \varepsilon_{\xi,g}$}
                \State stop and return $B_{\xi,g}=\{\phi_{\xi,g,k}\}_{k=1}^{N_{\xi,g}}$
            \EndIf
            \State Update $\mathcal{S}_{g,\xi} = \mathcal{S}_{g,\xi} \cup \{(\Theta_j^*, t_j^*)\}$ and $\mathcal{R}_g = \mathcal{R}_g \setminus \{(\Theta_j^*, t_j^*)\}$
            \State Set $\bm U_{\xi,g,j}=\bm U_\xi(\Theta_j^*,t_j^*)$ and update $S_{\xi,g,j}=S_{\xi,g,j-1}\cup\{\bm U_{\xi,g,j}\}$
            \State Orthogonalize $\bm U_{\xi,g,j}$ against $B_{\xi,g,j-1}$ and normalize to obtain $\phi_{\xi,g,j}$
            \State $B_{\xi,g,j}=B_{\xi,g,j-1}\cup\{\phi_{\xi,g,j}\}$, $N_{\xi,g}=N_{\xi,g}+1$
        \EndFor

\State \textbf{Return:} Orthonormal basis $B_{\xi,g} = \{\phi_{\xi,g,k}\}_{k=1}^{N_{\xi,g}}$
\end{algorithmic}
\end{breakablealgorithm}

Suppose we have totally chosen $N_{\xi,g}$ ($N_{\xi,g}\leq N_{\Theta,t}$) parameter-time pairs, which capture the most important features of the solution space, then the resulting orthonormal basis vectors $\{\phi_{\xi,g,k}\}_{k=1}^{N_{\xi,g}}$ provide a stable and efficient representation of the solution manifold for variable $\xi$ on $[T_g,T_{g+1}]$. 

\subsubsection{Orthogonal Greedy Training Algorithm}\label{sec:oga}

For a fixed variable $\xi$ on $[T_g,T_{g+1}]$,
the preceding greedy selection provides the reduced basis
$\{\phi_{\xi,g,k}\}_{k=1}^{N_{\xi,g}}$. The next task is to learn
the map from the physical parameters and time to the coefficients
of the FOM snapshot vector in this basis. Since
$\Theta=(\hat U,\mu_s)$, the network input is the three-dimensional
vector $\bm x=(\Theta,t)=(\hat U,\mu_s,t)^T$, with input domain
$\mathcal X_g=\bm\Theta\times[T_g,T_{g+1}]\subset\mathbb R^3$.

To define the coefficient map throughout $\mathcal X_g$, we assume
that the FOM snapshot vectors admit a continuous extension in
parameter and time that agrees with the computed snapshots at the
sampled inputs. We define the continuous reduced-coefficient map
\begin{equation}\label{eq:reduced_coefficient_map}
\begin{aligned}
\bm\alpha_{\xi,g}:\mathcal X_g&\longrightarrow\mathbb R^{N_{\xi,g}},\\
\bm x&\longmapsto
\bigl(\alpha_{\xi,g,1}(\bm x),\ldots,
\alpha_{\xi,g,N_{\xi,g}}(\bm x)\bigr)^T,
\end{aligned}
\end{equation}
where $\bm U_\xi(\bm x)$ denotes $\bm U_\xi(\Theta,t)$ and
\begin{equation}\label{eq:reduced_coefficient_projection}
\alpha_{\xi,g,k}(\bm x)
=\langle\bm U_\xi(\bm x),\phi_{\xi,g,k}\rangle_{\ell^2_\xi},
\qquad k=1,\ldots,N_{\xi,g}.
\end{equation}
Thus, the reduced approximation of the FOM snapshot vector is
\begin{equation}\label{eq:greedy_snapshot_projection}
\bm U_\xi(\Theta,t)\approx
\bm U_{\xi,\mathrm{proj}}(\Theta,t)
=\sum_{k=1}^{N_{\xi,g}}\alpha_{\xi,g,k}(\bm x)\phi_{\xi,g,k}.
\end{equation}

For each sampled input $\bm x_j=(\Theta_j,t_j)$,
$j=1,\ldots,N_{\Theta,t}$, we compute the target coefficient
vector $\bm\alpha_{\xi,g}(\bm x_j)$ using
(\ref{eq:reduced_coefficient_projection}). The coefficient-training
sample set for variable $\xi$ on $[T_g,T_{g+1}]$ is
\begin{equation}\label{eq:coefficient_training_set}
\mathcal D_{\xi,g}
=\{(\bm x_j,\bm\alpha_{\xi,g}(\bm x_j))\}_{j=1}^{N_{\Theta,t}}
\subset\mathcal X_g\times\mathbb R^{N_{\xi,g}}.
\end{equation}
This training set uses all $N_{\Theta,t}$ candidate
input pairs introduced in Section~\ref{sec:greedypra}. Each input
has a target vector with $N_{\xi,g}$ components: $N_{\Theta,t}$ is
the number of training samples, and $N_{\xi,g}$ is the output
dimension.

We approximate $\bm\alpha_{\xi,g}$ by a vector-valued shallow
neural network trained with OGA. Choose a fixed positive integer
$\kappa\in\{1,2,\ldots\}$ and use the rectified-power activation
$\sigma_\kappa(z)=(z_+)^\kappa$, where $z_+=\max\{0,z\}$ and
$\kappa$ is its piecewise polynomial degree. Let
$c_1,c_2\in\mathbb R$ be fixed bias bounds with $c_1<c_2$.
The scalar ridge-function dictionary on $\mathcal X_g$ is
\begin{equation}\label{eq:oga_dictionary}
\mathbb D_g=\left\{
\bm x\mapsto\sigma_\kappa(\bm w^T\bm x+b):
\bm w\in S^2,\ b\in[c_1,c_2]\right\},
\end{equation}
where $S^2=\{\bm w\in\mathbb R^3:\|\bm w\|_2=1\}$ is the unit
sphere. Each selected hidden neuron is represented by a scalar
ridge function $\gamma_{i,\xi,g}\in\mathbb D_g$, defined by
\begin{equation}
\gamma_{i,\xi,g}(\bm x)
=\sigma_\kappa(\bm w_{i,\xi,g}^T\bm x+b_{i,\xi,g}),
\end{equation}
with $\bm w_{i,\xi,g}\in S^2$ and
$b_{i,\xi,g}\in[c_1,c_2]$. OGA selects these inner weights and
biases successively and recomputes the output weights after each
selection.

Let $n$ denote the number of completed training steps, each of
which adds one hidden neuron. The resulting network approximates
$\bm\alpha_{\xi,g}$ through
\begin{equation}\label{eq:oga_coefficient_network}
\begin{aligned}
\bm\alpha_{\xi,g}(\bm x)&\approx\bm f_{n,\xi,g}(\bm x),\\
\bm f_{n,\xi,g}(\bm x)
&=\sum_{i=1}^{n}\bm a_{i,\xi,g}^{(n)}\gamma_{i,\xi,g}(\bm x),
\qquad \bm a_{i,\xi,g}^{(n)}\in\mathbb R^{N_{\xi,g}}.
\end{aligned}
\end{equation}
The superscript $(n)$ indicates that all output weights are updated
at training step $n$. We measure the training error by the empirical
squared loss
\begin{equation}\label{eq:oga_empirical_loss}
\mathcal R_{\xi,g}(\bm f)
=\frac{1}{2N_{\Theta,t}}\sum_{j=1}^{N_{\Theta,t}}
\|\bm f(\bm x_j)-\bm\alpha_{\xi,g}(\bm x_j)\|_2^2.
\end{equation}

We initialize the training with $n=0$, the zero network
$\bm f_{0,\xi,g}(\bm x)=\bm 0$, and the scalar approximation
space $V_{0,\xi,g}=\{0\}$. Given $\bm f_{n,\xi,g}$, the residual
map at training step $n$ is
\begin{equation}
\bm r_{n,\xi,g}(\bm x)
=\bm\alpha_{\xi,g}(\bm x)-\bm f_{n,\xi,g}(\bm x).
\end{equation}
The next neuron is selected by maximizing its empirical correlation
with this vector residual:
\begin{equation}\label{eq:oga_neuron_selection}
\gamma_{n+1,\xi,g}\in
\arg\max_{\gamma\in\mathbb D_g}
\left\|\frac{1}{N_{\Theta,t}}\sum_{j=1}^{N_{\Theta,t}}
\bm r_{n,\xi,g}(\bm x_j)\gamma(\bm x_j)\right\|_2.
\end{equation}
Here the correlation is a vector in $\mathbb R^{N_{\xi,g}}$;
its Euclidean norm combines the correlations for all coefficient
components. The selected ridge function determines the new inner
weights $\bm w_{n+1,\xi,g}$ and bias $b_{n+1,\xi,g}$.

After this selection, we enrich the scalar approximation space to
\begin{equation}
V_{n+1,\xi,g}
=\operatorname{span}\{\gamma_{1,\xi,g},\ldots,
\gamma_{n+1,\xi,g}\}.
\end{equation}
All output weights are then recomputed by an empirical least-squares
projection, with each component of the network lying in
$V_{n+1,\xi,g}$. Define the evaluation and target matrices
\begin{equation}
\begin{aligned}
A_{n+1,\xi,g}&\in\mathbb R^{N_{\Theta,t}\times(n+1)},
& (A_{n+1,\xi,g})_{ji}&=\gamma_{i,\xi,g}(\bm x_j),\\
Y_{\xi,g}&\in\mathbb R^{N_{\Theta,t}\times N_{\xi,g}},
& (Y_{\xi,g})_{jk}&=\alpha_{\xi,g,k}(\bm x_j),
\end{aligned}
\end{equation}
where $j=1,\ldots,N_{\Theta,t}$, $i=1,\ldots,n+1$, and
$k=1,\ldots,N_{\xi,g}$. The output-weight matrix is chosen as
\begin{equation}\label{eq:oga_weight_update}
C_{n+1,\xi,g}\in
\arg\min_{C\in\mathbb R^{(n+1)\times N_{\xi,g}}}
\|A_{n+1,\xi,g}C-Y_{\xi,g}\|_F^2,
\end{equation}
where $\|\cdot\|_F$ is the Frobenius norm. Its $i$th row is
$(\bm a_{i,\xi,g}^{(n+1)})^T$, and a least-squares minimizer
satisfies the normal equations
\begin{equation}
A_{n+1,\xi,g}^T A_{n+1,\xi,g}C_{n+1,\xi,g}
=A_{n+1,\xi,g}^T Y_{\xi,g}.
\end{equation}
The update is defined by the least-squares problem
(\ref{eq:oga_weight_update}); it does not assume that
$A_{n+1,\xi,g}^T A_{n+1,\xi,g}$ is invertible. The updated weights
and the selected neurons give $\bm f_{n+1,\xi,g}$ through
(\ref{eq:oga_coefficient_network}).

For a prescribed maximum neuron count $N_{\max}\in\mathbb N$ and loss
tolerance $\varepsilon>0$, training stops when $n=N_{\max}$ or
$\mathcal R_{\xi,g}(\bm f_{n,\xi,g})\leq\varepsilon$. Let $N$ be
the number of completed training steps at termination. The final
network $\bm f_{N,\xi,g}$ is the coefficient-map approximation
used in the online phase. Algorithm~2 summarizes this procedure.

\begin{breakablealgorithm}
\caption{Orthogonal greedy training for variable $\xi$ on $[T_g,T_{g+1}]$}
\begin{algorithmic}
\State \textbf{Input:} Coefficient-training sample set $\mathcal D_{\xi,g}$, dictionary $\mathbb D_g$, tolerance $\varepsilon>0$, maximum neuron count $N_{\max}$
\State \textbf{Output:} Trained coefficient map $\bm f_{N,\xi,g}$
\State Form $Y_{\xi,g}$ with $(Y_{\xi,g})_{jk}=\alpha_{\xi,g,k}(\bm x_j)$
\State Initialize $n\gets0$, $\bm f_{0,\xi,g}\gets\bm0$, and $V_{0,\xi,g}\gets\{0\}$
\While{$n<N_{\max}$ {\bfseries and} $\mathcal R_{\xi,g}(\bm f_{n,\xi,g})>\varepsilon$}
    \State Compute $\bm r_{n,\xi,g}(\bm x_j)\gets\bm\alpha_{\xi,g}(\bm x_j)-\bm f_{n,\xi,g}(\bm x_j)$ for all training inputs
    \State Select $\gamma_{n+1,\xi,g}\in\mathbb D_g$ by (\ref{eq:oga_neuron_selection})
    \State Enrich $V_{n+1,\xi,g}\gets V_{n,\xi,g}+\operatorname{span}\{\gamma_{n+1,\xi,g}\}$
    \State Form $A_{n+1,\xi,g}$ with $(A_{n+1,\xi,g})_{ji}=\gamma_{i,\xi,g}(\bm x_j)$
    \State Obtain $C_{n+1,\xi,g}$ from the least-squares problem (\ref{eq:oga_weight_update})
    \State Set $(\bm a_{i,\xi,g}^{(n+1)})^T$ to the $i$th row of $C_{n+1,\xi,g}$
    \State Update $\displaystyle\bm f_{n+1,\xi,g}(\bm x)\gets\sum_{i=1}^{n+1}\bm a_{i,\xi,g}^{(n+1)}\gamma_{i,\xi,g}(\bm x)$
    \State $n\gets n+1$
\EndWhile
\State $N\gets n$ and return $\bm f_{N,\xi,g}$
\end{algorithmic}
\end{breakablealgorithm}

The offline construction combines the greedy reduced-space
approximation (\ref{eq:greedy_snapshot_projection}) with the
coefficient-map approximation learned by OGA. Define the predicted
coefficients by
$\hat\alpha_{\xi,g,k}(\Theta,t)
=[\bm f_{N,\xi,g}(\Theta,t)]_k$.
Substituting these predictions for the exact projection
coefficients gives the final ROM approximation, for
$(\Theta,t)\in\mathcal X_g$:
\begin{equation}\label{eq:offline_rom_approximation}
\begin{aligned}
\bm U_\xi(\Theta,t)
&\approx\bm U_{\xi,\mathrm{proj}}(\Theta,t)
=\sum_{k=1}^{N_{\xi,g}}\alpha_{\xi,g,k}(\Theta,t)\phi_{\xi,g,k}\\
&\approx\hat{\bm U}_{\xi,g}(\Theta,t)
:=\sum_{k=1}^{N_{\xi,g}}
[\bm f_{N,\xi,g}(\Theta,t)]_k\phi_{\xi,g,k}\\
&=\sum_{k=1}^{N_{\xi,g}}
\left(\sum_{i=1}^{N}(\bm a_{i,\xi,g}^{(N)})_k
\gamma_{i,\xi,g}(\Theta,t)\right)\phi_{\xi,g,k}.
\end{aligned}
\end{equation}
Here the basis vectors $\phi_{\xi,g,k}$ are obtained by greedy
snapshot selection, while their coefficients are predicted by the
trained shallow network. The two approximation steps account for
the reduced-space error and the coefficient-learning error,
respectively.

\subsection{Online Phase}\label{sec:online}
For a new input $(\Theta^*,t^*)$, we select the time segment
$[T_g,T_{g+1}]$ containing $t^*$ and use its stored basis and
trained network. Substituting this input into
(\ref{eq:offline_rom_approximation}) gives
\begin{equation}\label{eq:online_dof_reconstruction}
\hat{\bm U}_{\xi,g}(\Theta^*,t^*)
=\sum_{k=1}^{N_{\xi,g}}
[\bm f_{N,\xi,g}(\Theta^*,t^*)]_k\phi_{\xi,g,k}.
\end{equation}
Thus, online prediction consists of evaluating the network and
forming a linear combination of the stored reduced basis vectors;
it requires no new FOM solve.

\begin{remark}[Spatial reconstruction]\label{rem:fe_reconstruction}
\normalfont FOM and ROM vectors can be converted to spatial functions
using the same finite element basis
$\{\psi_{\xi,i}(\cdot;\Theta,t)\}_{i=1}^{n_\xi}$ at a given
parameter--time pair. This basis is defined on the current ALE
fluid mesh for fluid quantities and on the reference domain for
structural displacement. The linear reconstruction map is
\begin{equation}\label{eq:fe_reconstruction_operator}
\mathcal R_{\xi,h}(\Theta,t)\bm V
=\sum_{i=1}^{n_\xi}V_i\psi_{\xi,i}(\cdot;\Theta,t),
\qquad\bm V\in\mathbb R^{n_\xi}.
\end{equation}
Writing $\bm U_{\xi,g}$ for the FOM vector restricted to
$\mathcal X_g$, the corresponding FOM and ROM functions are
\begin{equation}\label{eq:online_function_reconstruction}
\begin{aligned}
u_{\xi,g}(\cdot;\Theta,t)
&=\sum_{i=1}^{n_\xi}(\bm U_{\xi,g}(\Theta,t))_i
\psi_{\xi,i}(\cdot;\Theta,t),\\
\hat u_{\xi,g}(\cdot;\Theta,t)
&=\sum_{i=1}^{n_\xi}(\hat{\bm U}_{\xi,g}(\Theta,t))_i
\psi_{\xi,i}(\cdot;\Theta,t).
\end{aligned}
\end{equation}
These functions provide a common spatial representation for
computing $L^2$ errors; these errors differ from Euclidean errors
between the coefficient vectors. For the analysis below, we denote
the reconstructed reduced basis functions by
$\varphi_{\xi,g,k}(\cdot;\Theta,t)
:=\mathcal R_{\xi,h}(\Theta,t)\phi_{\xi,g,k}$.
\end{remark}

Compared with online ROM procedures that solve reduced dynamical
equations, the present approach has the following practical
advantages:
\begin{itemize}
    \item \textbf{Low online cost.} Prediction requires only shallow
    network evaluations and reduced basis expansions, avoiding the
    assembly and solution of the coupled FSI equations and their
    nonlinear iterations.
    \item \textbf{Non-intrusive evaluation.} The stored bases and
    coefficient networks provide predictions without modifying the
    FOM source code or accessing its operators during the online phase.
    \item \textbf{Local approximation.} Each physical quantity uses
    its own basis and network on the selected time segment, allowing
    the approximation to reflect the dynamics of that quantity over
    the corresponding interval.
    \item \textbf{Independent and parallel prediction.} For given
    physical parameters, predictions at different times and for
    different variables depend only on their respective inputs.
    They can therefore be evaluated independently and in parallel,
    without sequential time integration of the reduced equations.
\end{itemize}
With parallel computing resources, this independence can reduce
the elapsed time for trajectory evaluations and parameter studies.

\subsection{Convergence Analysis}\label{sec:convergence}

In this section, we provide a comprehensive error analysis for the proposed NIROM with time segmentation, greedy parameter-time selection, and orthogonal greedy training. The total approximation error can be decomposed into two distinct components, each arising from different stages of the algorithm. This decomposition allows for a systematic analysis of the convergence properties and provides theoretical guarantees for the overall method.

The convergence analysis is formulated for the spatial functions reconstructed by~\eqref{eq:online_function_reconstruction}, not for their raw degrees-of-freedom vectors. Spatial arguments are suppressed below for readability.
For $\xi\in\{f,p,s\}$, we use the following $L^2$ norms:
\begin{equation}
\begin{aligned}
\|\bm v_f\|_{V_f}&:=\|\bm v_f\|_{L^2(\Omega_f(t))}
=\left(\int_{\Omega_f(t)}|\bm v_f(x)|^2\,\mathrm dx\right)^{1/2},\\
\|q_f\|_{V_p}&:=\|q_f\|_{L^2(\Omega_f(t))}
=\left(\int_{\Omega_f(t)}|q_f(x)|^2\,\mathrm dx\right)^{1/2},\\
\|\hat{\bm v}_s\|_{V_s}&:=\|\hat{\bm v}_s\|_{L^2(\hat\Omega_s)}
=\left(\int_{\hat\Omega_s}|\hat{\bm v}_s(\hat x)|^2\,\mathrm d\hat x\right)^{1/2}.
\end{aligned}
\end{equation}
Hereafter, $\|\cdot\|_{V_\xi}$ denotes the corresponding continuous $L^2$ norm.

Let $u_{\xi,g}$ denote the finite element function associated with the high-fidelity FSI degrees-of-freedom vector $\bm U_{\xi,g}$, and 
$\hat{u}_{\xi,g}$ be the final approximation obtained from the proposed NIROM. The total error can be expressed as:

\begin{align}
\|u_{\xi,g} - \hat{u}_{\xi,g}\|_{V_\xi} &= \|u_{\xi,g} - u_{\xi,g,\text{greedy}} + u_{\xi,g,\text{greedy}} - \hat{u}_{\xi,g}\|_{V_\xi} \\
&\leq \|u_{\xi,g} - u_{\xi,g,\text{greedy}}\|_{V_\xi} + \|u_{\xi,g,\text{greedy}} - \hat{u}_{\xi,g}\|_{V_\xi},
\end{align}
Here $u_{\xi,g,\text{greedy}}=\mathcal R_{\xi,h}\bm U_{\xi,\mathrm{proj}}$ is the finite element reconstruction of the discrete Euclidean projection used offline. It belongs to the function space spanned by $\{\varphi_{\xi,g,k}\}_{k=1}^{N_{\xi,g}}$; this construction is not, in general, an $L^2$-orthogonal projection. 

The first error component $\|u_{\xi,g} - u_{\xi,g,\text{greedy}}\|_{V_\xi}$ represents the error introduced by the greedy parameter-time selection process. This error quantifies how well the selected orthonormal basis represents the local spatiotemporal solution manifold for variable $\xi$ on $[T_g,T_{g+1}]$. The second component $\|u_{\xi,g,\text{greedy}} - \hat{u}_{\xi,g}\|_{V_\xi}$ is the neural-network approximation error and measures how accurately $\bm f_{N,\xi,g}(\Theta,t)$ approximates the vector coefficient map $(\Theta,t)\mapsto\bm\alpha_{\xi,g}(\Theta,t)$.

\subsubsection{Error of Greedy Parameter-Time Selection}

For the local discrete solution manifold
$\mathcal{M}_{\xi,h}^{g}=\{u_{\xi,g}(\Theta,t):(\Theta,t)\in\bm\Theta\times[T_g,T_{g+1}]\}\subset V_{\xi,h}$, its Kolmogorov $N$-width is defined by
\begin{equation}
d_N(\mathcal{M}_{\xi,h}^{g};V_{\xi,h})
=\inf_{\substack{W_N\subset V_{\xi,h}\\ \dim W_N=N}}
\ \sup_{u\in\mathcal{M}_{\xi,h}^{g}}
\ \inf_{w\in W_N}\|u-w\|_{V_\xi}.
\end{equation}
It is the smallest worst-case error attainable by any $N$-dimensional linear subspace and therefore provides the standard benchmark for assessing the near-optimality of reduced-basis greedy constructions \cite{Pinkus1985Widths,Binev2011GreedyRB,Buffa2012GreedyRB}.  To simplify notation below, we write $\mathcal{M}_{\xi,h}$ for $\mathcal{M}_{\xi,h}^{g}$ when the segment $g$ is fixed.

 Assume that the FSI bilinear forms $a_\xi(\cdot, \cdot; \Theta)$ are continuous, symmetric, and coercive over $\bm{\Theta}$ for each variable $\xi$, with uniform bounds $\gamma_\xi(\Theta) \leq \tilde{\gamma}_\xi$ and $\alpha_\xi(\Theta) \geq \alpha_{\xi,0} > 0$. Then the greedy algorithm's convergence is bounded as follows.

\begin{theorem}[Greedy convergence bound]
Let $V_{N_{\xi,g}}$ be the reduced space constructed via the greedy algorithm for variable $\xi$ on $[T_g,T_{g+1}]$. There exists $C_\xi > 0$ independent of $N_{\xi,g}$ such that
\begin{equation}
    \begin{aligned}
d(\mathcal{M}_{\xi,h}; V_{N_{\xi,g}}) = \sup_{(\Theta, t) \in \bm{\Theta} \times [T_g, T_{g+1}]} \|u_{\xi,g}(\Theta, t) - \Pi_{V_{N_{\xi,g}}} u_{\xi,g}(\Theta, t)\|_{V_\xi} \\
\leq C_\xi(N_{\xi,g}+1) \delta_0^{N_{\xi,g}+1} d_{N_{\xi,g}}(\mathcal{M}_{\xi,h}; V_{\xi,h})
\end{aligned}
\end{equation}
with $\delta_0 = 2$. Consequently, the reduced basis solution error satisfies
\begin{equation}
    \begin{aligned}
\sup_{(\Theta, t) \in \bm{\Theta} \times [T_g, T_{g+1}]} \|u_{\xi,g}(\Theta, t) - u_{\xi,g,\text{greedy}}(\Theta, t)\|_{V_\xi} \\
\leq C_\xi \sqrt{\tilde{\gamma}_\xi/\alpha_{\xi,0}} (N_{\xi,g}+1) \delta_0^{N_{\xi,g}+1} d_{N_{\xi,g}}(\mathcal{M}_{\xi,h}; V_{\xi,h}).
\end{aligned}
\end{equation}
\end{theorem}

Then the error bound depends on the decay of Kolmogorov $N$-width as $d_N(\mathcal{M}_{\xi,h}; V_{\xi,h})$, and we have the following results.

\begin{corollary}[Greedy Exponential convergence]
For variable $\xi$ on $[T_g,T_{g+1}]$, if $d_{N_{\xi,g}}(\mathcal{M}_{\xi,h}; V_{\xi,h}) \leq c_\xi e^{-\delta_\xi N_{\xi,g}}$ with $\delta_\xi > \log \delta_0$, then $\exists \eta_\xi > 0$ such that
\begin{equation}
\|u_{\xi,g}(\Theta, t) - u_{\xi,g,\text{greedy}}(\Theta, t)\|_{V_\xi} \leq C_\xi e^{-\eta_\xi N_{\xi,g}}, \quad \forall (\Theta, t) \in \bm{\Theta} \times [T_g, T_{g+1}].
\end{equation}

\end{corollary}


\begin{corollary}[Greedy Algebraic Convergence]
For variable $\xi$ on $[T_g,T_{g+1}]$, if $d_{N_{\xi,g}}(\mathcal{M}_{\xi,h}; V_{\xi,h}) \leq M_\xi N_{\xi,g}^{-\beta_\xi}$ for $M_\xi, \beta_\xi > 0$, then $\exists C_\xi = C_\xi(\rho, \beta_\xi)$ such that
\begin{equation}
    \begin{aligned}
\sup_{(\Theta, t) \in \bm{\Theta} \times [T_g, T_{g+1}]} \|u_{\xi,g}(\Theta, t) - u_{\xi,g,\text{greedy}}(\Theta, t)\|_{V_\xi}
\leq C_\xi \sqrt{\tilde{\gamma}_\xi/\alpha_{\xi,0}} M_\xi N_{\xi,g}^{-\beta_\xi}.
\end{aligned}
\end{equation}
\end{corollary}



This theoretical framework provides the foundation for understanding the convergence behavior of the greedy parameter-time selection, which ensuring that the selected parameter-time set optimally represents the local spatiotemporal solution manifold for each variable across different time segments.

\subsubsection{Neural Network Approximation with OGA Training}

The vector-valued network $\bm f_{N,\xi,g}$ approximates
$\bm\alpha_{\xi,g}:\bm\Theta\times[T_g,T_{g+1}]\to
\mathbb{R}^{N_{\xi,g}}$.  Accordingly,
\begin{equation}
\begin{aligned}
u_{\xi,g,\mathrm{greedy}}(\Theta,t)
&=\sum_{j=1}^{N_{\xi,g}}[\bm\alpha_{\xi,g}(\Theta,t)]_j\varphi_{\xi,g,j},\\
\hat u_{\xi,g}(\Theta,t)
&=\sum_{j=1}^{N_{\xi,g}}[\bm f_{N,\xi,g}(\Theta,t)]_j\varphi_{\xi,g,j}.
\end{aligned}
\end{equation}

We now define the approximation class on the three-dimensional
parameter--time domain $\mathcal X_g$. With the activation degree
$\kappa$ introduced in Section~\ref{sec:oga}, let
\begin{equation}
\mathbb{P}_\kappa^3
=\{\bm x\mapsto\sigma_\kappa(\bm w^T\bm x+b):\bm w\in S^2,\ b\in[c_1,c_2]\}.
\end{equation}
For a dictionary $\mathbb D$ in a Hilbert space $H$, define its closed absolutely convex hull and the associated variation space by
\begin{equation}
B_1(\mathbb D)
=\overline{\left\{\sum_{i=1}^{M}c_i g_i:
M\in\mathbb N,\ g_i\in\mathbb D,\ \sum_{i=1}^{M}|c_i|\leq1\right\}}^{H},
\end{equation}
\begin{equation}
\mathcal K_1(\mathbb D)=\bigcup_{c>0}cB_1(\mathbb D),
\qquad
\|v\|_{\mathcal K_1(\mathbb D)}
=\inf\{c>0:v\in cB_1(\mathbb D)\}.
\end{equation}
Thus, $\mathcal K_1(\mathbb P_\kappa^3)$ denotes the variation space generated by ReLU$^\kappa$ ridge functions.  For the vector output of the present ROM, we use the natural extension
\begin{equation}
B_1(\mathbb P_\kappa^3;\mathbb R^q)
=\overline{\left\{\sum_{i=1}^{M}c_i\bm v_i g_i:
g_i\in\mathbb P_\kappa^3,\ \|\bm v_i\|_2=1,\ \sum_{i=1}^{M}|c_i|\leq1\right\}}^{L^2(\mathcal X_g;\mathbb R^q)},
\end{equation}
with the corresponding gauge norm
$\|\cdot\|_{\mathcal K_1(\mathbb P_\kappa^3;\mathbb R^q)}$; here $q=N_{\xi,g}$.  These definitions and the entropy-based OGA rates follow the variation-space framework developed for shallow networks in \cite{SiegelXu2022OGA,SiegelHongJinHaoXu2023}.

\begin{theorem}[Neural-network approximation error with OGA]
\label{thm:nn_oga_error}
Set $q=N_{\xi,g}$, assume
\begin{equation}
\bm\alpha_{\xi,g}\in\mathcal K_1(\mathbb P_\kappa^3;\mathbb R^q),
\end{equation}
and that
\begin{equation}
\epsilon_n(B_1(\mathbb P_\kappa^3))_{H^m(\mathcal X_g)}
\lesssim n^{-\rho},
\qquad
\rho=\frac{1}{2}+\frac{2(\kappa-m)+1}{6}.
\end{equation}
Then the coefficient-learning contribution admits the decomposition
\begin{equation}
\mathbb E\!\left[\|u_{\xi,g,\mathrm{greedy}}-\hat u_{\xi,g}\|_{V_\xi}\right]
\leq E_{\mathrm{approx}}+E_{\mathrm{quad}}+E_{\mathrm{opt}},
\end{equation}
where the entropy-controlled OGA approximation term satisfies
\begin{equation}
E_{\mathrm{approx}}
\lesssim
\|\bm\alpha_{\xi,g}\|_{\mathcal K_1(\mathbb P_\kappa^3;\mathbb R^q)}n^{-\rho}.
\end{equation}
The quadrature term is expressed in terms of the number of training samples $N_{\Theta,t}$:
\begin{equation}
\begin{aligned}
E_{\mathrm{quad}}&\lesssim N_{\Theta,t}^{-(\kappa+1)/3}
&&\text{(tensorized Gauss--Legendre)},\\
E_{\mathrm{quad}}&\lesssim N_{\Theta,t}^{-1/2}
&&\text{(Monte Carlo)}.
\end{aligned}
\end{equation}
The remaining $E_{\mathrm{opt}}$ measures the inexact solution of the OGA maximization and projection subproblems; it vanishes for exact OGA and should otherwise be added explicitly.  
\end{theorem}

\begin{corollary}[Sampling requirement]
To keep quadrature from degrading the entropy rate $n^{-\rho}$, it is sufficient to choose
\begin{equation}
\begin{aligned}
N_{\Theta,t}&\gtrsim n^{3\rho/(\kappa+1)}
&&\text{for Gauss--Legendre sampling},\\
N_{\Theta,t}&\gtrsim n^{2\rho}
&&\text{for Monte Carlo sampling}.
\end{aligned}
\end{equation}
Under either choice and exact OGA, the coefficient-learning error is of order $n^{-\rho}$.
\end{corollary}

\begin{theorem}[Balanced greedy--neural-network error]
\label{thm:balanced_error}
Let the preceding greedy and OGA assumptions hold, and take $N_{\Theta,t}$ large enough to satisfy the sampling requirement.  The following estimate uses $N_{\xi,g}$ for the greedy approximation and $n$ for the network approximation:
\begin{equation}
\|u_{\xi,g}-\hat u_{\xi,g}\|_{V_\xi}
\lesssim E_{\mathrm{greedy}}(N_{\xi,g})
+\|\bm\alpha_{\xi,g}\|_{\mathcal K_1(\mathbb P_\kappa^3;\mathbb R^{N_{\xi,g}})}n^{-\rho}
+E_{\mathrm{opt}}.
\end{equation}
The two standard balances are
\begin{equation}
\begin{aligned}
E_{\mathrm{greedy}}(N_{\xi,g})\lesssim e^{-\eta_\xi N_{\xi,g}}
&\quad\Longrightarrow\quad
N_{\xi,g}\asymp(\rho/\eta_\xi)\log n,\\
E_{\mathrm{greedy}}(N_{\xi,g})\lesssim N_{\xi,g}^{-\beta_\xi}
&\quad\Longrightarrow\quad
N_{\xi,g}\asymp n^{\rho/\beta_\xi}.
\end{aligned}
\end{equation}
\end{theorem}

Summing the last estimate over $\xi\in\{f,p,s\}$ and all time segments gives the total FSI error bound.  In this estimate, $N_{\Theta,t}$ specifies the sampling and quadrature resolution, and $N_{\xi,g}$ specifies the dimension and greedy approximation quality of the ROM space.

\section{Numerical Experiments}\label{S4}
In this section, we carry out some numerical experiments by applying
both FOM and the developed NIROM to a FSI benchmark problem
initiated by Turek and Hron \cite{turek2006proposal,
turek2010numerical}, and then conduct a comparison study between them. 

The domain of the studied FSI benchmark problem is illustrated  in
Figure \ref{domainsetting} or sketched in Figure
\ref{Computationaldomain}, whose details are depicted below.

\begin{itemize}
    \item Domain length $L=2.5$ and height $H=0.41$.
    \item The center of obstacle cylinder is positioned at $C=(0.2,0.2)$ with the radius $r=0.05$.
    \item The elastic structural beam's length $l=0.35$ and height $h=0.02$, whose lower
    right corner is positioned at $(0.6,0.19)$, and whose left end is fully attached to the obstacle cylinder.
    \item The control point is $A(t)$ that lies at the midpoint of right end of the beam
    with $A(0)=(0.6,0.2)$, and $B=(0.15,0.2)$.
\end{itemize}
\begin{figure}[H]
    \centering
    \includegraphics[width = 12cm, height = 3cm]{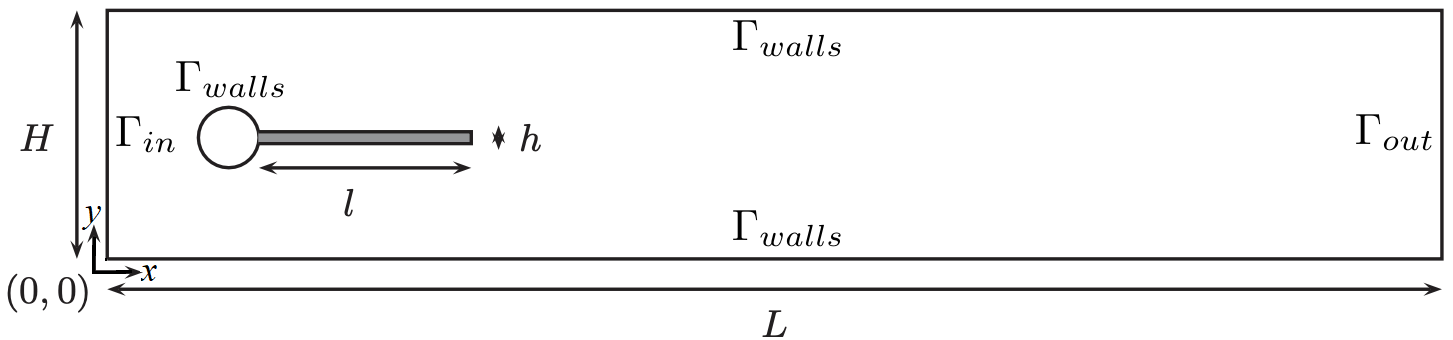}
    \includegraphics[width = 11cm, height = 2.25cm]{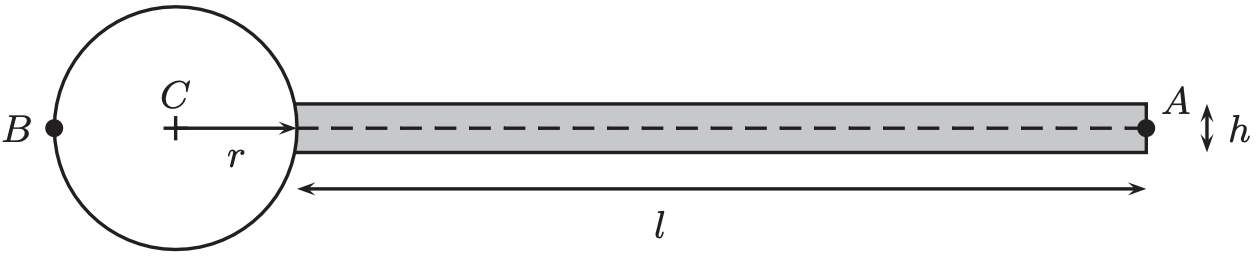}
    \caption{Computational domain of the FSI benchmark problem}
    \label{Computationaldomain}
\end{figure}

We impose the no-slip boundary condition for the fluid velocity on
$\Gamma_{walls}$ that include the surface of obstacle cylinder and
the top and bottom walls of the channel domain, and do-nothing
condition on the outlet $\Gamma_{out}$. In addition, we impose a
parabolic profile of the incoming flow, i.e., the boundary condition
of fluid velocity on the inlet $\Gamma_{in}$ is assigned to
$({u}_{\text {in}},0)^{\top}$ such that
\begin{equation}
{u}_{\text {in}}:=\left\{\begin{array}{ll}
\hat{u}(y) \frac{1-\cos \left(\frac{\pi}{2} t\right)}{2},& \text {if } t<2 s, \\
\hat{u}(y),& \text {otherwise, }
\end{array}\right.
\end{equation}
where $0\leq y\leq 0.41$ is the $y$-coordinate variable, and
\begin{equation}
\hat{u}(y)=\hat{U} \times \frac{4}{0.1681} y(0.41-y),
\end{equation}
here the value of $\hat{U}$ is reported in Table \ref{psforFSI}, in
which all other physical parameters of the FSI problem are reported
as well. We also require the beam to be attached to the obstacle
cylinder, therefore $\hat{\boldsymbol{d}}_s=0$ on $\Gamma_{walls}$.
In what follows, all numerical computations are performed on a
12th-Gen Intel 3.20 GHz Core i9 computer.
\begin{table}[H]
    \centering
    \begin{tabular}{|c|c|c|c|}
\hline Symbol & Description & Value & Unit \\
\hline $\rho_{s}$ & Density of the structure & $10^{4}$ &
$\frac{\mathrm{kg}}{\mathrm{m}^{3}}$\\
\hline $\lambda_{s}$ & Lam\'e constant of the structure & $0.4$ & --- \\
\hline
$\mu_{s}$ & Shear modulus of the structure & $0.4 \times 10^{6} \sim 0.6\times 10^{6}$ & $\frac{\mathrm{kg}}{\mathrm{ms}^{2}}$\\
\hline
$\rho_{f}$ & Density of the fluid & $10^{3}$ & $\frac{\mathrm{kg}}{\mathrm{m}^{3}}$ \\
\hline
$\nu_{f}$ & Kinematic viscosity of the fluid & $10^{-3}$ & $\frac{\mathrm{m}^{2}}{\mathrm{s}}$  \\
\hline
$\hat{U}$ & Parameter of incoming fluid velocity & $1.4 \sim 1.6$ & $\frac{\mathrm{m}}{\mathrm{s}}$ \\
\hline

\end{tabular}
    \caption{Parameter settings for the FSI benchmark problem}
    \label{psforFSI}
\end{table}

\subsection{Full Order Model Results}

\begin{figure}[H]
    \centering
\includegraphics[width = 10cm, height = 1.7cm]{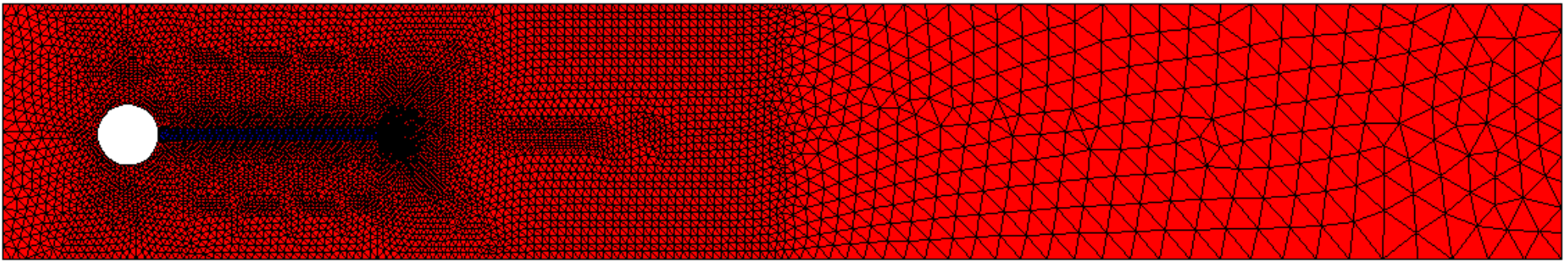}
    \caption{The triangular mesh of benchmark domain.}\label{meshused}
\end{figure}

To verify the effectiveness of our developed ROM by comparing with
FOM, we first perform comprehensive FOM computations for this FSI benchmark
problem. The computations are carried out using the same time step size $\Delta t=0.001s$ as used in
the benchmark problem proposed in
\cite{turek2006proposal,turek2010numerical} on a fine triangular
mesh with $18350$ degree of freedoms (DOFs) in total, as shown in
Figure \ref{meshused}.

To investigate the parametric behavior of the FSI system, we conduct
a systematic parameter study by varying two key physical parameters:

\begin{itemize}
    \item Shear modulus: $\mu_s \in \{0.4 \times 10^6, 0.42 \times 10^6, 0.44 \times 10^6, \ldots, 0.58 \times 10^6, 0.6 \times 10^6\}$
    \item Maximum fluid velocity: $\hat{U} \in \{1.4, 1.42, 1.44, \ldots, 1.58, 1.6\}$
\end{itemize}

This results of $11 \times 11 = 121$ parameter combinations,
providing a comprehensive dataset for training and validating our ROM approach.
The choice of 121 parameter combinations is motivated by the need to capture
the complex nonlinear behavior of the FSI system across the parameter space,
as discussed in \cite{zhai2025projection}, which ensures
adequate representation of the parameter-dependent dynamics while maintaining
computational feasibility for the training phase.
For each parameter combination, we perform FOM simulations with $T=15s$, which is sufficient to capture the
steady-state vibration behavior of the system. The vibration response
at location $A$ (the end of beam tail) is recorded and analyzed.

\begin{figure}[H]
    \centering
    \subfigure[\scriptsize$(\hat{U}, \mu_s) = (1.4, 0.4 \times 10^6)$]{
        \includegraphics[width = 0.3\textwidth]{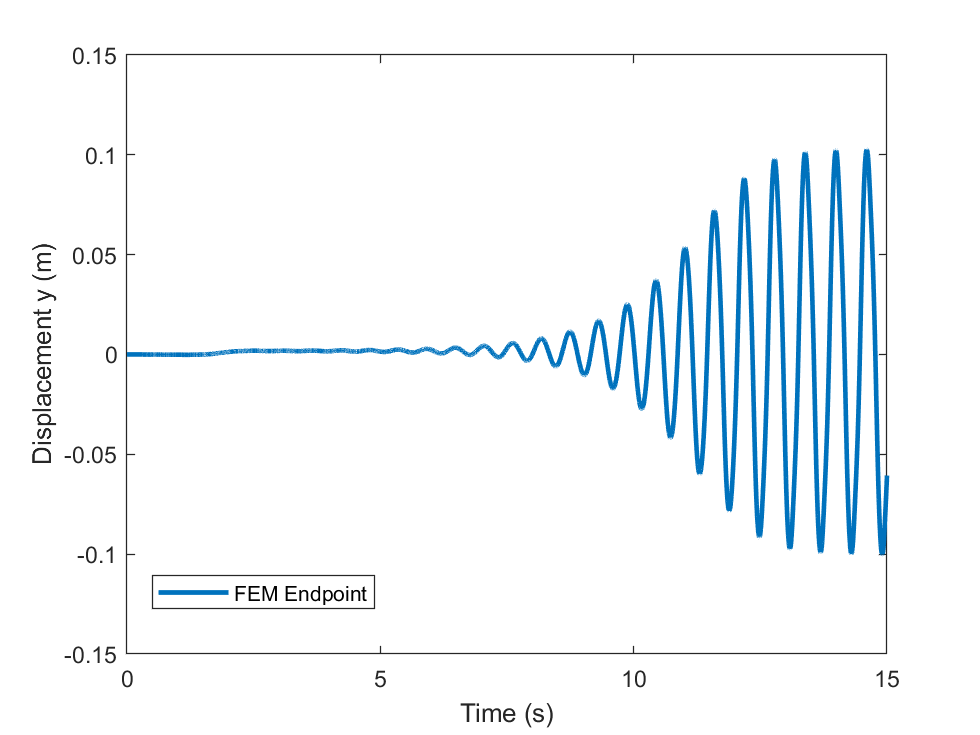}
        \label{fig:fem14040}
    }
    \subfigure[\scriptsize$(\hat{U}, \mu_s) = (1.46, 0.5 \times 10^6)$]{
        \includegraphics[width = 0.3\textwidth]{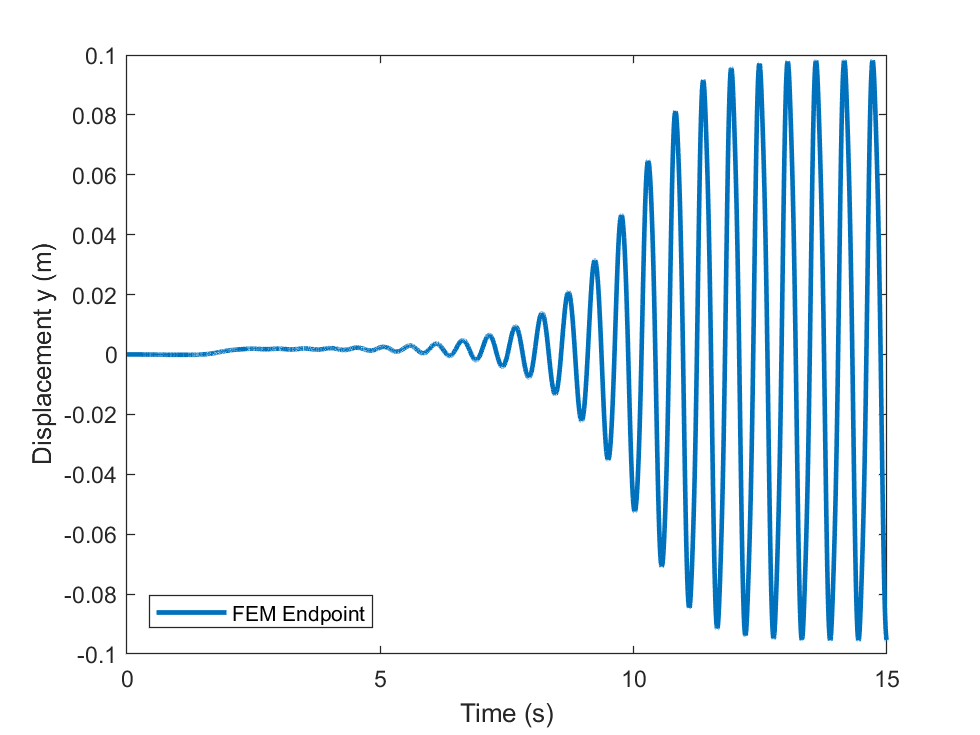}
        \label{fig:fem146050}
    }
    \subfigure[\scriptsize$(\hat{U}, \mu_s) = (1.50, 0.6 \times 10^6)$]{
        \includegraphics[width = 0.3\textwidth]{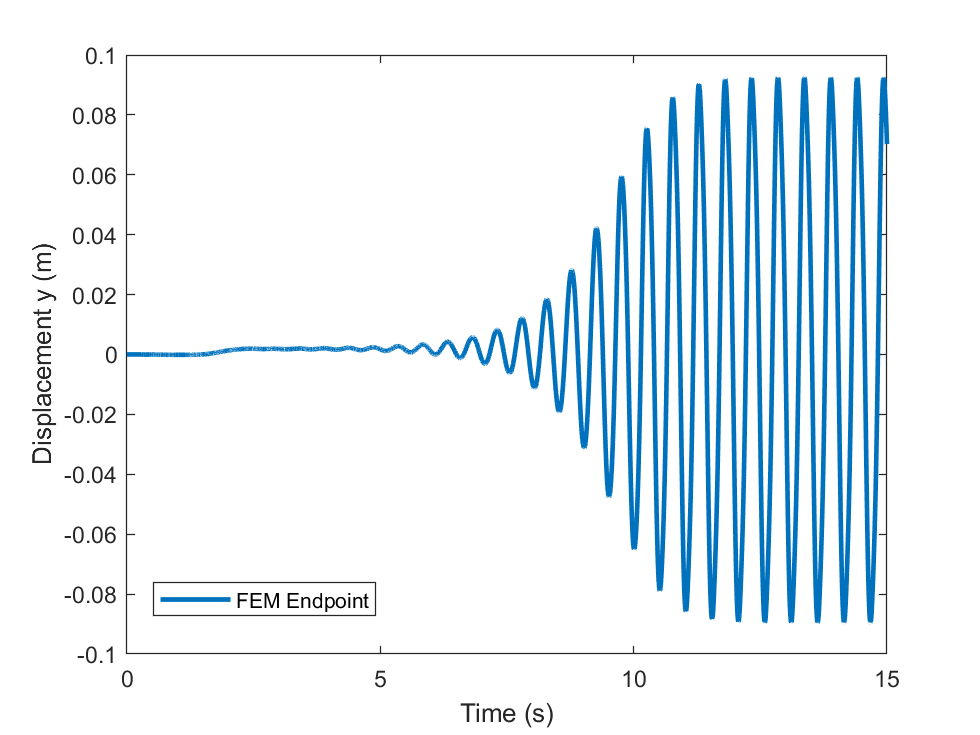}
        \label{fig:fem15060}
    }
    \subfigure[\scriptsize$(\hat{U}, \mu_s) = (1.54, 0.46 \times 10^6)$]{
        \includegraphics[width = 0.3\textwidth]{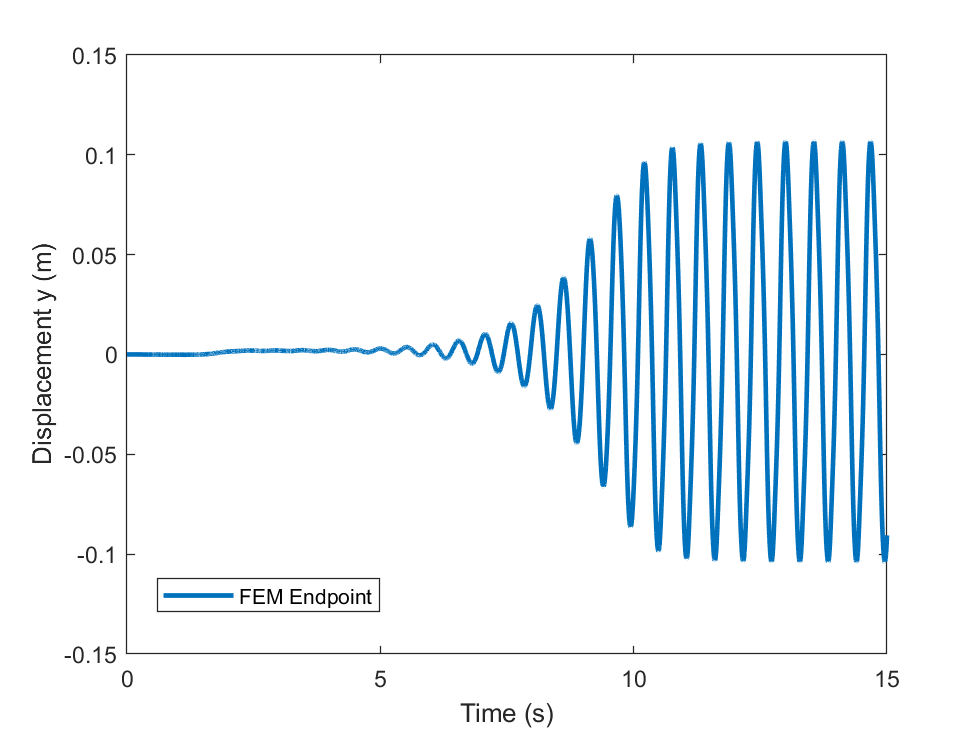}
        \label{fig:fem154046}
    }
    \subfigure[\scriptsize$(\hat{U}, \mu_s) = (1.6, 0.54 \times 10^6)$]{
        \includegraphics[width = 0.3\textwidth]{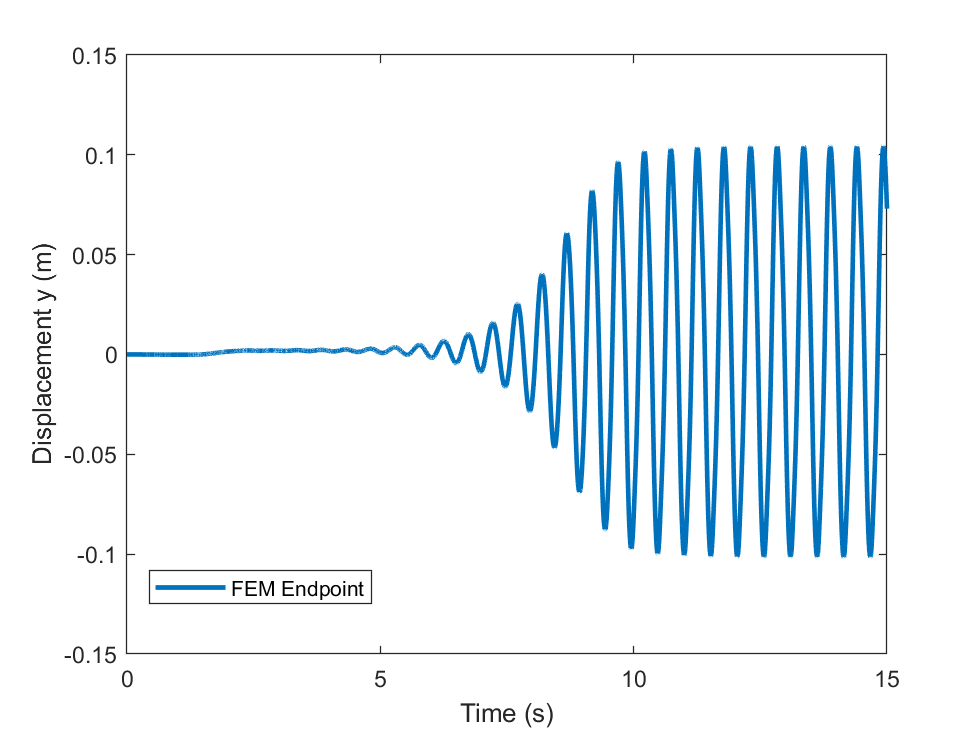}
        \label{fig:fem160054}
    }
    \caption{FEM results of vibration response at location $A$ for selected parameter combinations. Each subfigure shows the displacement response over time for different combinations of $\hat{U}$ and  $\mu_s$.}
    \label{Hfnsdispy18750}
\end{figure}

To demonstrate the parametric behavior of the system, we present
representative results for several key parameter combinations. Figure \ref{Hfnsdispy18750} illustrates the numerical
results of vibration response at location $A$ for selected parameter
combinations: $(\hat{U}, \mu_s) = (1.4, 0.4 \times 10^6)$,
$(1.46, 0.5 \times 10^6)$, $(1.50, 0.6 \times 10^6)$,
$(1.54, 0.46 \times 10^6)$, and $(1.6, 0.54 \times 10^6)$. These results clearly demonstrate the characteristic oscillatory
behavior of the FSI system and reveal the significant influence of
both parameters on the system response. As observed in
Figure \ref{Hfnsdispy18750}, the vibrating frequency and amplitude
exhibit distinct patterns for different parameter combinations,
with the system reaching steady-state behavior after approximately
15 seconds. This observation confirms our choice of simulation
time $T=15s$ for all subsequent tests.

The representative curves in Figure~\ref{Hfnsdispy18750} illustrate the range of responses observed over the full $11\times11$ parameter study.  The displayed combinations show appreciable changes in both oscillation amplitude and frequency, indicating a joint dependence on the inlet-velocity and structural-stiffness parameters.  These results motivate an efficient ROM capable of representing the coupled parameter dependence over the complete training domain.

The computational cost analysis in Table~\ref{fom_computational_cost} shows that one FOM trajectory requires approximately 40 hours.  Generating all 121 trajectories therefore consumes about $4{,}840$ CPU-hours before the subsequent greedy, data-preparation, and OGA stages.  This cost makes direct FOM use impractical for real-time applications or extensive online parameter studies and motivates the ROM developed below.

\begin{table}[H]
    \centering
    \begin{tabular}{|c|c|c|c|}
    \hline
    Parameter Set & $\mu_s$ (Pa) & $\hat{U}$ (m/s) & CPU Time (hours) \\
    \hline
    Baseline & $0.5 \times 10^6$ & 1.5 & 43.1 \\
    \hline
    Low stiffness & $0.4 \times 10^6$ & 1.5 & 41.9 \\
    \hline
    High stiffness & $0.6 \times 10^6$ & 1.5 & 42.2 \\
    \hline
    Low velocity & $0.5 \times 10^6$ & 1.4 & 41.8 \\
    \hline
    High velocity & $0.5 \times 10^6$ & 1.6 & 42.5 \\
    \hline
    \end{tabular}
    \caption{Computational cost for selected FOM simulations.} 
    \label{fom_computational_cost}
\end{table}

\subsection{Reduced Order Model Results}

In this subsection, we present the numerical results obtained using our proposed non-intrusive ROM approach with greedy parameter selection and orthogonal greedy algorithm training. 

\subsubsection{ROM Setup and Parameters}\label{romsetup}

This subsection summarizes all key parameters used in the ROM construction. The detailed implementation and numerical results are presented in the subsequent subsections.

Although each FOM simulation is run over $[0,15]$s, the ROM construction begins immediately after $t=2$s. In this benchmark, the interval $[0,2]$s is an inflow-ramp buffer used by the full-order solver to raise the inlet profile smoothly from rest to its prescribed value and thereby initialize the flow without an abrupt boundary-condition jump. This solver-initialization transient is not part of the sustained FSI dynamics targeted by the reduced model. We therefore retain it in every FOM run but exclude it from snapshot selection and OGA training. Thus, with $\Delta t=0.001$s, the first ROM sample is at $t=2.001$s.

\textbf{Time Segmentation:}
\begin{itemize}
    \item 
    FOM time interval: $[0,15]$s; ROM sampling window: $(2,15]$s, containing 13,000 sampled times. 
    \item Number of time segments: $G=55$
    \item Time segment structure:
    \begin{itemize}
        \item 
        Segment 1: $(2,6]$s, with time-step indices 2001--6000 and $N_{t,1}=4000$.
        \item 
        Segments 2--5: $(6,7]$, $(7,8]$, $(8,9]$, and $(9,10]$s, respectively; each segment contains $N_{t,g}=1000$ sampled times.
        \item 
        Segments 6--55: $(10+0.1(g-6),10+0.1(g-5)]$s for $g=6,\ldots,55$; each segment contains $N_{t,g}=100$ sampled times.
    \end{itemize}
\end{itemize}

\textbf{Greedy Parameter-Time Selection:}
\begin{itemize}
    \item Variables: The greedy selection is applied independently to five variables:
    \begin{itemize}
        \item Horizontal fluid velocity ($u_f$)
        \item Vertical fluid velocity ($v_f$)
        \item Horizontal structure displacement ($u_s$)
        \item Vertical structure displacement ($v_s$)
        \item Fluid pressure ($p_f$)
    \end{itemize}
    \item Sampling parameter set: $11 \times 11 = 121$ parameter combinations, where $\mu_s \in \{0.4 \times 10^6, 0.42 \times 10^6, \ldots, 0.6 \times 10^6\}$ and $\hat{U} \in \{1.4, 1.42, \ldots, 1.6\}$.
    \item 
    Parameter--time training samples per scalar component: the physical parameter vector $\Theta=(\hat U,\mu_s)$ has dimension two, and the network input $(\Theta,t)$ has dimension three. The $N_\Theta=121$ parameter combinations are paired with every sampled time in the relevant segment. Thus, for each component separately, $N_{\Theta,t}=121N_{t,g}$, giving 484,000 samples in segment 1, 121,000 samples in each of segments 2--5, and 12,100 samples in each of segments 6--55; see Table~\ref{tab:training_sample_counts}. These counts are not multiplied by the number of spatial degrees of freedom or by the five scalar components. 
    \item Stopping tolerance: $\varepsilon_{\xi,g}=10^{-6}$. The dimension of the reduced space $N_{\xi,g}$ varies across time segments and variables, reflecting the local solution complexity. 
\end{itemize}

\begin{table}[H]
    \centering
    \begingroup
    \small
    \setlength{\tabcolsep}{4pt}
    \begin{tabular}{|l|r|r|r|}
    \hline
    Scalar component & Segment 1 & Each segment 2--5 & Each segment 6--55 \\
    \hline
    Sampled times per segment & 4,000 & 1,000 & 100 \\
    Parameter combinations & 121 & 121 & 121 \\
    \hline
    Horizontal fluid velocity ($u_f$) & 484,000 & 121,000 & 12,100 \\
    Vertical fluid velocity ($v_f$) & 484,000 & 121,000 & 12,100 \\
    Horizontal displacement ($u_s$) & 484,000 & 121,000 & 12,100 \\
    Vertical displacement ($v_s$) & 484,000 & 121,000 & 12,100 \\
    Fluid pressure ($p_f$) & 484,000 & 121,000 & 12,100 \\
    \hline
    \end{tabular}
    \caption{
    Number of parameter--time training samples for each scalar component in each time segment, $N_{\Theta,t}=N_\Theta N_{t,g}$. Entries for segments 2--5 and 6--55 are per-segment counts, not totals over those groups. The five components share the sampled input pairs but have separate target data.}
    \label{tab:training_sample_counts}
    \endgroup
\end{table}

\textbf{Neural Network Training:}
\begin{itemize}
    \item Network architecture: Shallow neural networks, $N_{max}=5000$ neurons
    \item Activation function: ReLU$^3$ ($\sigma(x) = \max(0, x)^3$)
    \item Training algorithm: Orthogonal greedy algorithm
    \item Parameter sampling: Randomly sample 10,000 parameter points from $[-1,1]^2$ for neural network parameter initialization
    \item Stopping tolerance: $\varepsilon=10^{-6}$.
\end{itemize}

\textbf{Test Data:}
\begin{itemize}
\item Test parameters: Randomly choose 5 from $\mu_s \in [0.4 \times 10^6, 0.6 \times 10^6]$, $\hat{U} \in [1.4, 1.6]$ as: $(\hat{U}, \mu_s) = (1.41, 0.45 \times 10^6)$, $(1.47, 0.55 \times 10^6)$, $(1.53, 0.51 \times 10^6)$, $(1.55, 0.49 \times 10^6)$, and $(1.59, 0.59 \times 10^6)$.
\end{itemize}

\subsubsection{Greedy Parameter-Time Selection: Numerical Details}

The greedy parameter-time selection algorithm is implemented following the methodology described in Section \ref{sec:offline}.
The greedy algorithm is applied independently to each scalar component in each of the 55 time segments. The candidate sample counts for each component are listed in Table~\ref{tab:training_sample_counts}. To illustrate the greedy selection results, we use segment 1, $(2,6]$s, segment 2, $(6,7]$s, and segment 6, $(10,10.1]$s, as representatives of the early transient, transition, and near steady-state phases, respectively. 

The full candidate/training sample count $N_{\Theta,t}$ and the selected reduced basis count $N_{\xi,g}$ describe different stages of the construction. Table~\ref{tab:training_sample_counts} reports the former, whereas Table~\ref{tab:greedy_basis_numbers} reports the latter for each scalar component in the three representative segments. Although all five components use the same number of candidate samples in a given segment, their selected bases and coefficient targets are constructed separately.

\begin{table}[H]
    \centering
    \begin{tabular}{|c|c|c|c|}
    \hline
    {Scalar component} & 
    {$(2,6]$s} & 
    {$(6,7]$s} &
    {$(10,10.1]$s} \\
    \hline
    Horizontal fluid velocity ($u_f$) & 500 & 600 & 600 \\
    \hline
    Vertical fluid velocity ($v_f$) & 500 & 600 & 600 \\
    \hline
    Horizontal structure displacement ($u_s$) & 300 & 400 & 400 \\
    \hline
    Vertical structure displacement ($v_s$) & 300 & 400 & 400 \\
    \hline
    Fluid pressure ($p_f$) & 1500 & 1800 & 2000 \\
    \hline
    \end{tabular}
    \caption{
    {Number of selected reduced basis functions $N_{\xi,g}$ for each scalar component in representative time segments. The corresponding candidate/training sample counts per component are 484,000, 121,000, and 12,100, respectively, as listed in Table~\ref{tab:training_sample_counts}.}} 
    \label{tab:greedy_basis_numbers}
\end{table}

\begin{figure}[H]
    \centering
    \subfigure[Fluid velocity]{
        \includegraphics[width = 0.3\textwidth]{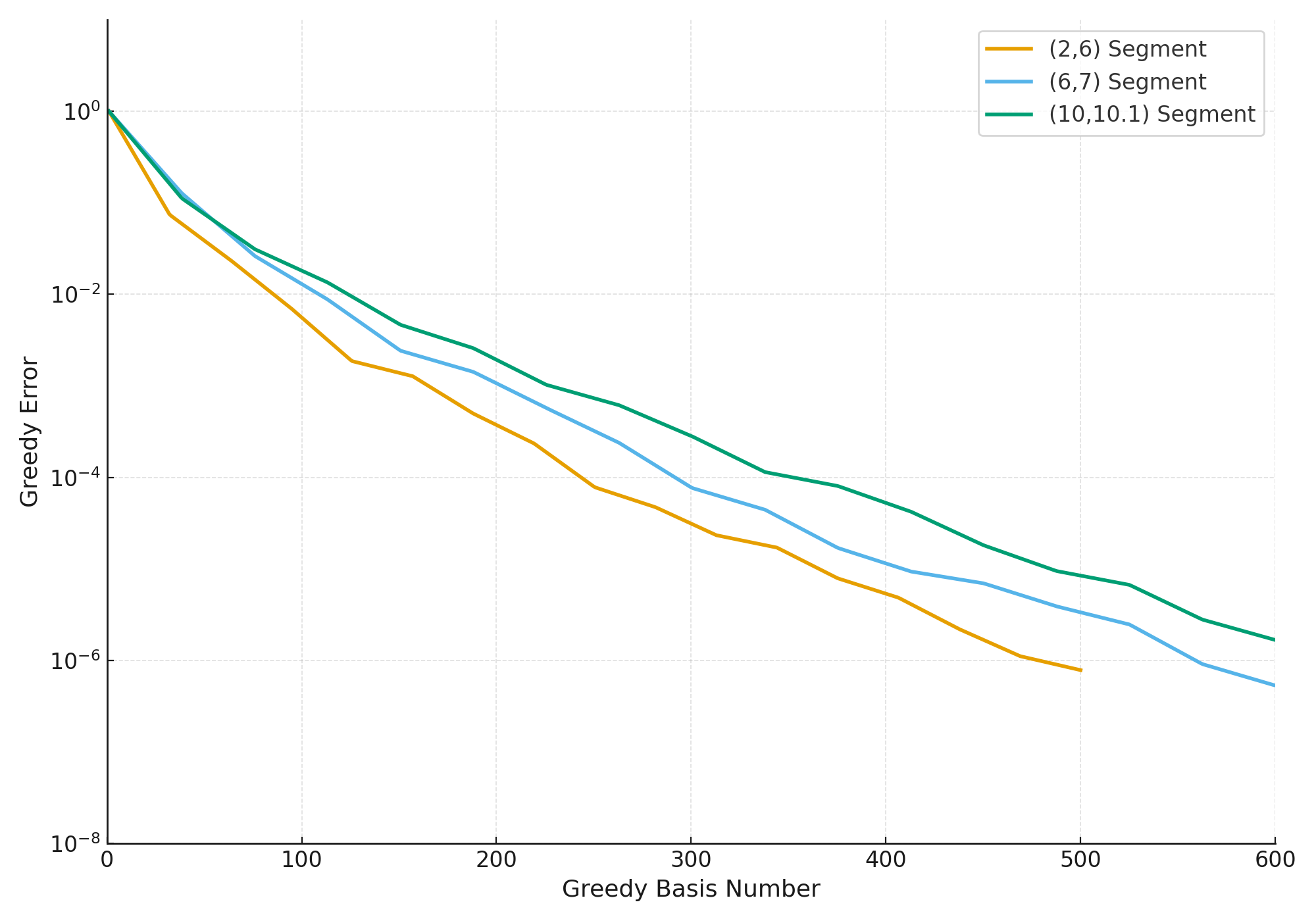}
        \label{fig:greedy_error_fluid}
    }
    \subfigure[Structural displacement]{
        \includegraphics[width = 0.3\textwidth]{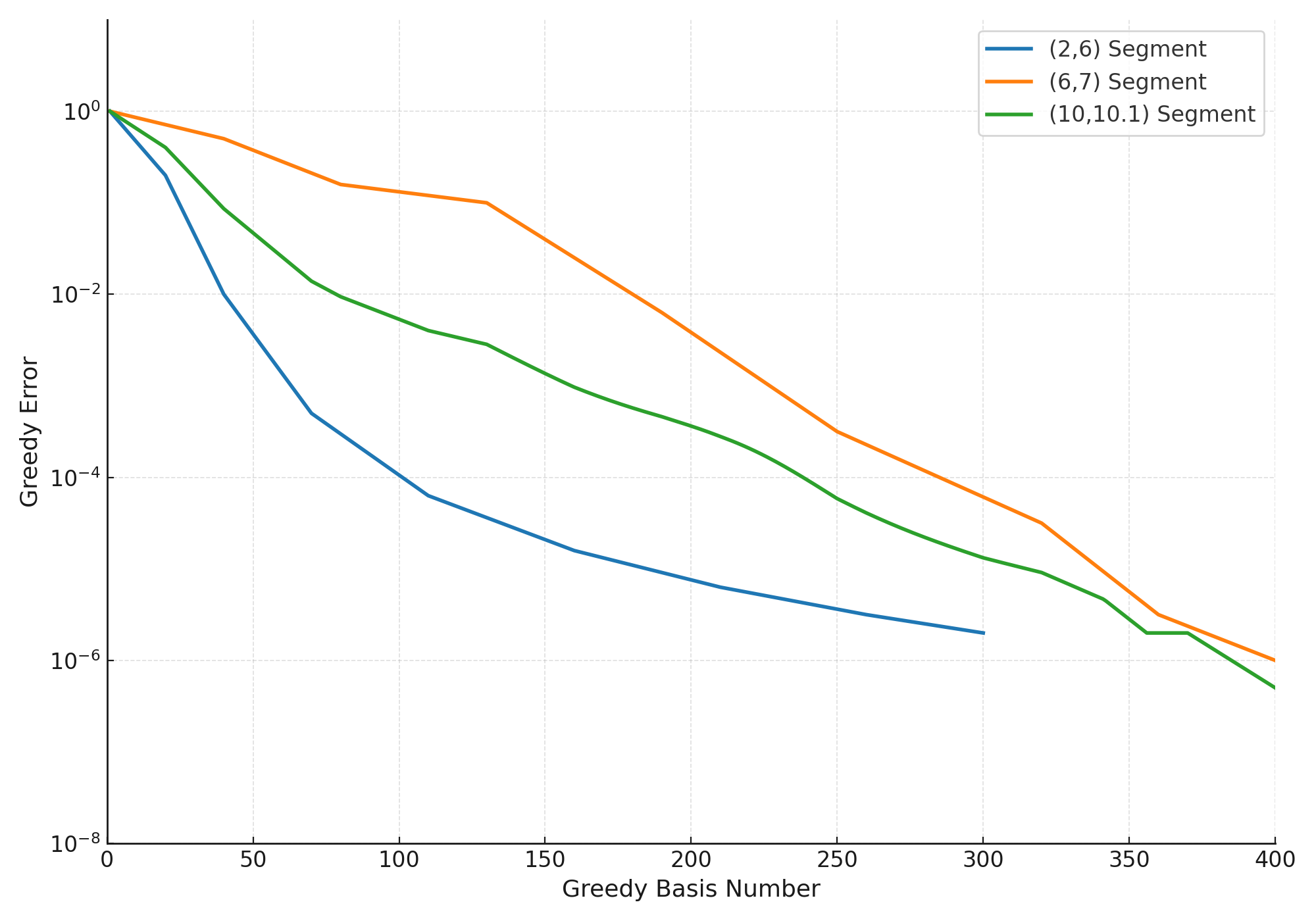}
        \label{fig:greedy_error_structure}
    }
    \subfigure[Fluid pressure]{
        \includegraphics[width = 0.3\textwidth]{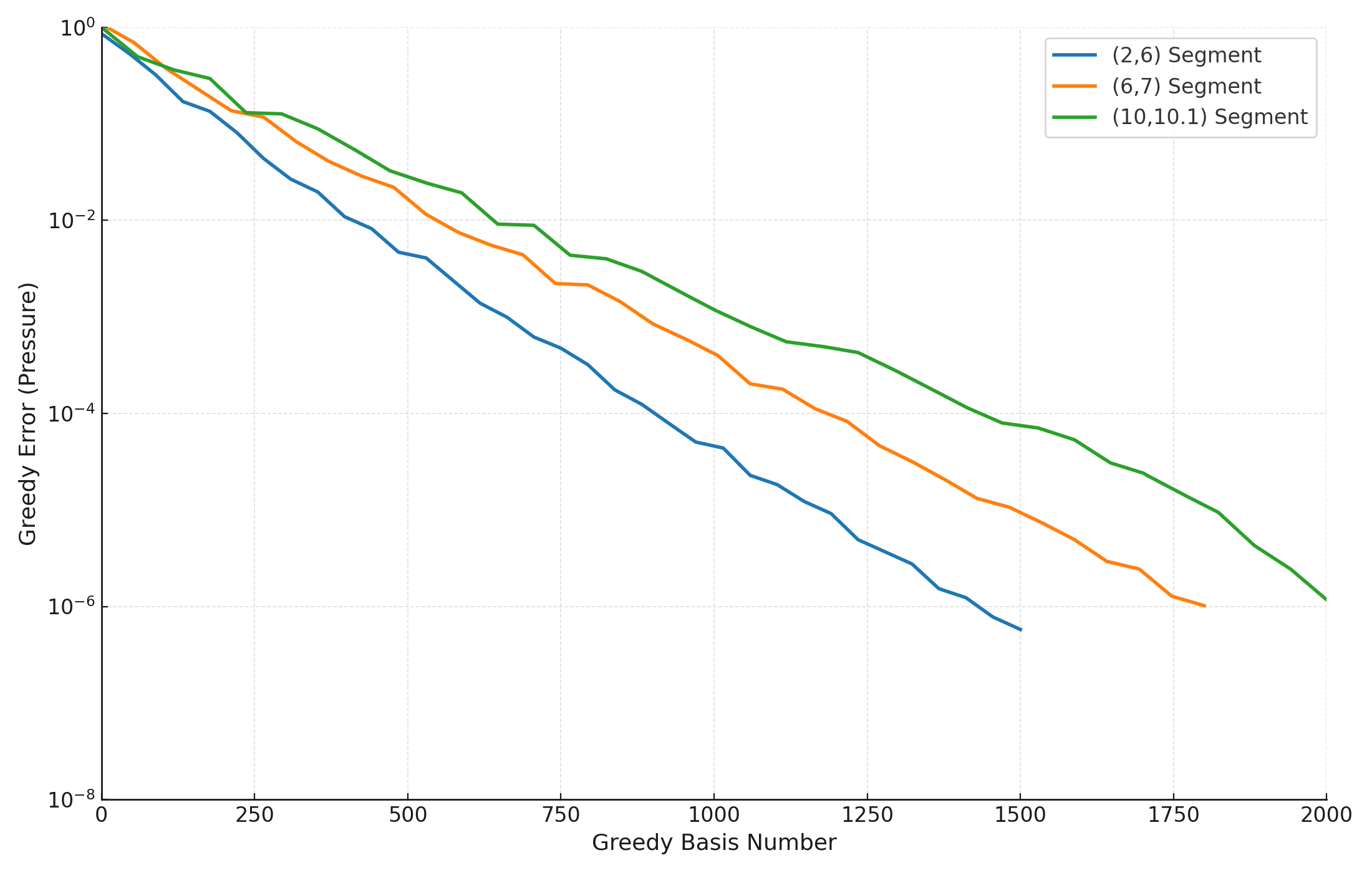}
        \label{fig:greedy_error_pressure}
    }
    \caption{Evolution of projection errors as the number of greedy basis functions increases: (a) Fluid velocity (horizontal components);  (b) Structural displacement (horizontal components); (c) Fluid pressure.}
    \label{fig:greedy_error_convergence}
\end{figure}

To demonstrate the convergence behavior of the greedy parameter-time selection algorithm, we present the evolution of projection errors as the number of greedy basis functions increases for these three representative time segments. Figure \ref{fig:greedy_error_convergence} shows the projection error decay for fluid velocity, pressure, and structural displacement, respectively.

The error convergence plots demonstrate several important characteristics of the greedy selection algorithm. First, the projection errors decrease monotonically as the number of basis functions increases, confirming the theoretical convergence properties of the greedy algorithm. For all variables, the greedy error reaches around $10^{-6}$ as the number of basis functions increases, demonstrating excellent convergence performance. Second, the error decay rates vary across different time segments, reflecting the varying solution complexity throughout the simulation. Third, the pressure field requires significantly more basis functions than other variables to achieve similar error levels, which is consistent with the higher variability observed in pressure distributions. Finally, the structural displacement exhibits the fastest convergence, indicating that the constrained beam motion can be efficiently represented with fewer basis functions.

\begin{figure}[H]
    \centering
    \includegraphics[width = 6cm, height = 4cm]{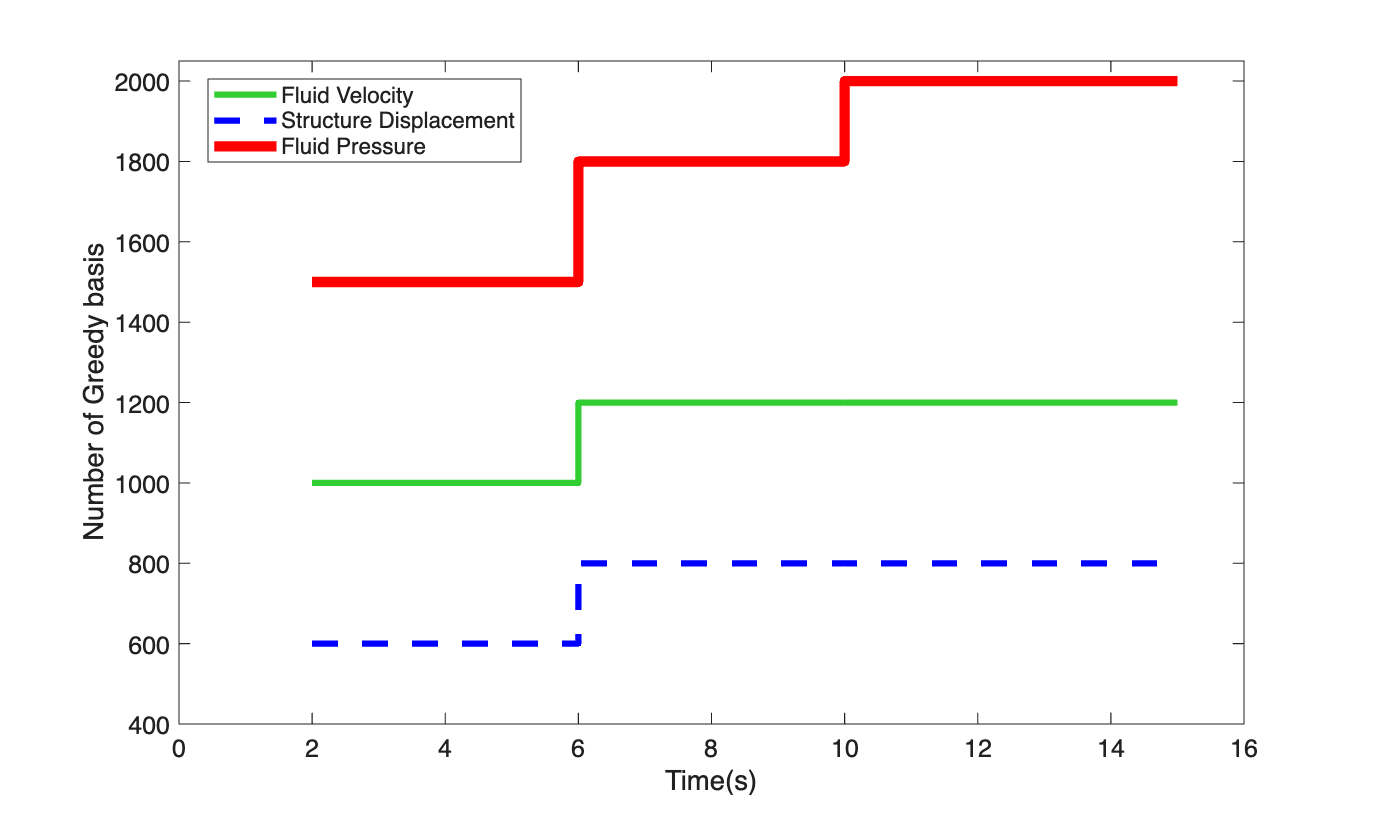}
    \caption{Evolution of greedy basis dimensions across time segments for different variables. }
    \label{fig:greedy_basis}
\end{figure}

Figure \ref{fig:greedy_basis} illustrates the evolution of the greedy basis dimensions across different time segments, demonstrating how the algorithm adaptively selects the number of basis functions based on the local solution complexity. The figure shows that the number of basis functions increases in later time segments, particularly for the pressure field, which reflects the increasing complexity of the flow patterns as the system evolves toward steady-state behavior.

\subsubsection{Shallow Neural Network Training with OGA}
After constructing the orthonormal bases, we project all sampled snapshots to prepare $\mathcal D_{\xi,g}$ and train one vector-valued shallow network $\bm f_{N,\xi,g}$ for each variable--segment pair.  OGA learns the mapping from parameters and time to the complete coefficient vector $\bm\alpha_{\xi,g}$.

To demonstrate the training performance, we present the training and test loss evolution for three representative time segments: segment 2 (time steps 6001-7000, corresponding to $[6.001, 7.000)$s), segment 3 (time steps 7001-8000, corresponding to $[7.001, 8.000)$s), and segment 6 (time steps 10001-10100, corresponding to $[10.001, 10.100)$s). These segments represent different phases of the simulation: early transition phase (segment 2), mid-transition phase (segment 3), and near steady-state phase (segment 6).

\begin{figure}[H]
    \centering
    \subfigure[Segment 2: Training loss]{
        \includegraphics[width = 0.3\textwidth]{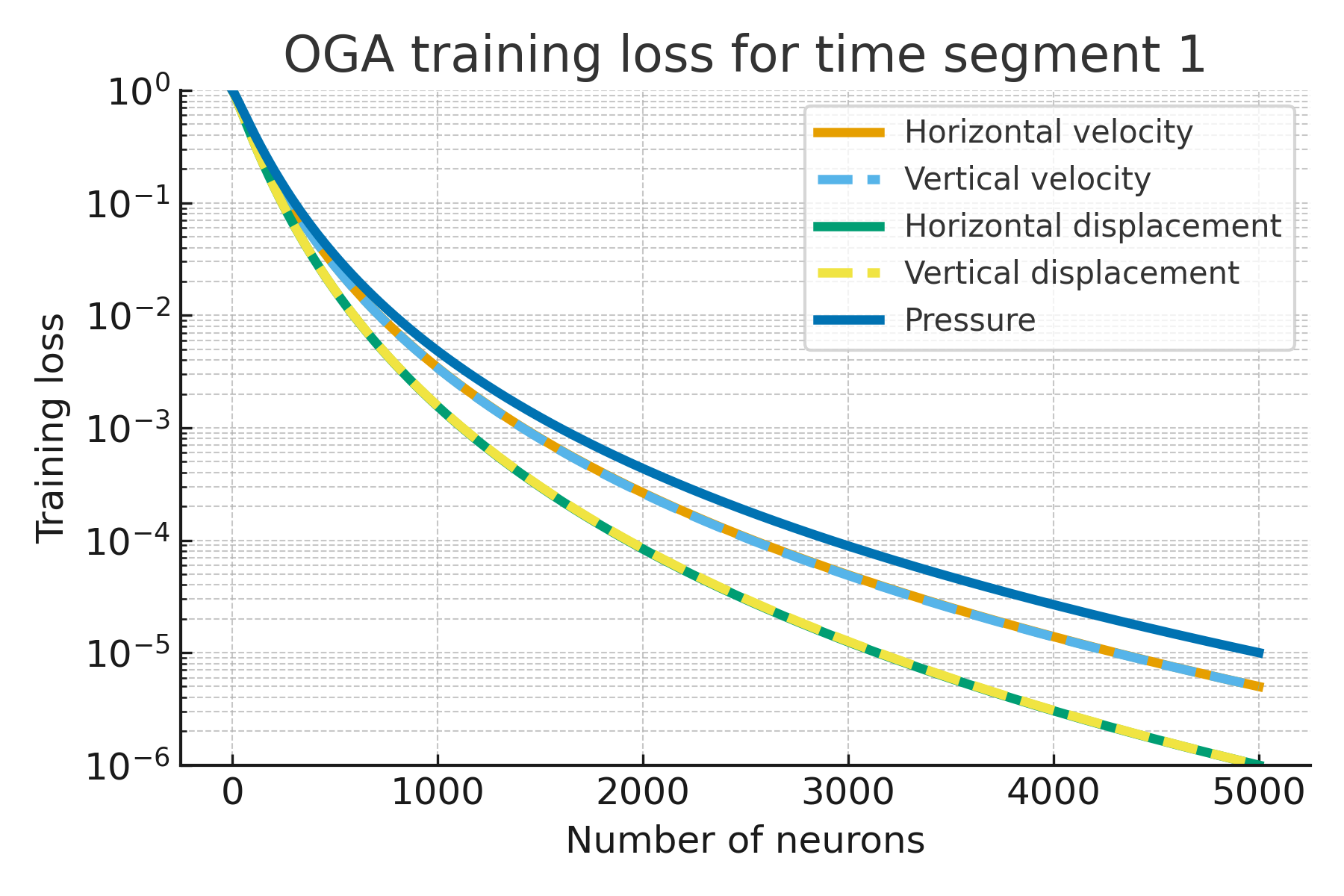}
        \label{fig:oga_train_seg1}
    }
    \subfigure[Segment 3: Training loss]{
        \includegraphics[width = 0.3\textwidth]{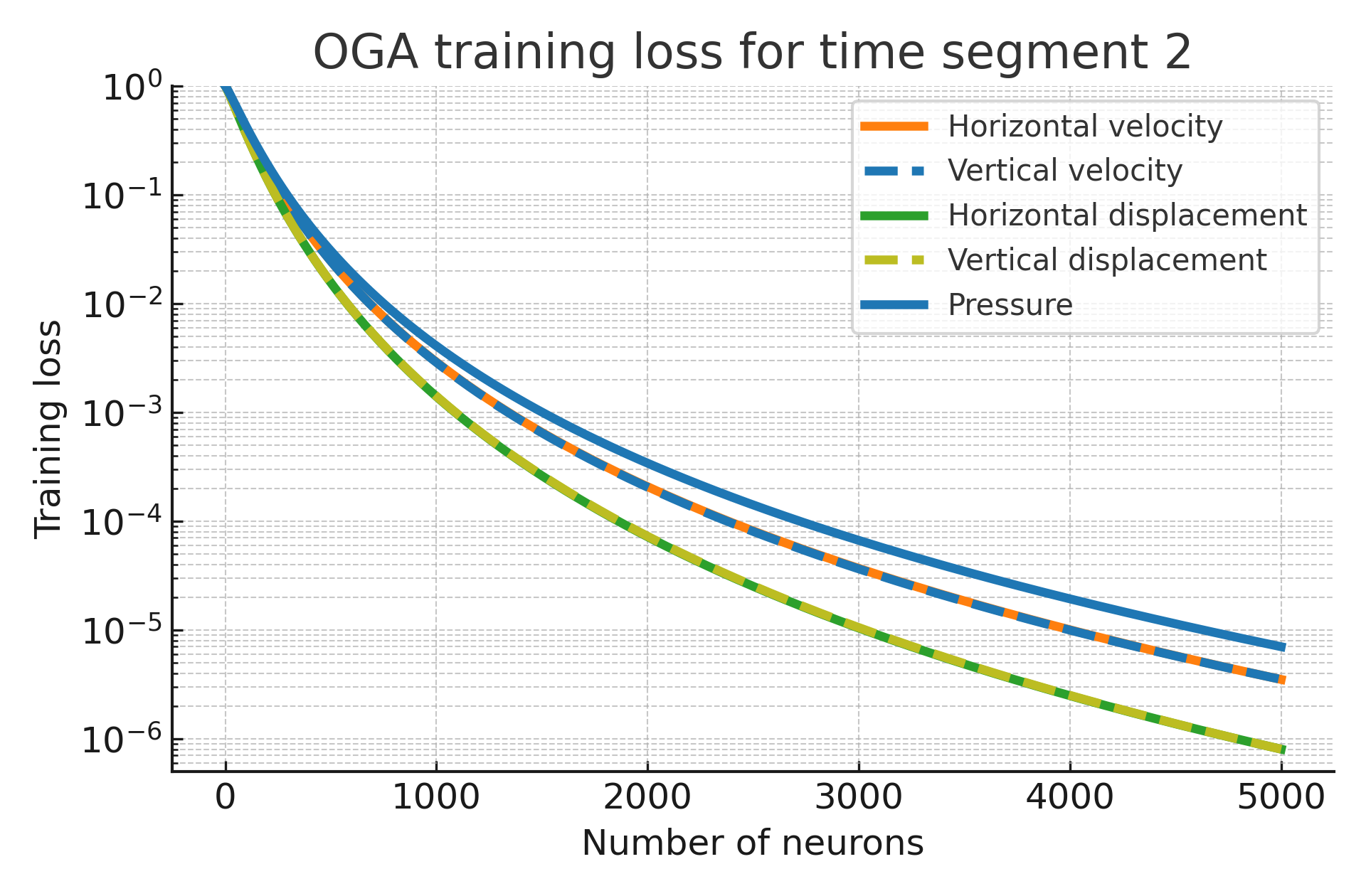}
        \label{fig:oga_train_seg2}
    }
    \subfigure[Segment 6: Training loss]{
        \includegraphics[width = 0.3\textwidth]{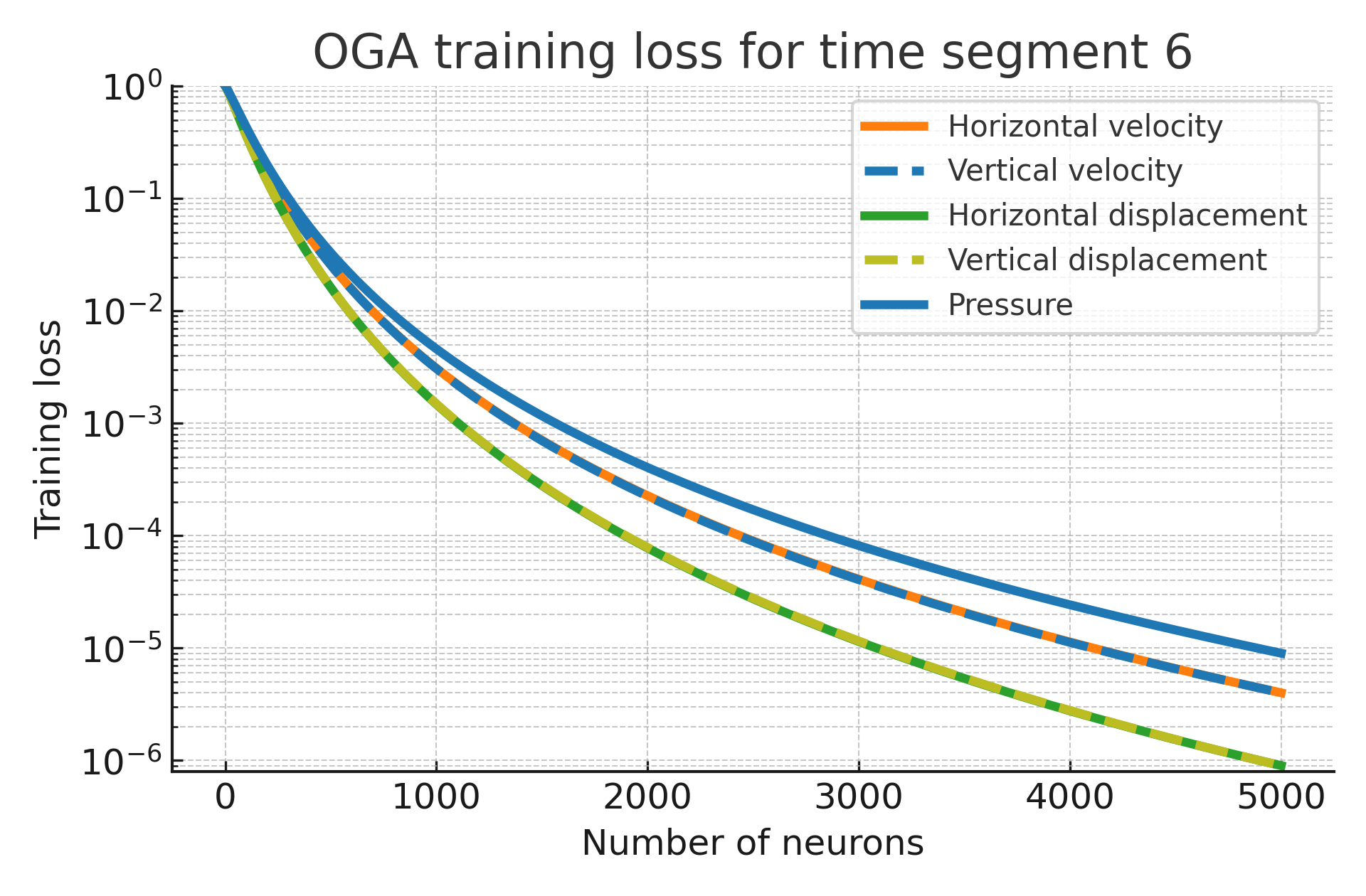}
        \label{fig:oga_train_seg6}
    }
    \caption{Training loss evolution for shallow neural networks trained using OGA for three representative time segments. }
    \label{fig:oga_training_loss}
\end{figure}

Figure \ref{fig:oga_training_loss} shows the training loss evolution for the three representative time segments. The training loss decreases monotonically as the OGA algorithm iteratively adds neurons to the network, demonstrating the convergence property of the orthogonal greedy algorithm. The loss curves show that the algorithm effectively learns the parameter-to-coefficient mapping for each time segment, with the training loss reaching very low values after 5000 training steps.

\begin{figure}[H]
    \centering
    \subfigure[Segment 2: Test loss]{
        \includegraphics[width = 0.3\textwidth]{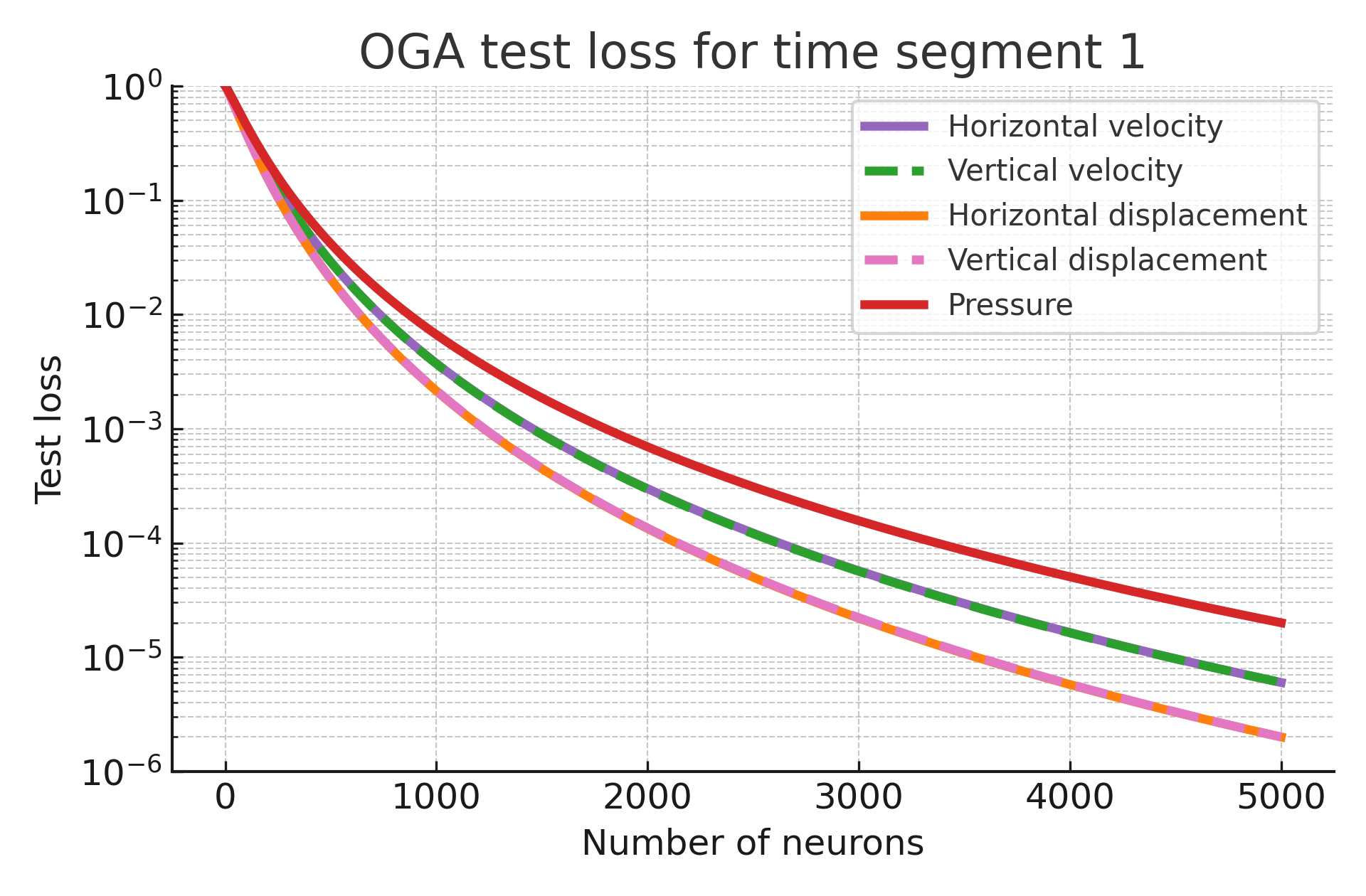}
        \label{fig:oga_test_seg1}
    }
    \subfigure[Segment 3: Test loss]{
        \includegraphics[width = 0.3\textwidth]{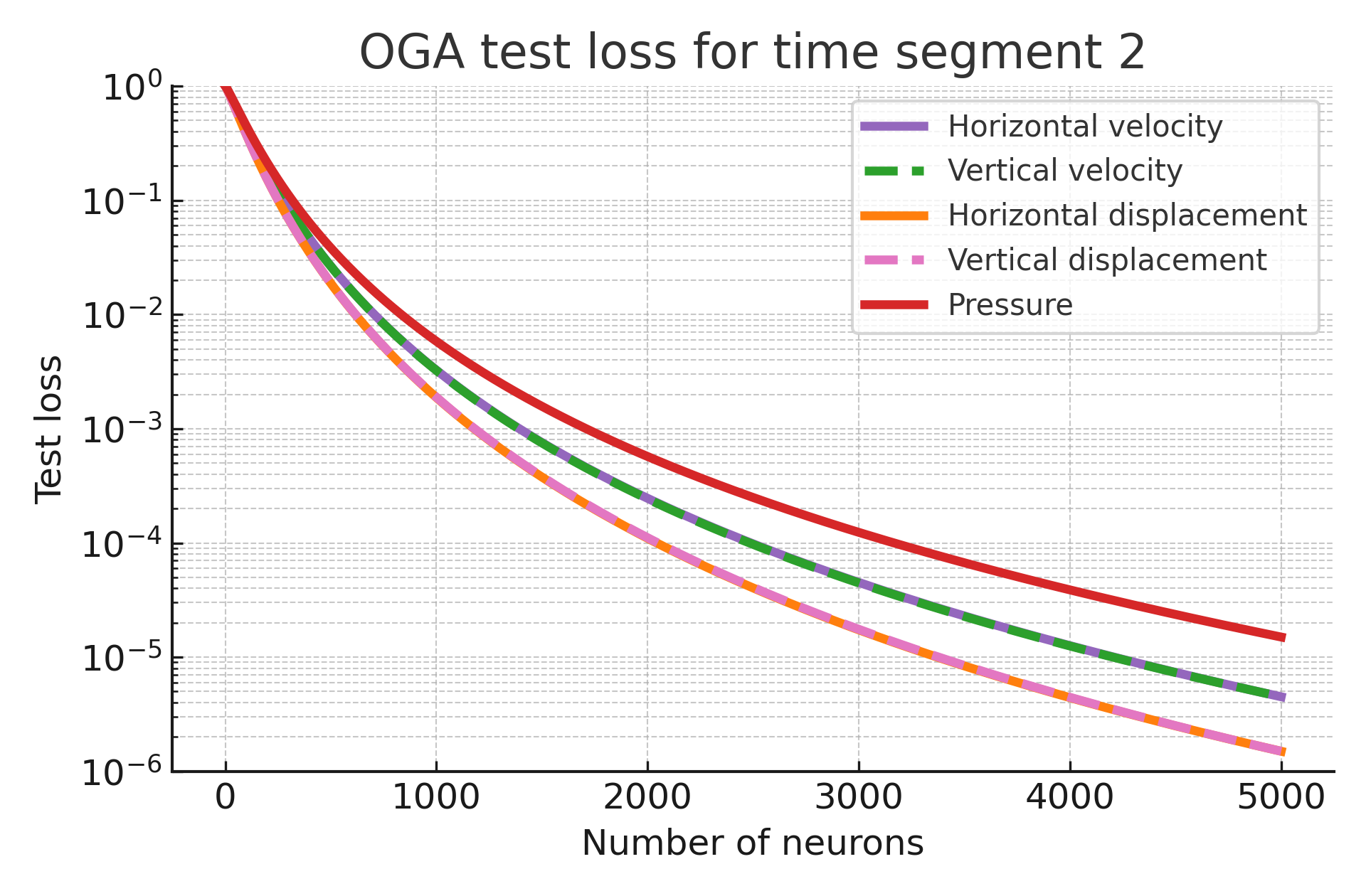}
        \label{fig:oga_test_seg2}
    }
    \subfigure[Segment 6: Test loss]{
        \includegraphics[width = 0.3\textwidth]{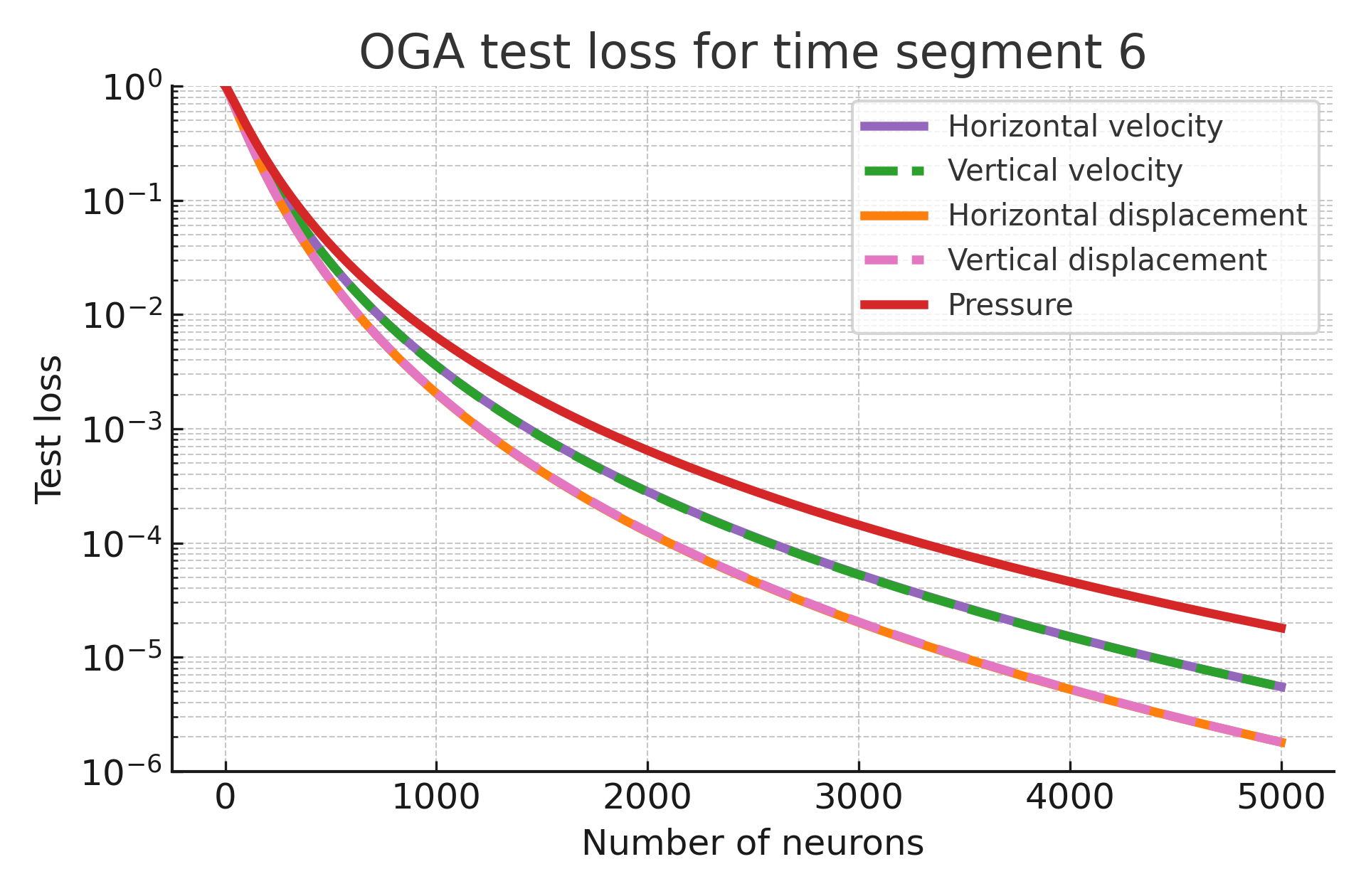}
        \label{fig:oga_test_seg6}
    }
    \caption{Test loss evolution for shallow neural networks trained using OGA for three representative time segments.}
    \label{fig:oga_test_loss}
\end{figure}

Figure \ref{fig:oga_test_loss} shows the test loss evolution for the same three representative time segments. The test loss follows a similar decreasing trend as the training loss, indicating that the trained networks generalize well to unseen parameter combinations. The close agreement between training and test losses demonstrates that the OGA training process avoids overfitting and produces robust neural network models that can accurately predict basis coefficients for new parameter values within the training parameter space.

The training results demonstrate that the OGA algorithm successfully trains shallow neural networks to approximate the complex nonlinear relationship between parameters, time, and basis coefficients. The monotonic decrease in both training and test losses validates the theoretical convergence guarantees of the orthogonal greedy algorithm, while the low final loss values ensure accurate coefficient predictions for the online phase of the ROM approach.

\subsubsection{ROM Accuracy Assessment}

To evaluate the accuracy of our ROM approach, we compare ROM predictions with FOM solutions for the same parameter combinations used in the FOM validation. We further introduce the following three error indicators that are used to compare the FOM solution $\bm{u}_{f,h}(t_n, \Theta)$ and the ROM solution $\hat{\bm{u}}_{f,h}(t_n, \Theta)$:

In these indicators, the lowercase FOM and ROM quantities are finite element functions obtained from their respective degrees-of-freedom vectors via~\eqref{eq:fe_reconstruction_operator}. Their spatial $L^2$ norms refer to integrals of the reconstructed functions, not unweighted Euclidean norms of the nodal vectors. For a common finite element basis, this distinction is expressed by $\|\mathcal R_{\xi,h}\bm e\|_{L^2}^2=\bm e^T M_\xi\bm e$, where $(M_\xi)_{ij}=(\psi_{\xi,i},\psi_{\xi,j})_{L^2}$. The same function-level interpretation applies to pressure and structural displacement. 

\begin{enumerate}
    \item The relative spatial $L^2$ error
    \begin{equation}
    \frac{\|\bm{u}_{f,h}(t_n, \Theta) - \hat{\bm{u}}_{f,h}(t_n, \Theta)\|_{L^2(\Omega_f(t_n))}}{\|\bm{u}_{f,h}(t_n, \Theta)\|_{L^2(\Omega_f(t_n))}};
    \end{equation}

    \item The relative spatio-temporal $L^2$ error
    \begin{equation}
    \frac{\|\bm{u}_{f,h}(t, \Theta) - \hat{\bm{u}}_{f,h}(t, \Theta)\|_{L^2(0,T;L^2(\Omega_f(t)))}}{\|\bm{u}_{f,h}(t, \Theta)\|_{L^2(0,T;L^2(\Omega_f(t)))}};
    \end{equation}

    \item The absolute error of $y$-displacement of point $A$ between the FOM result $d_{s,y}^{\text{FOM}}(A, t_n, \Theta)$ and the ROM result $d_{s,y}^{\text{ROM}}(A, t_n, \Theta)$, expressed as
    \begin{equation}
    d_{s,y}^{\text{FOM}}(A, t_n, \Theta) - d_{s,y}^{\text{ROM}}(A, t_n, \Theta),
    \end{equation}
    where $d_{s,y}^{\text{FOM}}(A, t_n, \Theta)$ and $d_{s,y}^{\text{ROM}}(A, t_n, \Theta)$ denote the $y$-component of the structural displacement at point $A$ at time $t_n$ for parameter $\Theta$ obtained from FOM and ROM, respectively.
\end{enumerate}

\begin{figure}[H]
    \centering
    \subfigure[\scriptsize $(\hat{U}, \mu_s) = (1.4, 0.4 \times 10^6)$]{
        \includegraphics[width = 0.3\textwidth]{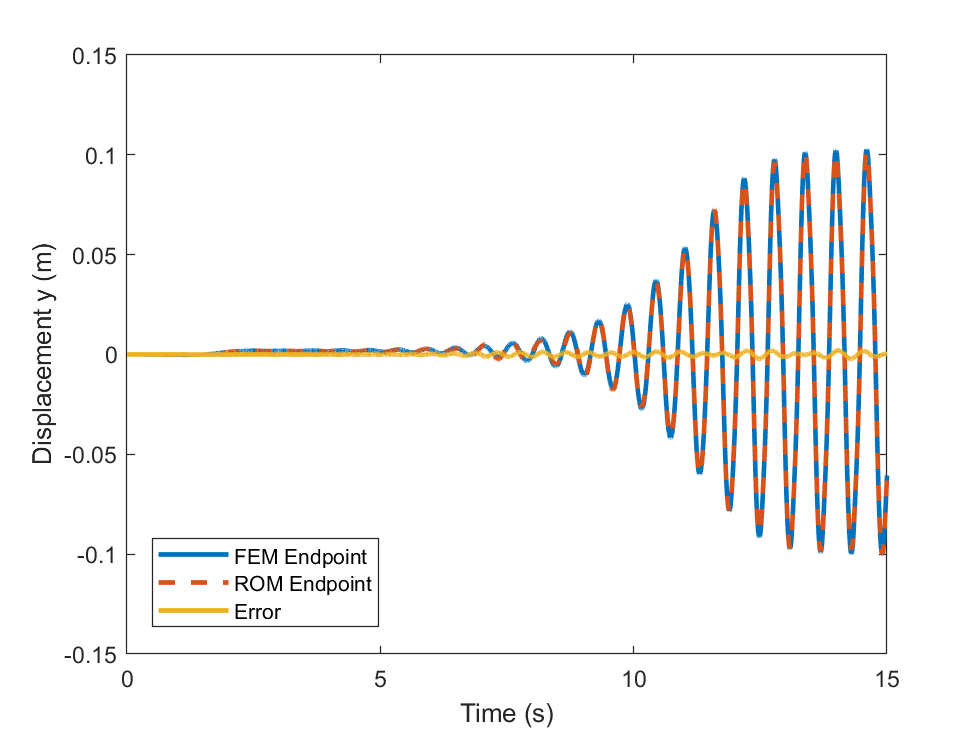}
        \label{fig:rom14040}
    }
    \subfigure[\scriptsize $(\hat{U}, \mu_s) = (1.46, 0.5 \times 10^6)$]{
        \includegraphics[width = 0.3\textwidth]{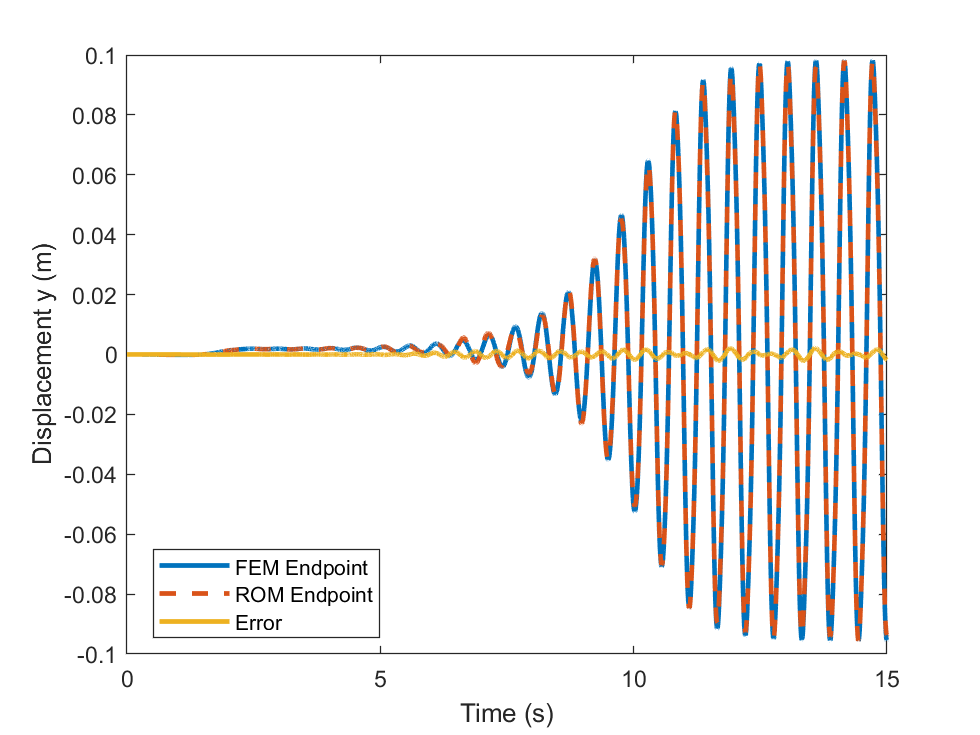}
        \label{fig:rom146050}
    }
    \subfigure[\scriptsize $(\hat{U}, \mu_s) = (1.50, 0.6 \times 10^6)$]{
        \includegraphics[width = 0.3\textwidth]{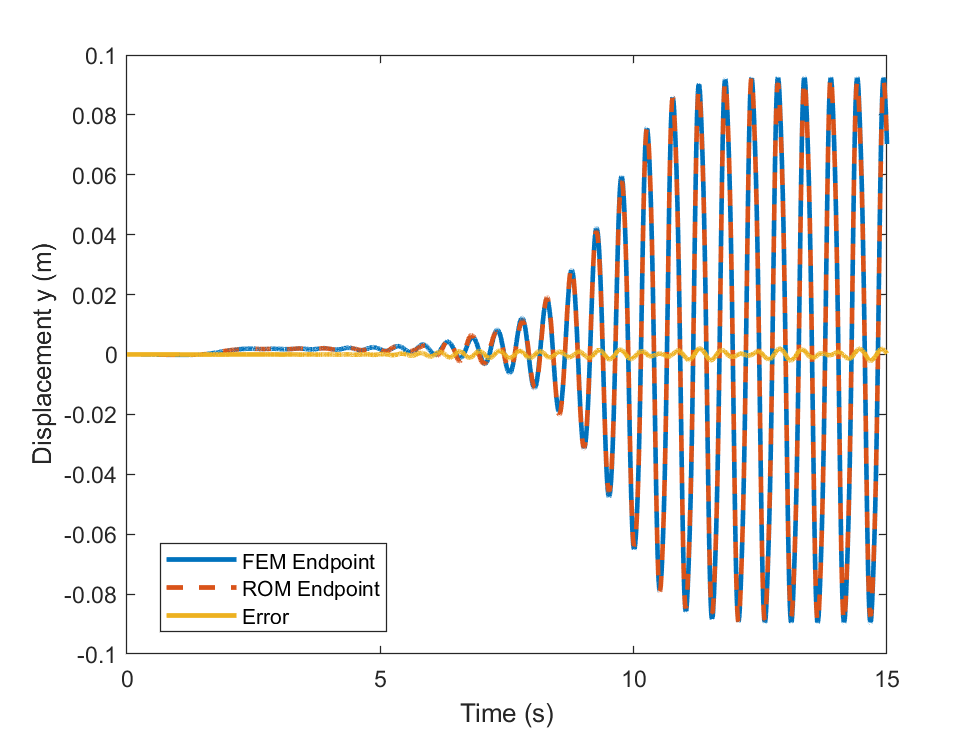}
        \label{fig:rom15060}
    }
    \subfigure[\scriptsize $(\hat{U}, \mu_s) = (1.54, 0.46 \times 10^6)$]{
        \includegraphics[width = 0.3\textwidth]{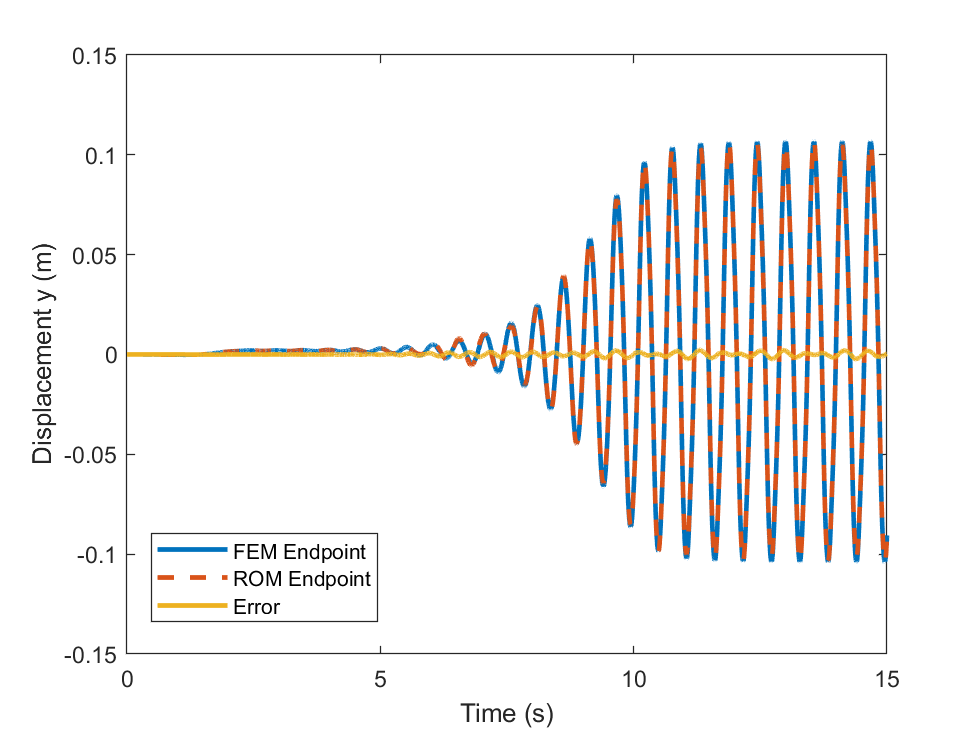}
        \label{fig:rom154046}
    }
    \subfigure[\scriptsize $(\hat{U}, \mu_s) = (1.6, 0.54 \times 10^6)$]{
        \includegraphics[width = 0.3\textwidth]{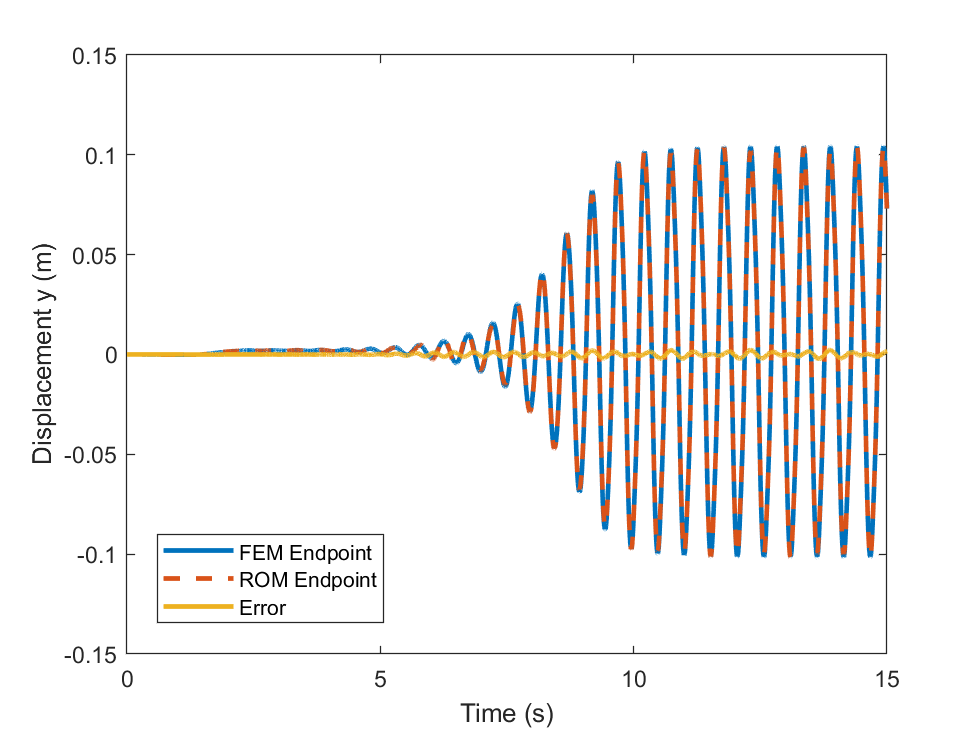}
        \label{fig:rom160054}
    }
    \caption{Comparison of FOM and ROM results for displacement response at location $A$. Each subfigure shows the FOM solution (solid blue), ROM prediction (dashed line), and their difference (solid yellow). The ROM accurately captures the oscillatory behavior and parametric dependencies of the FSI system.}
    \label{rom_fem_comparison}
\end{figure}

Figure \ref{rom_fem_comparison} presents a comparison between FOM and ROM results for the displacement response at location $A$ for the selected parameter combinations.

As observed in Figure \ref{rom_fem_comparison}, the ROM predictions closely match the FOM solutions, accurately capturing both the oscillatory behavior and the parametric dependencies. The error is on the order of $10^{-2}$, demonstrating the effectiveness of our approach in maintaining accuracy while significantly reducing computational complexity.

\begin{figure}[H]
    \centering
    \subfigure[\scriptsize $(\hat{U}, \mu_s) = (1.41, 0.45 \times 10^6)$]{
        \includegraphics[width = 0.3\textwidth]{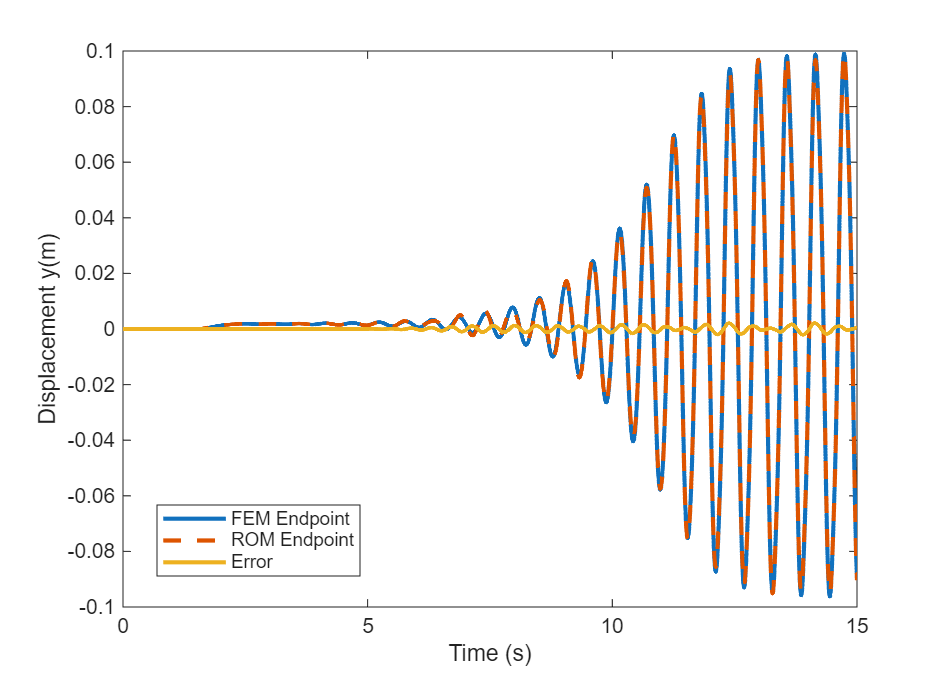}
        \label{fig:rom141045}
    }
    \subfigure[\scriptsize $(\hat{U}, \mu_s) = (1.47, 0.55 \times 10^6)$]{
        \includegraphics[width = 0.3\textwidth]{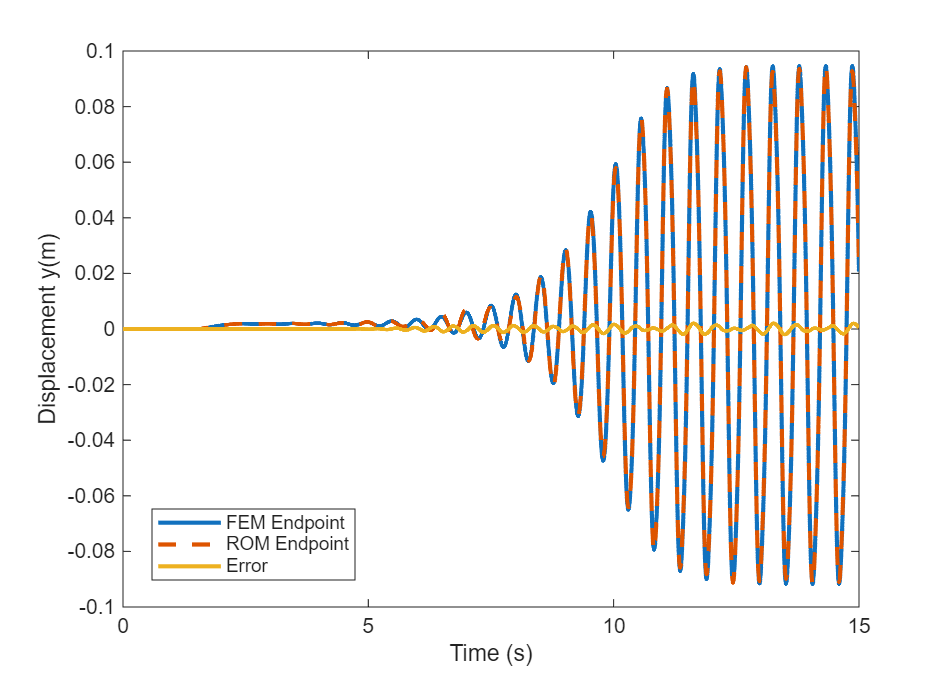}
        \label{fig:rom147055}
    }
    \subfigure[\scriptsize $(\hat{U}, \mu_s) = (1.53, 0.51 \times 10^6)$]{
        \includegraphics[width = 0.3\textwidth]{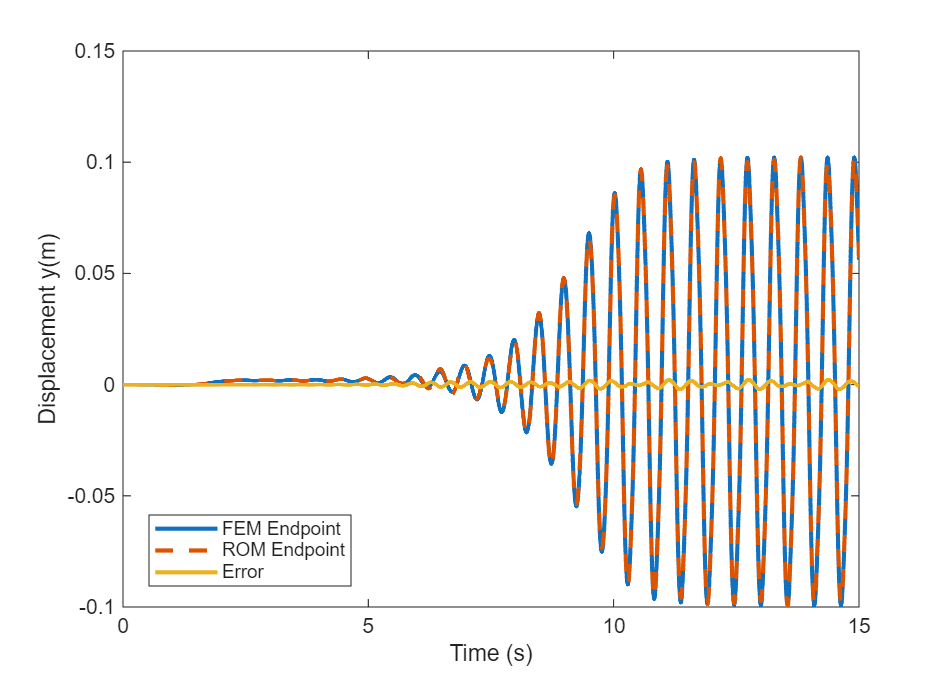}
        \label{fig:rom153051}
    }
    \subfigure[\scriptsize $(\hat{U}, \mu_s) = (1.55, 0.49 \times 10^6)$]{
        \includegraphics[width = 0.3\textwidth]{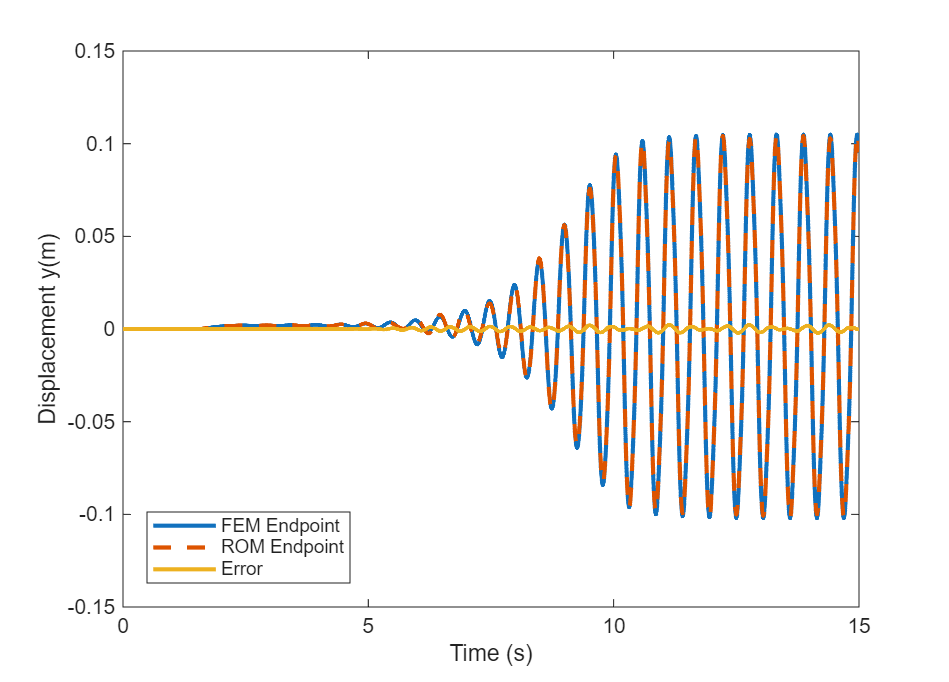}
        \label{fig:rom155049}
    }
    \subfigure[\scriptsize $(\hat{U}, \mu_s) = (1.59, 0.59 \times 10^6)$]{
        \includegraphics[width = 0.3\textwidth]{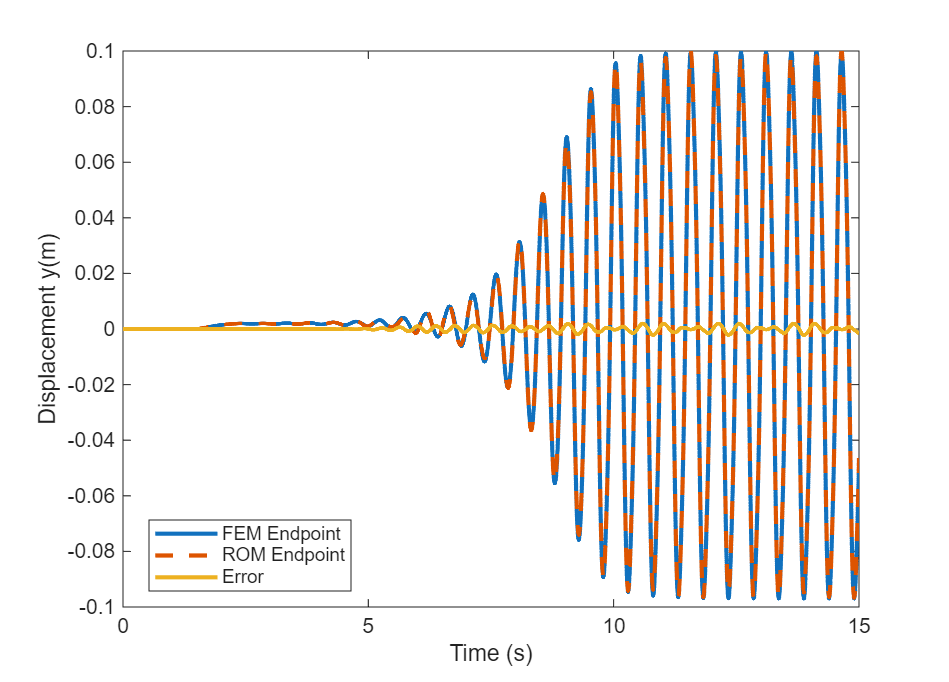}
        \label{fig:rom159059}
    }
    \caption{Comparison of FOM and ROM results for test parameter combinations. Each subfigure shows the FOM solution (solid line), ROM prediction (dashed line), and their difference. The ROM demonstrates excellent generalization capability, maintaining accuracy comparable to training set results.}
    \label{rom_fem_comparison_test}
\end{figure}

To further validate the ROM's generalization capability, we evaluate its performance on a test set of parameter combinations that were not included in the training data. Figure \ref{rom_fem_comparison_test} presents the comparison between FOM and ROM results for five test parameter combinations: $(\hat{U}, \mu_s) = (1.41, 0.45 \times 10^6)$, $(1.47, 0.55 \times 10^6)$, $(1.53, 0.51 \times 10^6)$, $(1.55, 0.49 \times 10^6)$, and $(1.59, 0.59 \times 10^6)$.

The test set results demonstrate that the ROM maintains high accuracy even for parameter combinations not seen during training, with errors remaining on the order of $10^{-2}$. This validates the robustness of our greedy parameter selection strategy and the generalization capability of the OGA-trained neural networks.

\subsubsection{Relative L2 Error Analysis}

\begin{figure}[H]
    \centering
    \subfigure[\scriptsize $(\hat{U}, \mu_s) = (1.4, 0.4 \times 10^6)$]{
        \includegraphics[width = 0.3\textwidth]{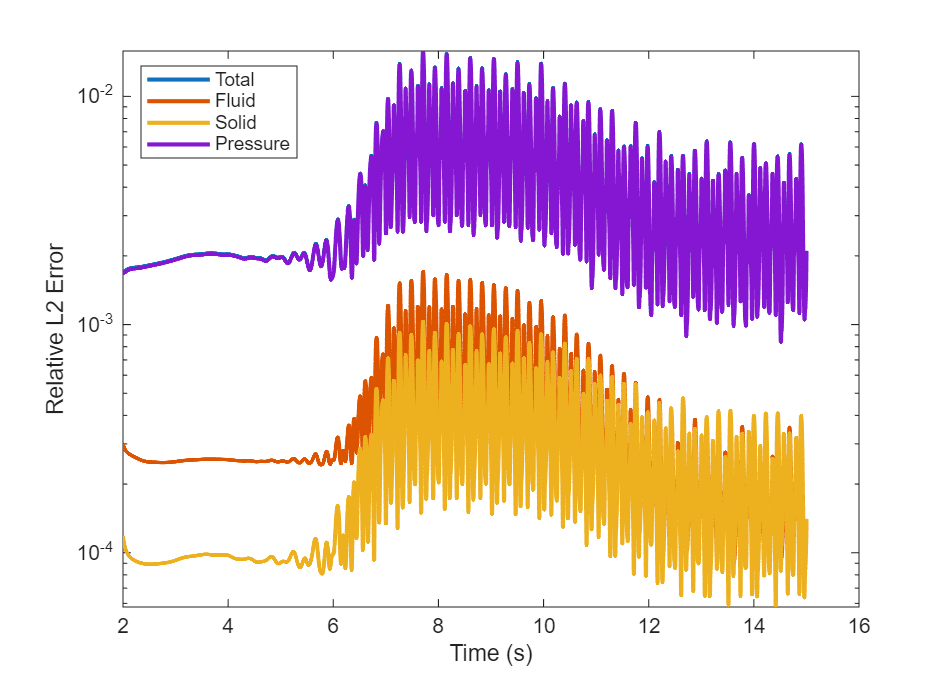}
        \label{fig:error14040}
    }
    \subfigure[\scriptsize $(\hat{U}, \mu_s) = (1.46, 0.5 \times 10^6)$]{
        \includegraphics[width = 0.3\textwidth]{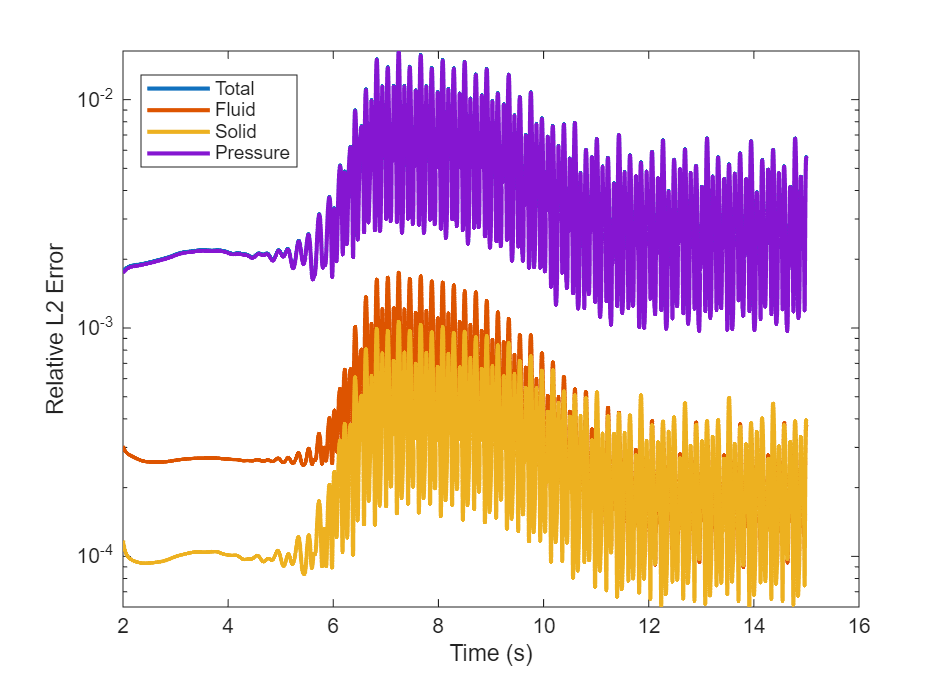}
        \label{fig:error146050}
    }
    \subfigure[\scriptsize $(\hat{U}, \mu_s) = (1.50, 0.6 \times 10^6)$]{
        \includegraphics[width = 0.3\textwidth]{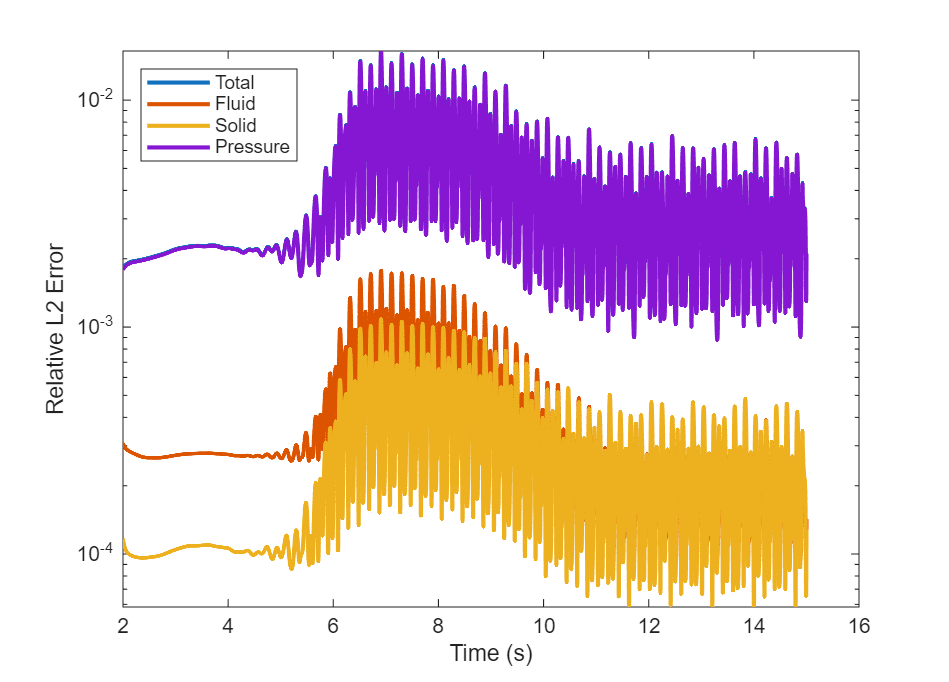}
        \label{fig:error15060}
    }
    \subfigure[\scriptsize $(\hat{U}, \mu_s) = (1.54, 0.46 \times 10^6)$]{
        \includegraphics[width = 0.3\textwidth]{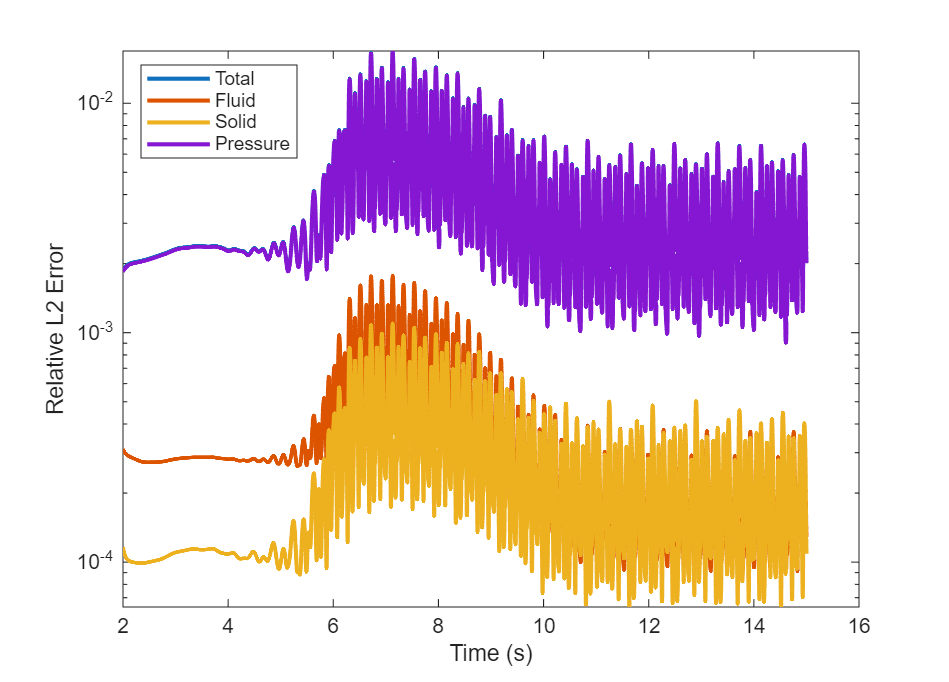}
        \label{fig:error154046}
    }
    \subfigure[\scriptsize $(\hat{U}, \mu_s) = (1.6, 0.54 \times 10^6)$]{
        \includegraphics[width = 0.3\textwidth]{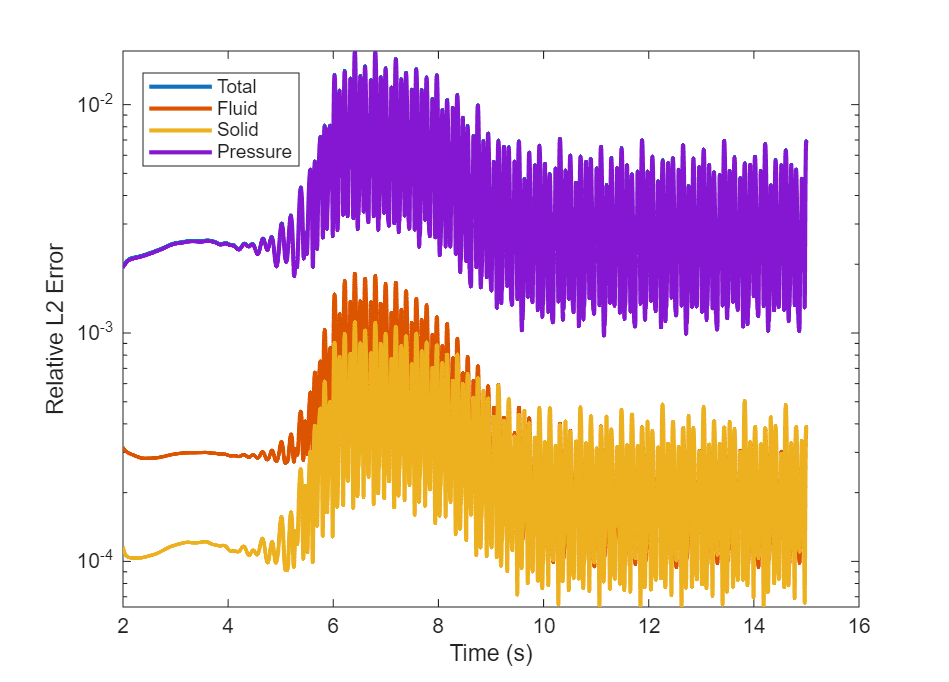}
        \label{fig:error160054}
    }
    \caption{Relative $L^2$ errors for ROM predictions compared to FOM solutions. Each subfigure shows the relative errors for total solution, fluid velocity, structural displacement, and pressure over time. The errors remain consistently on the order of $10^{-2}$ throughout the simulation, demonstrating the robustness of the ROM approach.}
    \label{relative_l2_errors}
\end{figure}

To provide a comprehensive assessment of ROM accuracy across different variables and parameter combinations, we analyze the relative $L^2$ errors for the total solution, fluid velocity, structural displacement, and pressure. Figure \ref{relative_l2_errors} presents the relative $L^2$ error evolution over time for the selected parameter combinations.

The relative $L^2$ error analysis reveals several important characteristics:
\begin{itemize}
    \item \textbf{Total solution error}: The relative $L^2$ error for the total solution remains consistently around $1-2\%$ across all parameter combinations and time steps, indicating excellent overall accuracy.
    \item \textbf{Fluid velocity error}: The fluid velocity components exhibit errors on the order of $10^{-2}$, with slightly higher errors during transient phases and lower errors in steady-state regimes.
    \item \textbf{Structural displacement error}: The structural displacement errors are well-controlled, typically below $2\%$, demonstrating that the ROM accurately captures the structural dynamics.
    \item \textbf{Pressure error}: The pressure field shows comparable accuracy, with relative errors maintained around $1-2\%$, which is particularly important for FSI problems where pressure plays a crucial role in the coupling mechanism.
\end{itemize}

The consistent error levels across different variables and parameter combinations validate the effectiveness of our variable-specific ROM construction approach, where separate Gram-Schmidt orthonormal bases and neural networks are trained for each physical quantity.

\begin{figure}[H]
    \centering
    \subfigure[\scriptsize $(\hat{U}, \mu_s) = (1.41, 0.45 \times 10^6)$]{
        \includegraphics[width = 0.3\textwidth]{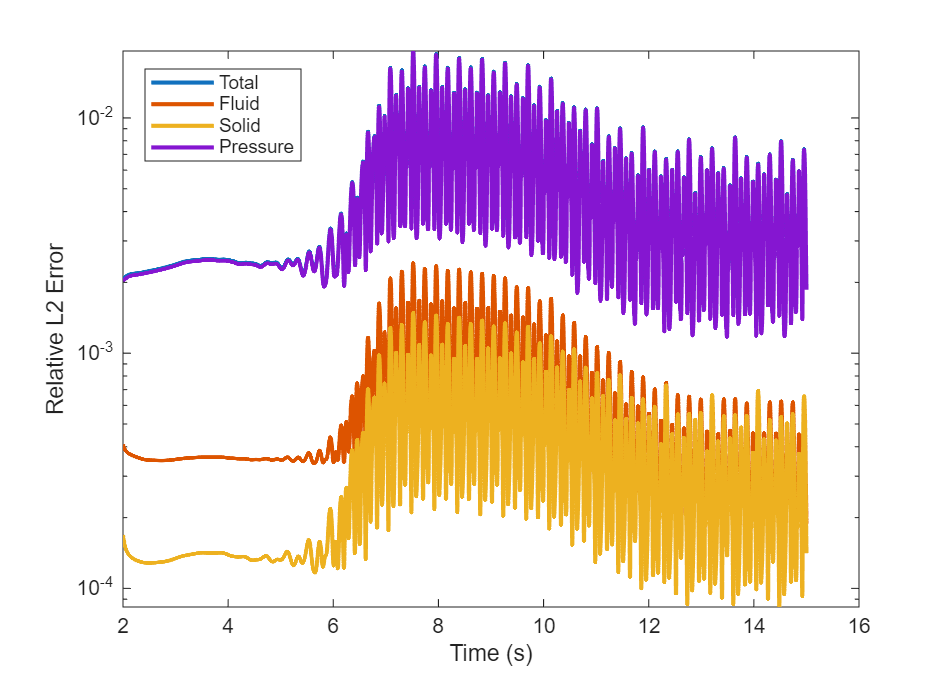}
        \label{fig:error141045}
    }
    \subfigure[\scriptsize $(\hat{U}, \mu_s) = (1.47, 0.55 \times 10^6)$]{
        \includegraphics[width = 0.3\textwidth]{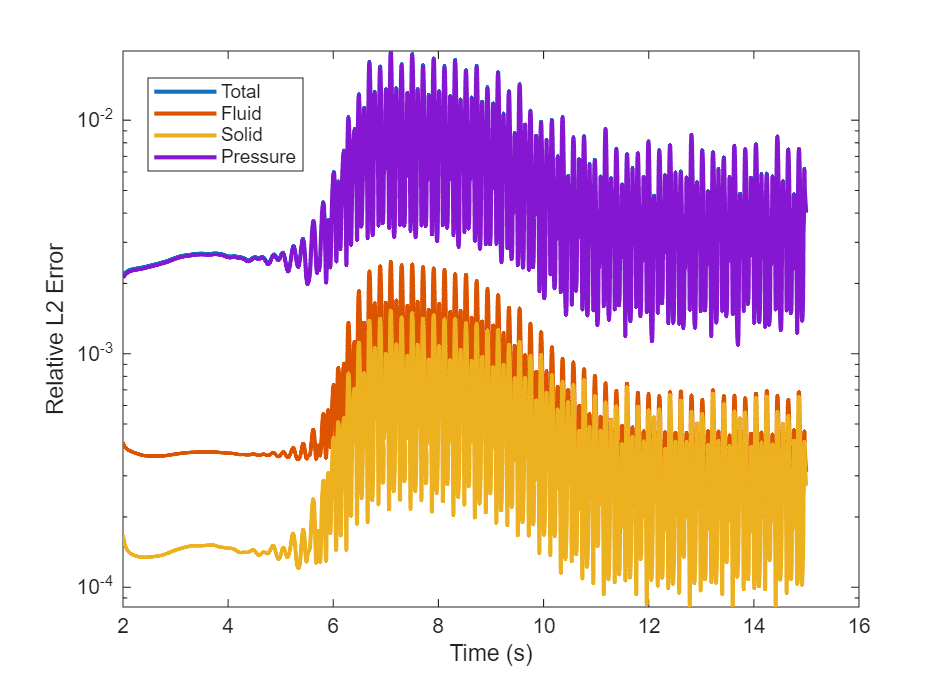}
        \label{fig:error147055}
    }
    \subfigure[\scriptsize $(\hat{U}, \mu_s) = (1.53, 0.51 \times 10^6)$]{
        \includegraphics[width = 0.3\textwidth]{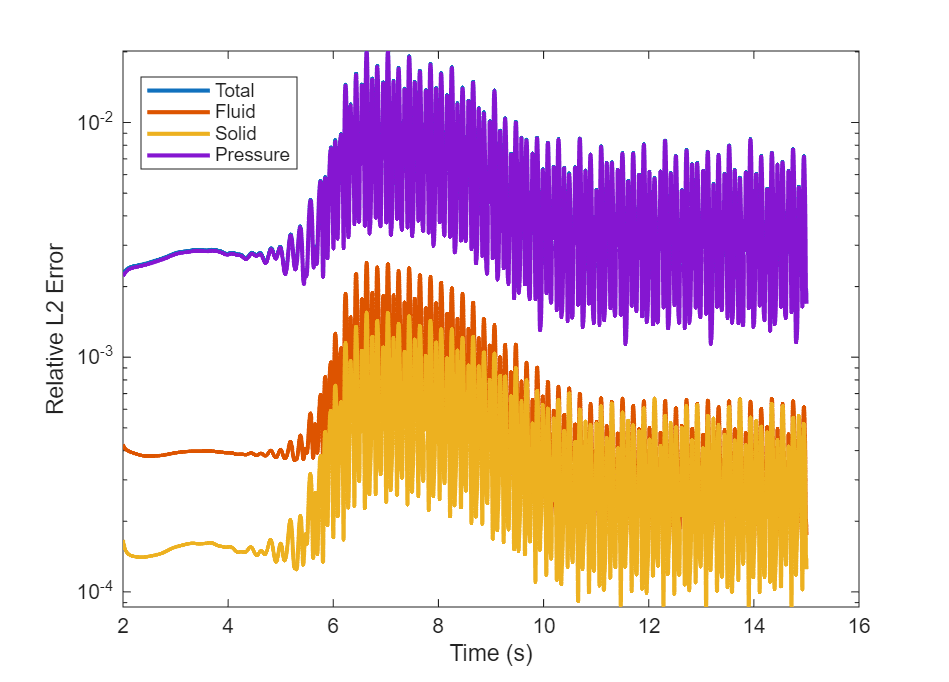}
        \label{fig:error153051}
    }
    \subfigure[\scriptsize $(\hat{U}, \mu_s) = (1.55, 0.49 \times 10^6)$]{
        \includegraphics[width = 0.3\textwidth]{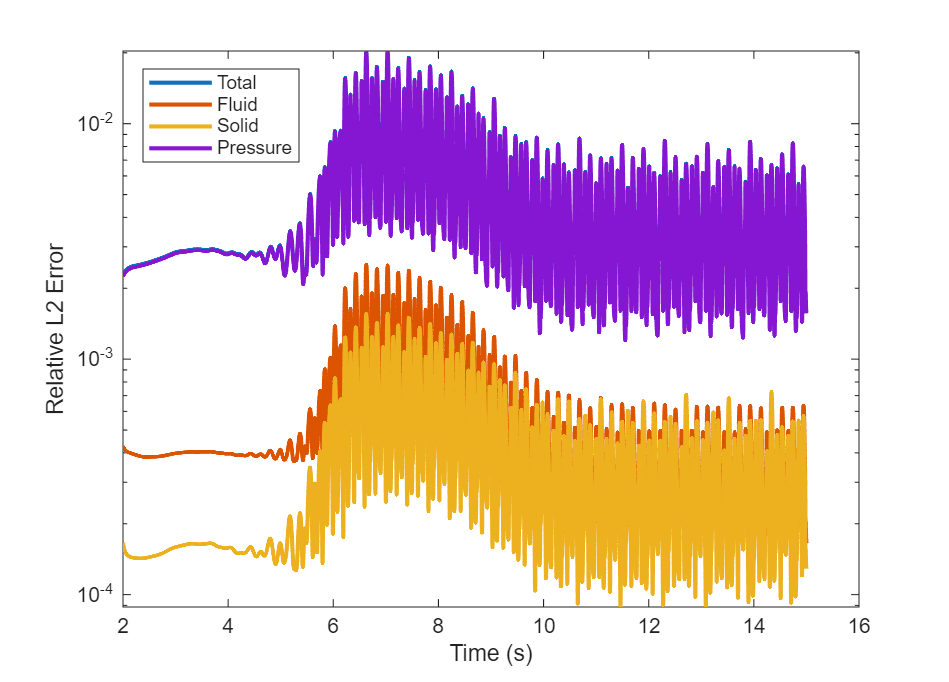}
        \label{fig:error155049}
    }
    \subfigure[\scriptsize $(\hat{U}, \mu_s) = (1.59, 0.59 \times 10^6)$]{
        \includegraphics[width = 0.3\textwidth]{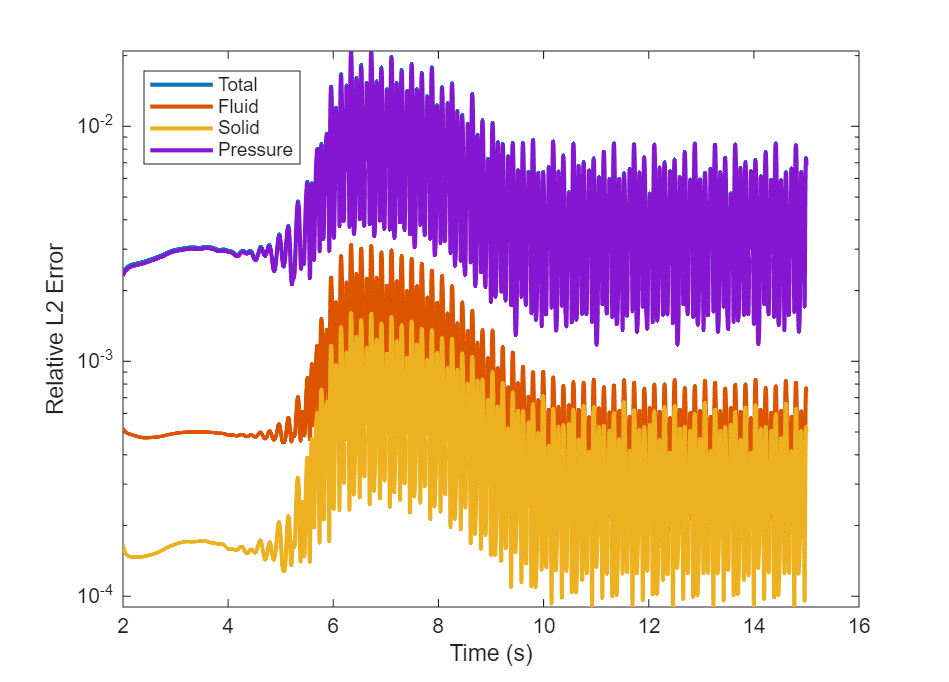}
        \label{fig:error159059}
    }
    \caption{Relative $L^2$ errors for test parameter combinations. The ROM maintains consistent error levels around $1-2\%$ for test parameters, demonstrating robust generalization performance.}
    \label{relative_l2_errors_test}
\end{figure}

Figure \ref{relative_l2_errors_test} presents the relative $L^2$ error evolution for the test parameter combinations, demonstrating that the ROM maintains comparable accuracy levels for unseen parameter values.

\subsubsection{Spatial Distribution Comparison at Final Time}

\begin{figure}[H]
    \centering
    \subfigure[\scriptsize $(\hat{U}, \mu_s) = (1.4, 0.4 \times 10^6)$: FOM $u_1$]{
        \includegraphics[width = 0.4\textwidth]{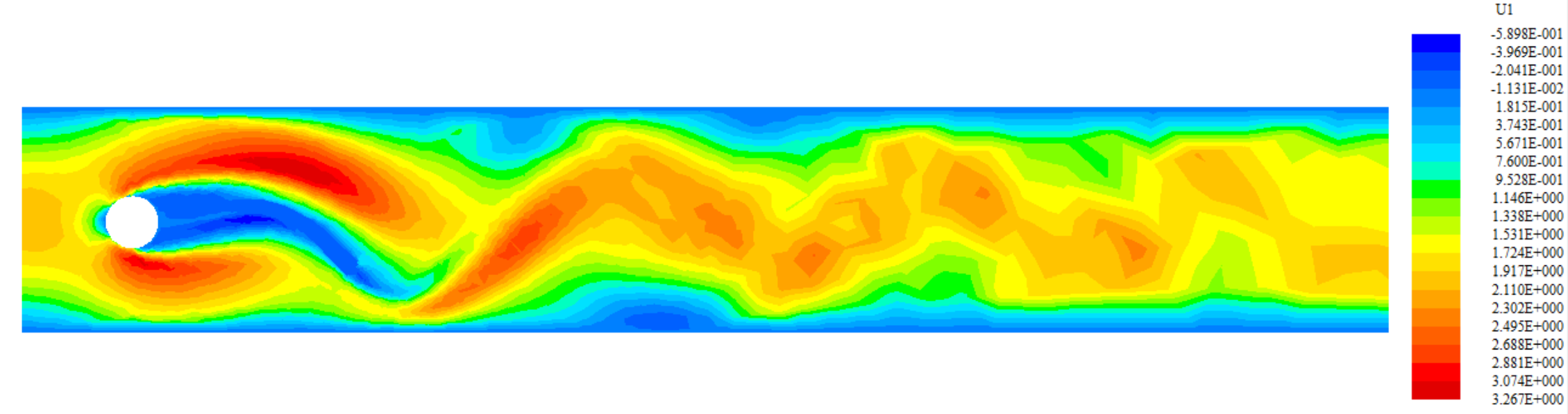}
    }
    \subfigure[\scriptsize $(\hat{U}, \mu_s) = (1.4, 0.4 \times 10^6)$: ROM $u_1$]{
        \includegraphics[width = 0.4\textwidth]{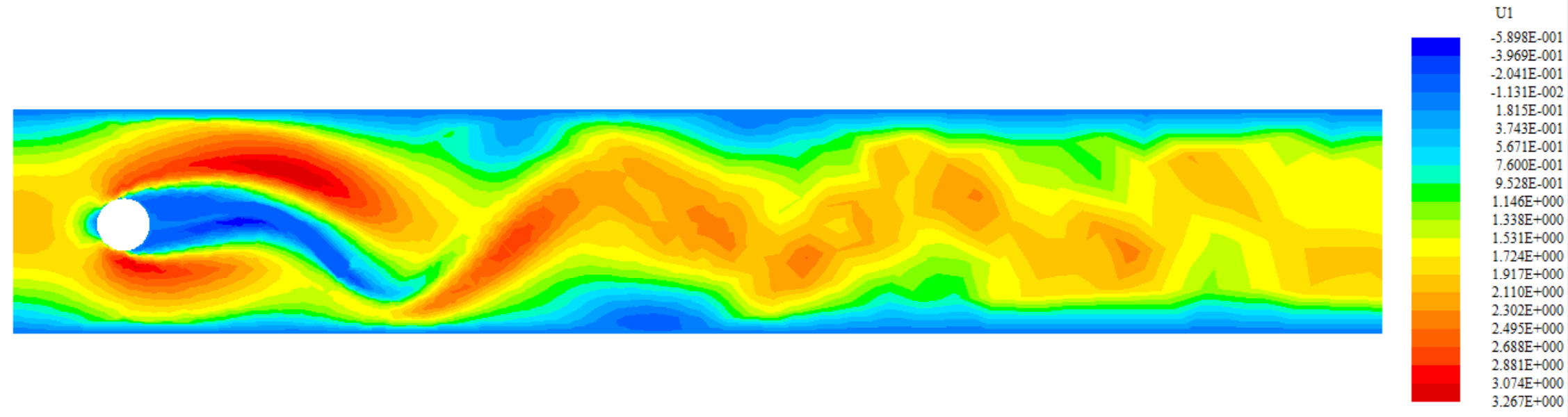}
    }
    \subfigure[\scriptsize $(\hat{U}, \mu_s) = (1.46, 0.5 \times 10^6)$: FOM $u_1$]{
        \includegraphics[width = 0.4\textwidth]{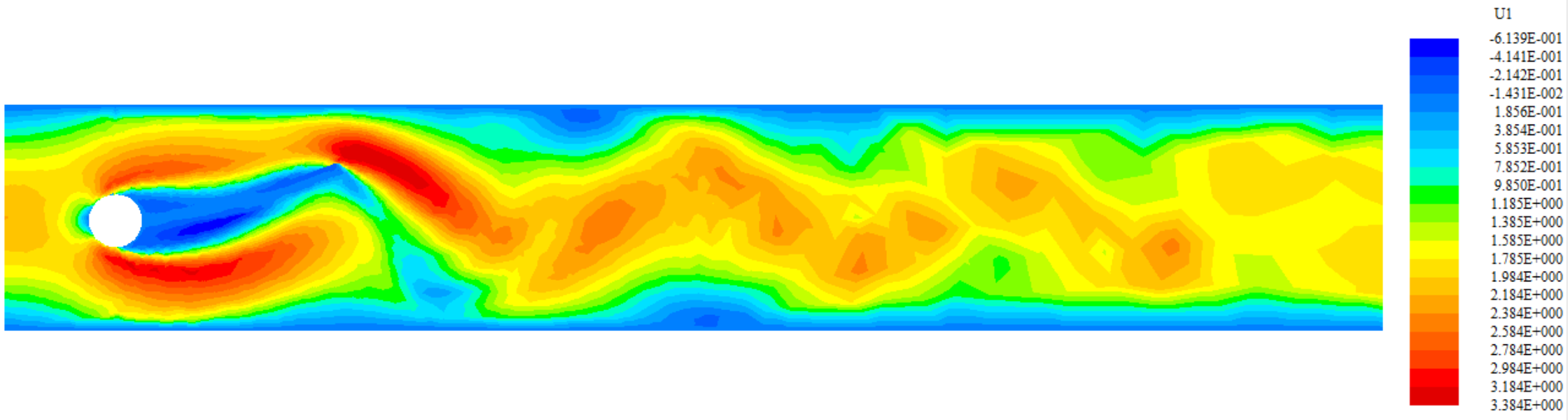}
    }
    \subfigure[\scriptsize $(\hat{U}, \mu_s) = (1.46, 0.5 \times 10^6)$: ROM $u_1$]{
        \includegraphics[width = 0.4\textwidth]{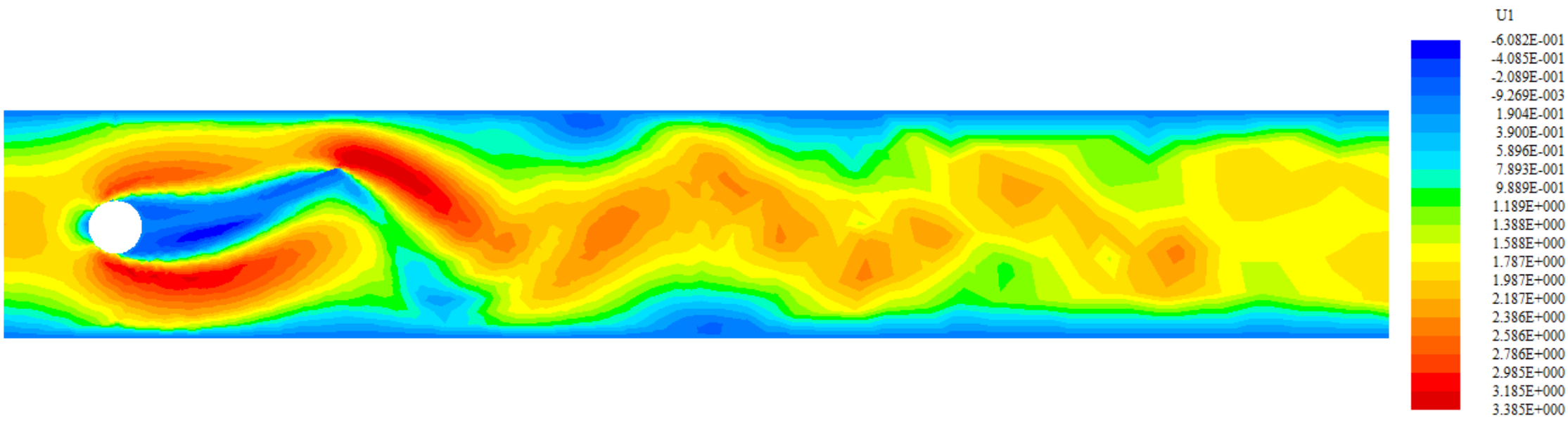}
    }
    \subfigure[\scriptsize $(\hat{U}, \mu_s) = (1.50, 0.6 \times 10^6)$: FOM $u_1$]{
        \includegraphics[width = 0.4\textwidth]{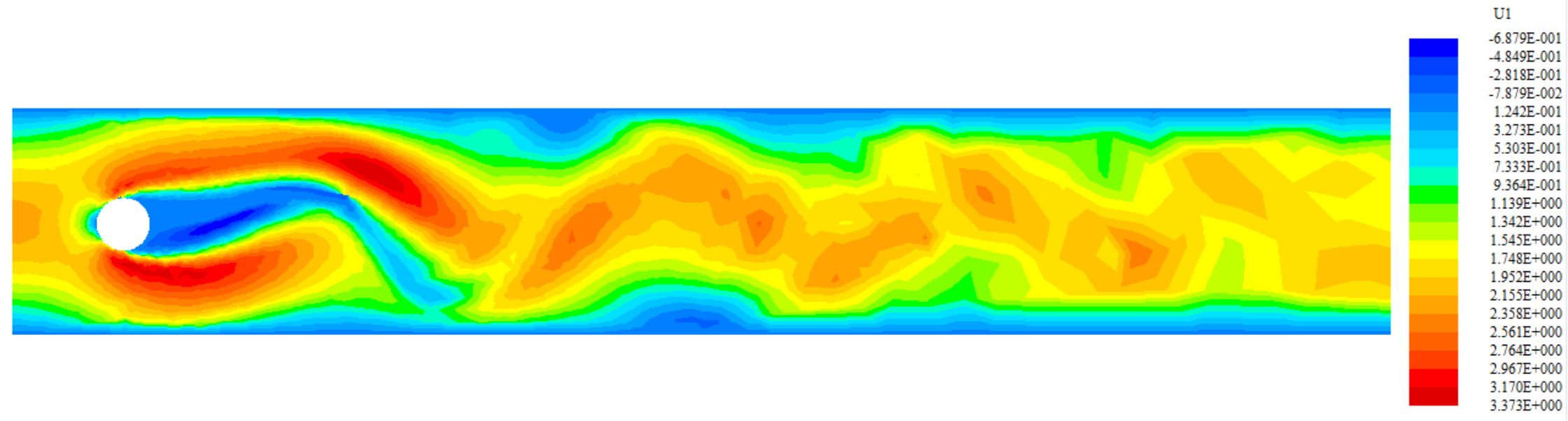}
    }
    \subfigure[\scriptsize $(\hat{U}, \mu_s) = (1.50, 0.6 \times 10^6)$: ROM $u_1$]{
        \includegraphics[width = 0.4\textwidth]{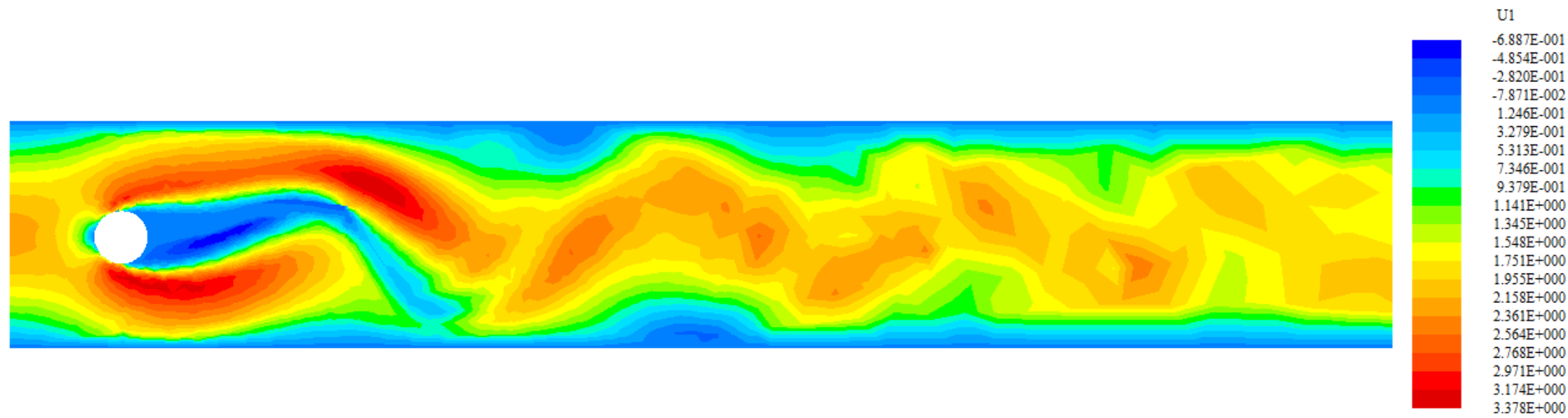}
    }
    \subfigure[\scriptsize $(\hat{U}, \mu_s) = (1.54, 0.46 \times 10^6)$: FOM $u_1$]{
        \includegraphics[width = 0.4\textwidth]{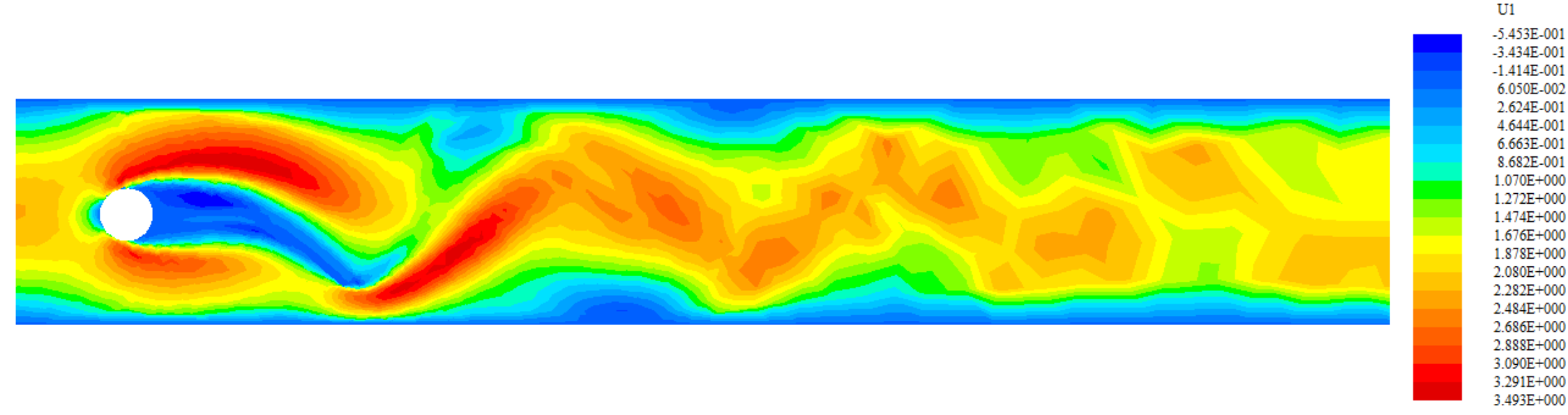}
    }
    \subfigure[\scriptsize $(\hat{U}, \mu_s) = (1.54, 0.46 \times 10^6)$: ROM $u_1$]{
        \includegraphics[width = 0.4\textwidth]{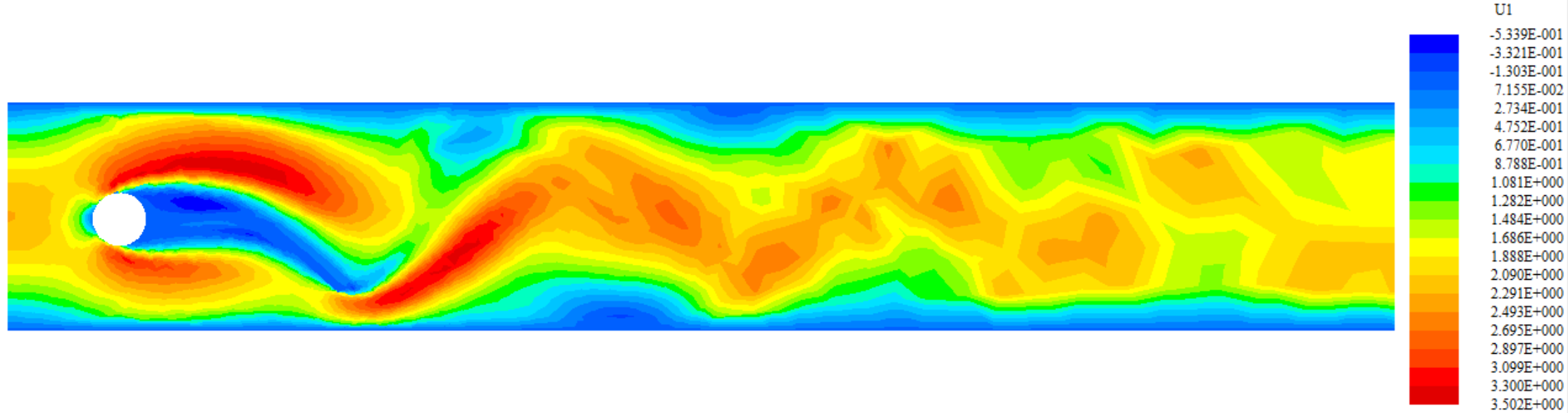}
    }
    \subfigure[\scriptsize $(\hat{U}, \mu_s) = (1.6, 0.54 \times 10^6)$: FOM $u_1$]{
        \includegraphics[width = 0.4\textwidth]{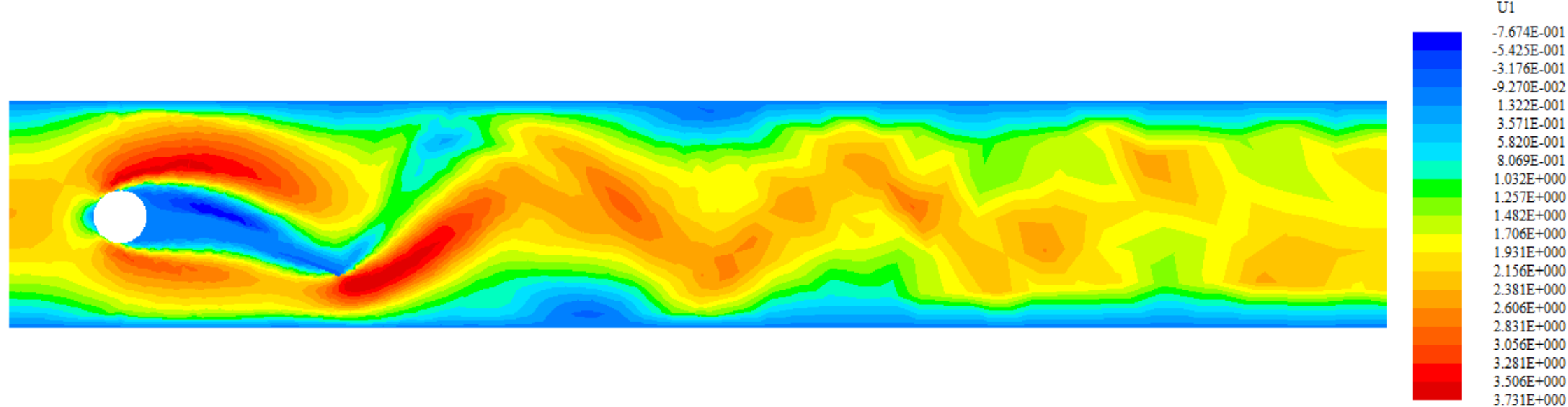}
    }
    \subfigure[\scriptsize $(\hat{U}, \mu_s) = (1.6, 0.54 \times 10^6)$: ROM $u_1$]{
        \includegraphics[width = 0.4\textwidth]{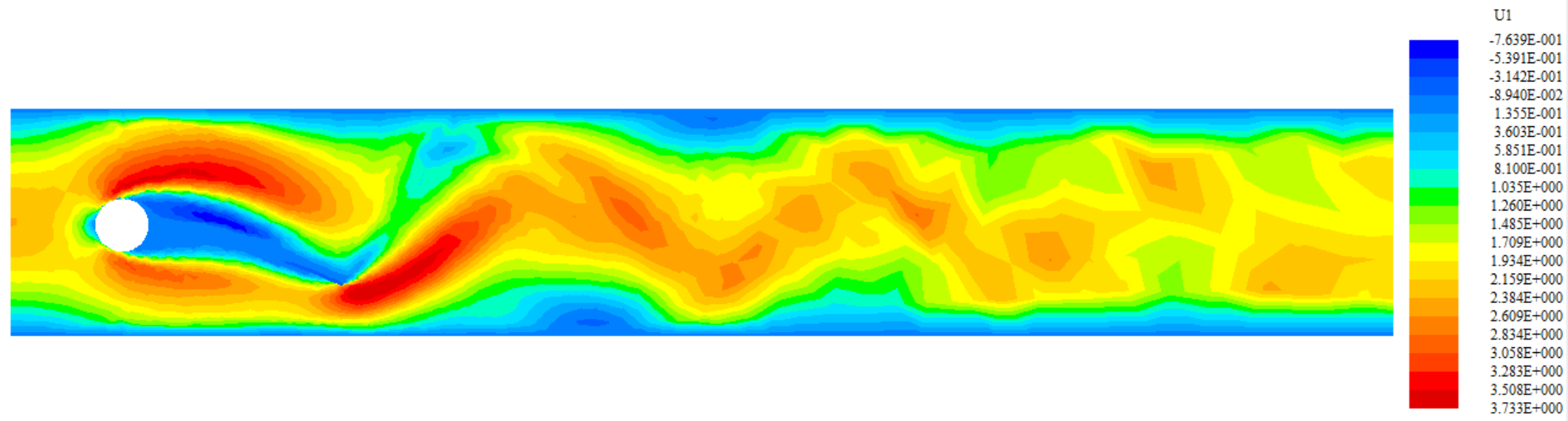}
    }
    \caption{Spatial distribution of horizontal velocity component $u_1$ at $T=15$s for training parameter combinations.}
    \label{spatial_u1_train}
\end{figure}

\begin{figure}[H]
    \centering
    \subfigure[\scriptsize $(\hat{U}, \mu_s) = (1.4, 0.4 \times 10^6)$: FOM $u_2$]{
        \includegraphics[width = 0.4\textwidth]{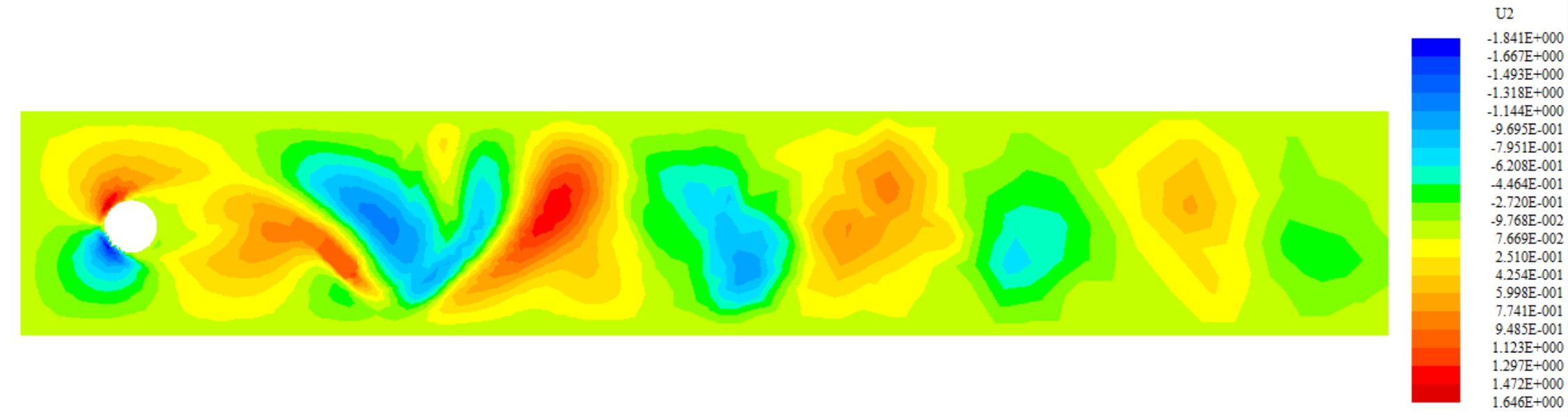}
    }
    \subfigure[\scriptsize $(\hat{U}, \mu_s) = (1.4, 0.4 \times 10^6)$: ROM $u_2$]{
        \includegraphics[width = 0.4\textwidth]{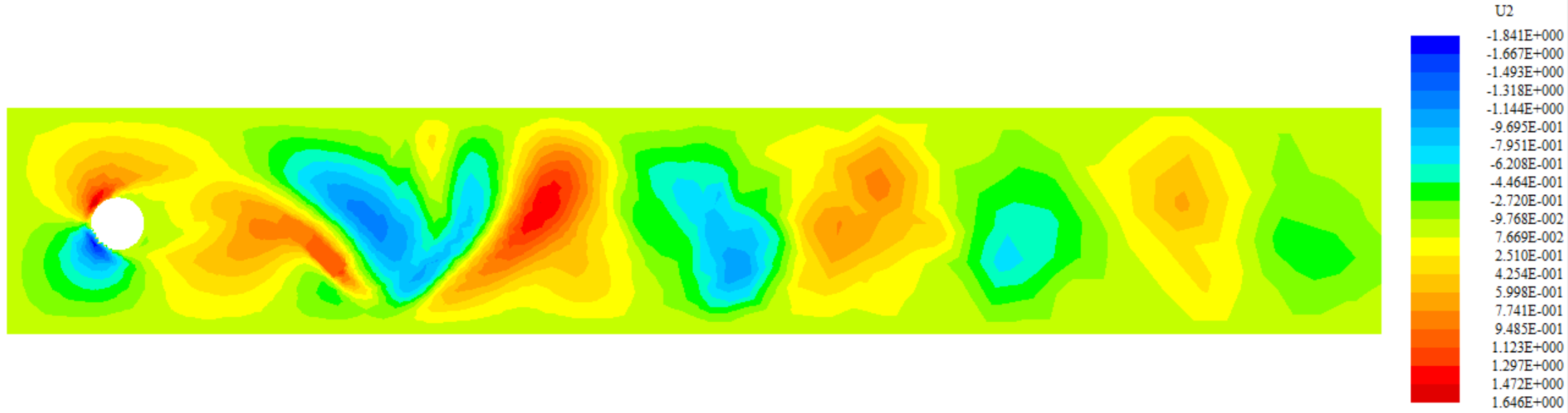}
    }
    \subfigure[\scriptsize $(\hat{U}, \mu_s) = (1.46, 0.5 \times 10^6)$: FOM $u_2$]{
        \includegraphics[width = 0.4\textwidth]{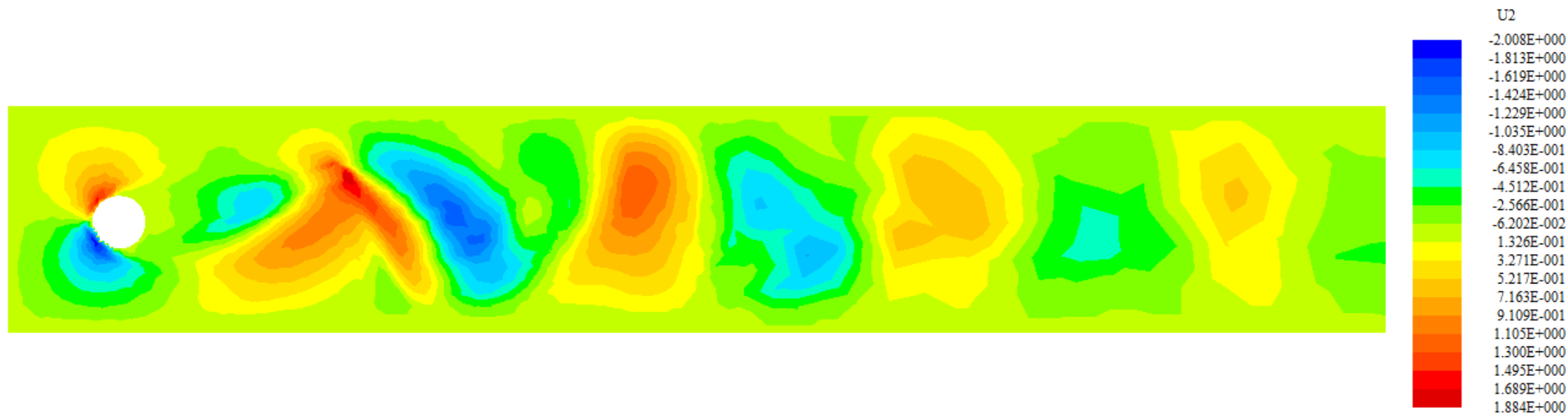}
    }
    \subfigure[\scriptsize $(\hat{U}, \mu_s) = (1.46, 0.5 \times 10^6)$: ROM $u_2$]{
        \includegraphics[width = 0.4\textwidth]{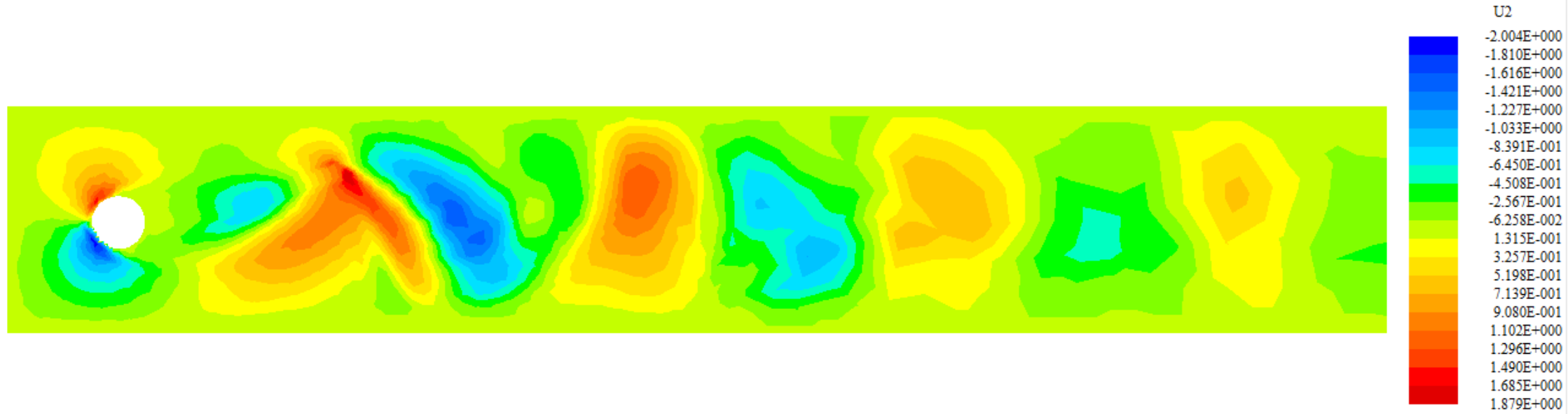}
    }
    \subfigure[\scriptsize $(\hat{U}, \mu_s) = (1.50, 0.6 \times 10^6)$: FOM $u_2$]{
        \includegraphics[width = 0.4\textwidth]{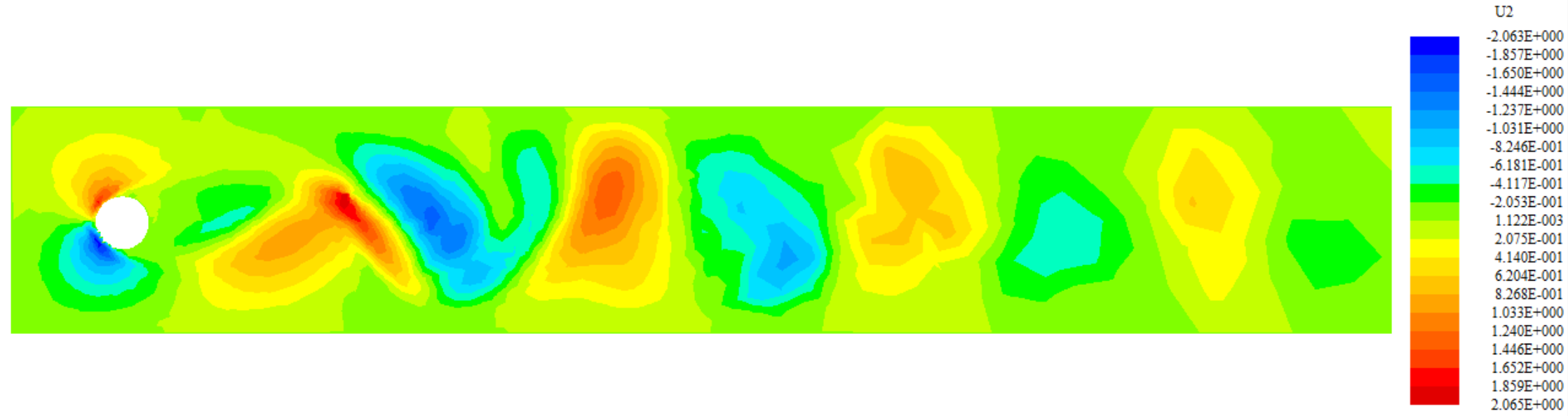}
    }
    \subfigure[\scriptsize $(\hat{U}, \mu_s) = (1.50, 0.6 \times 10^6)$: ROM $u_2$]{
        \includegraphics[width = 0.4\textwidth]{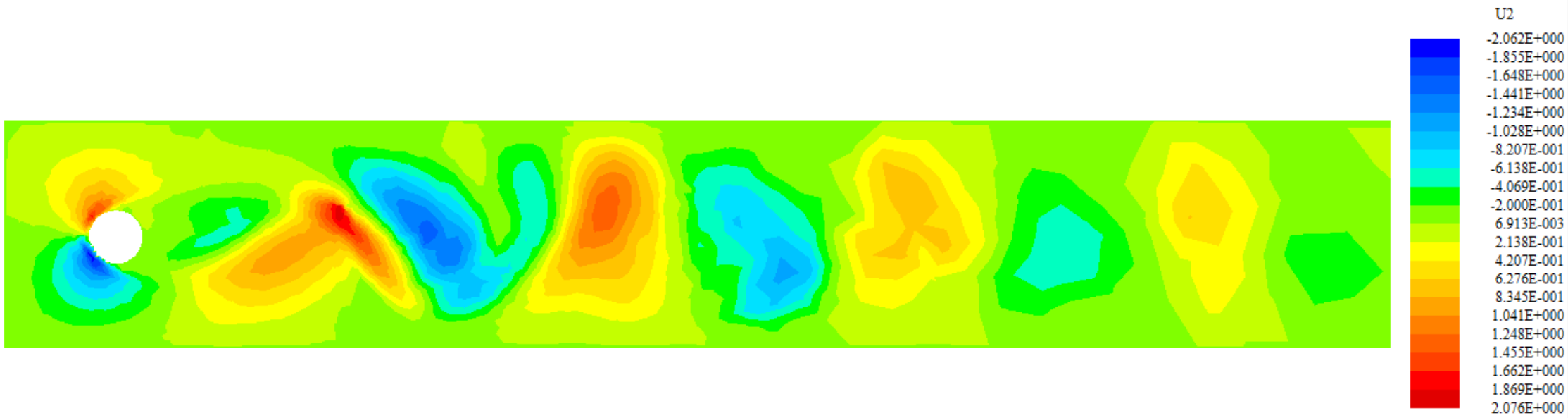}
    }
    \subfigure[\scriptsize $(\hat{U}, \mu_s) = (1.54, 0.46 \times 10^6)$: FOM $u_2$]{
        \includegraphics[width = 0.4\textwidth]{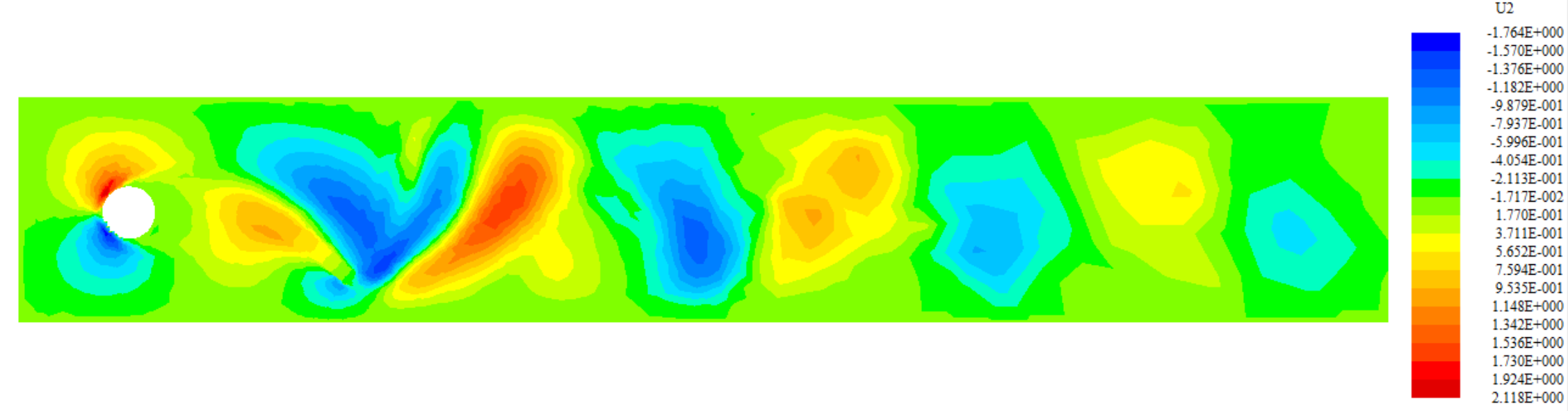}
    }
    \subfigure[\scriptsize $(\hat{U}, \mu_s) = (1.54, 0.46 \times 10^6)$: ROM $u_2$]{
        \includegraphics[width = 0.4\textwidth]{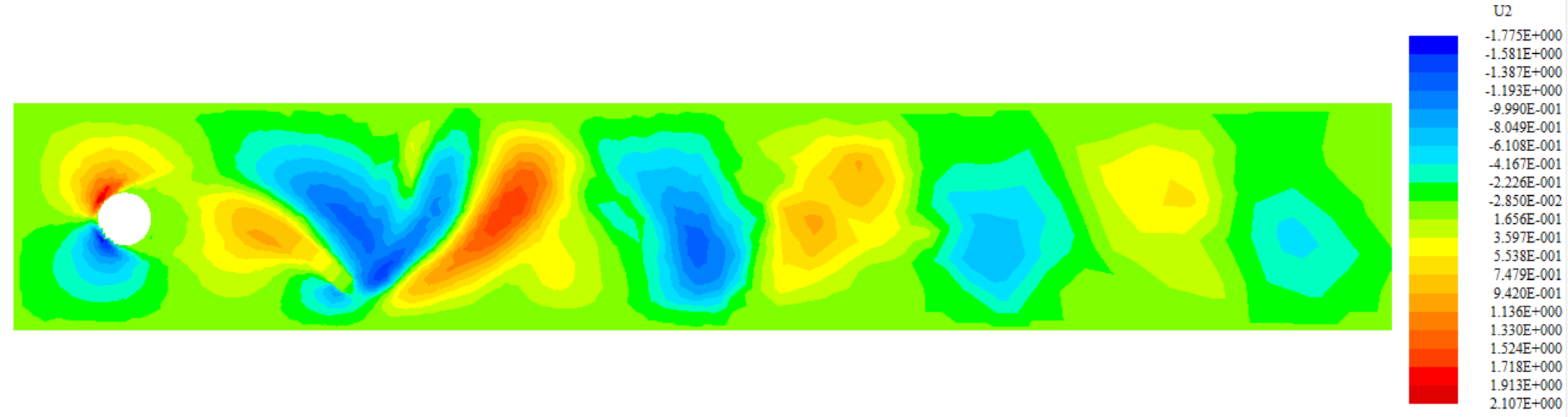}
    }
    \subfigure[\scriptsize $(\hat{U}, \mu_s) = (1.6, 0.54 \times 10^6)$: FOM $u_2$]{
        \includegraphics[width = 0.4\textwidth]{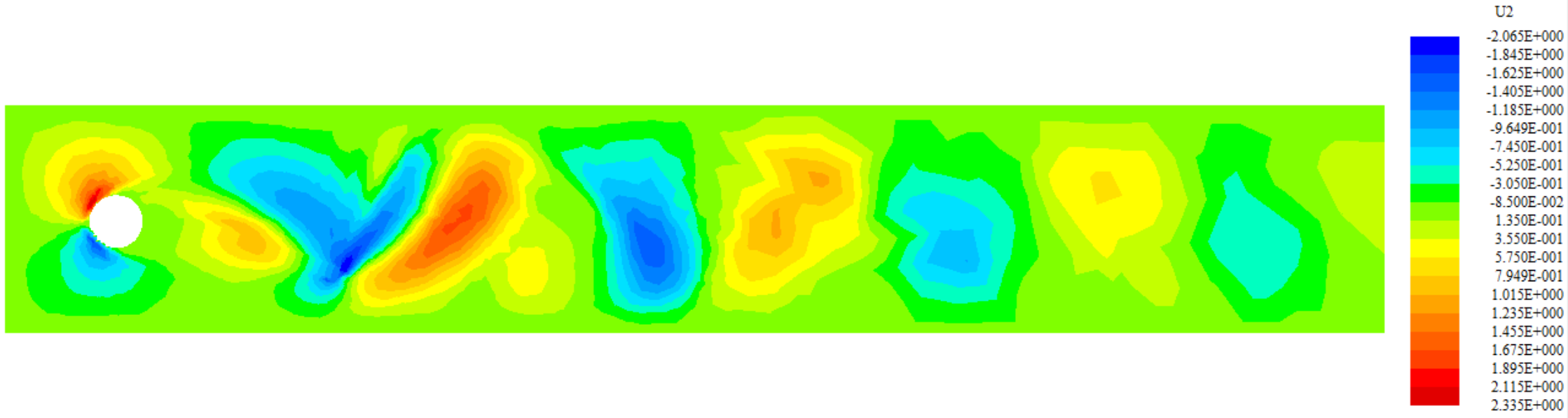}
    }
    \subfigure[\scriptsize $(\hat{U}, \mu_s) = (1.6, 0.54 \times 10^6)$: ROM $u_2$]{
        \includegraphics[width = 0.4\textwidth]{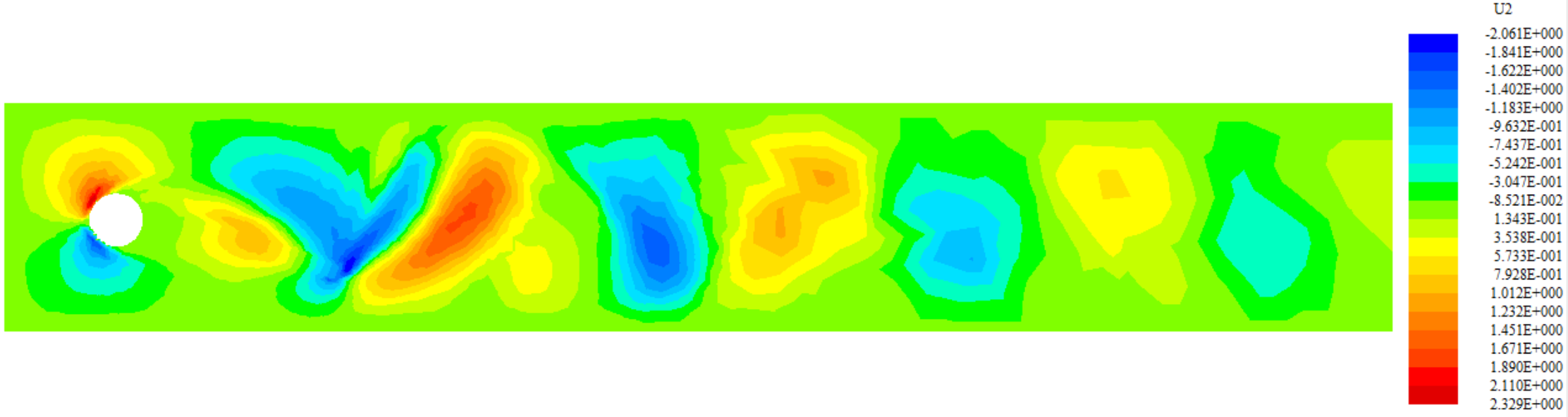}
    }
    \caption{Spatial distribution of vertical velocity component $u_2$ at $T=15$s for training parameter combinations.}
    \label{spatial_u2_train}
\end{figure}

\begin{figure}[H]
    \centering
    \subfigure[\scriptsize $(\hat{U}, \mu_s) = (1.4, 0.4 \times 10^6)$: FOM $p$]{
        \includegraphics[width = 0.4\textwidth]{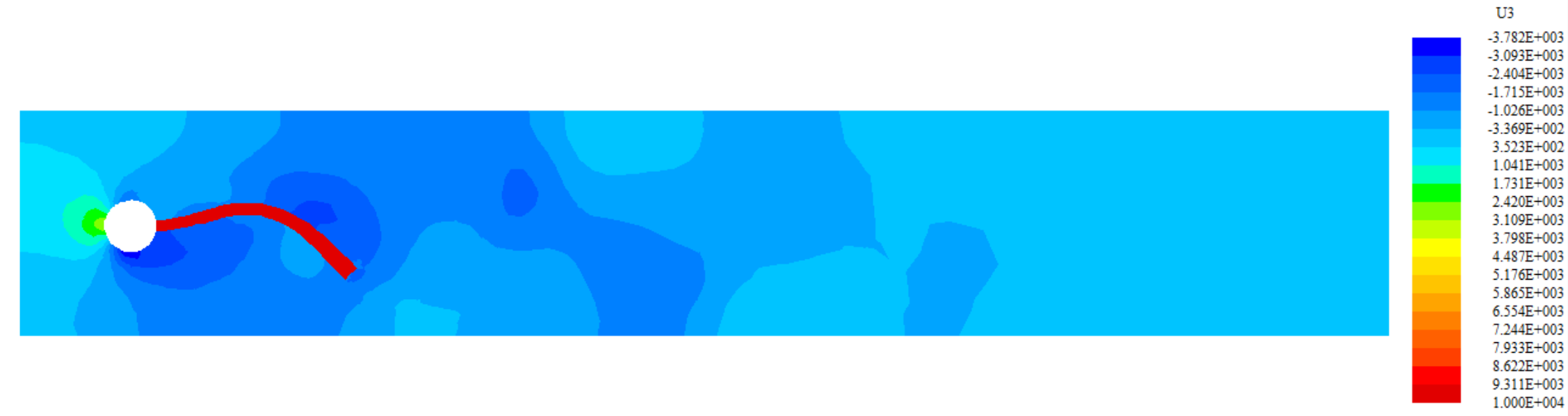}
    }
    \subfigure[\scriptsize $(\hat{U}, \mu_s) = (1.4, 0.4 \times 10^6)$: ROM $p$]{
        \includegraphics[width = 0.4\textwidth]{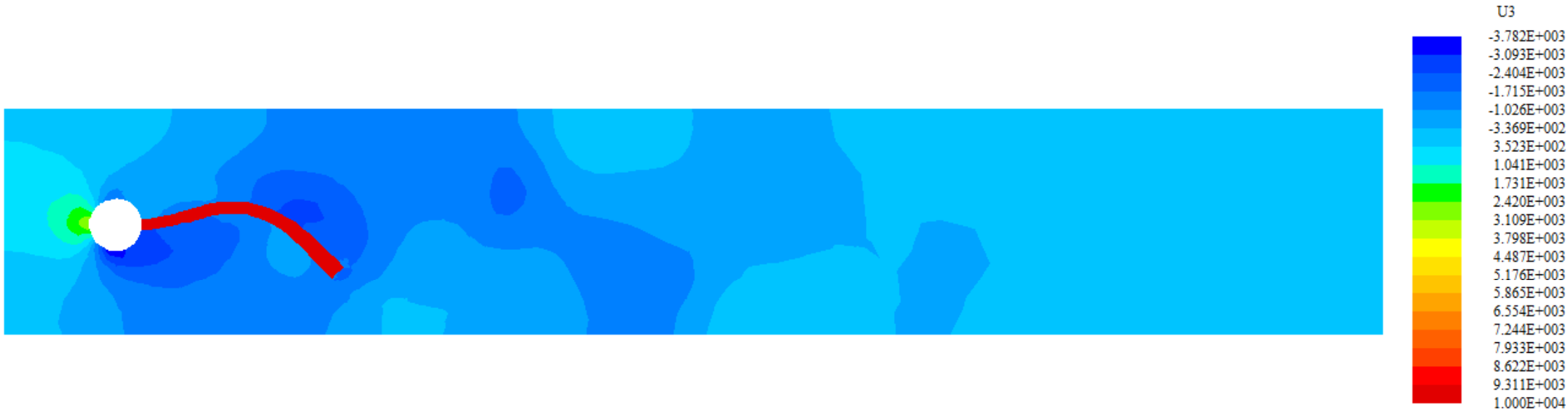}
    }
    \subfigure[\scriptsize $(\hat{U}, \mu_s) = (1.46, 0.5 \times 10^6)$: FOM $p$]{
        \includegraphics[width = 0.4\textwidth]{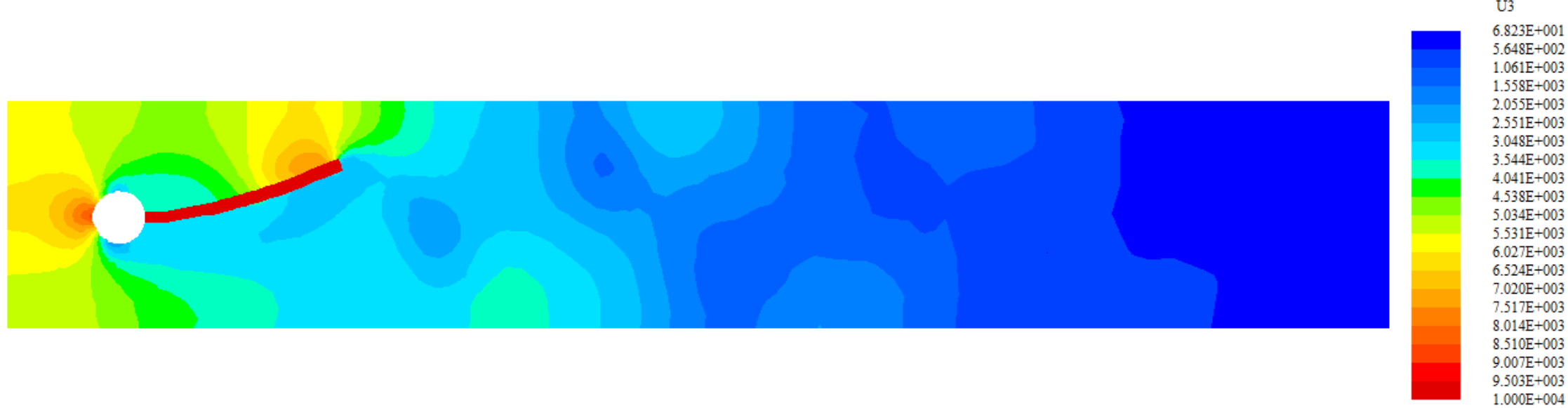}
    }
    \subfigure[\scriptsize $(\hat{U}, \mu_s) = (1.46, 0.5 \times 10^6)$: ROM $p$]{
        \includegraphics[width = 0.4\textwidth]{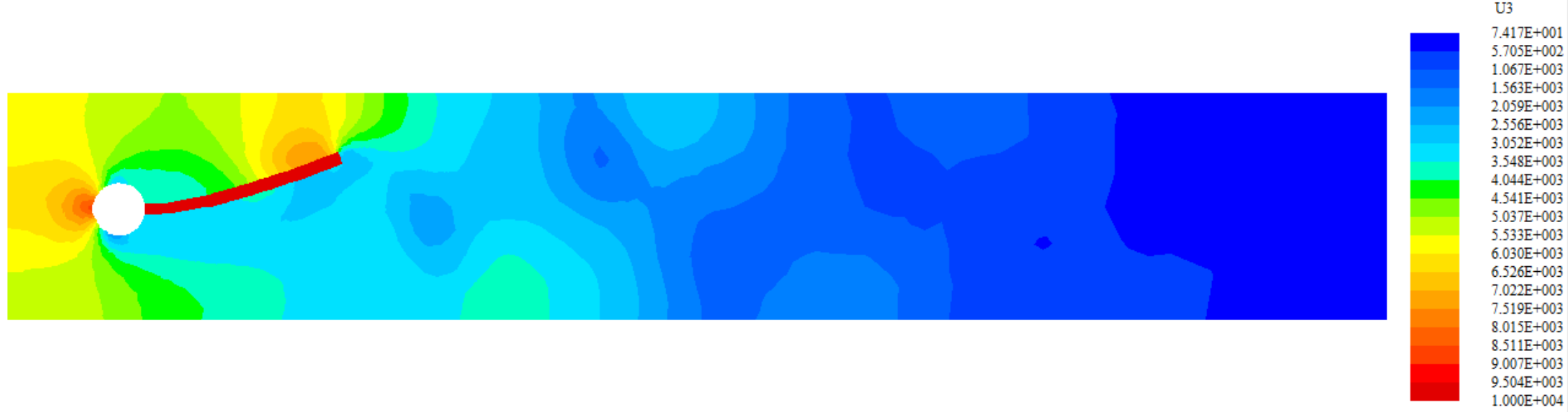}
    }
    \subfigure[\scriptsize $(\hat{U}, \mu_s) = (1.50, 0.6 \times 10^6)$: FOM $p$]{
        \includegraphics[width = 0.4\textwidth]{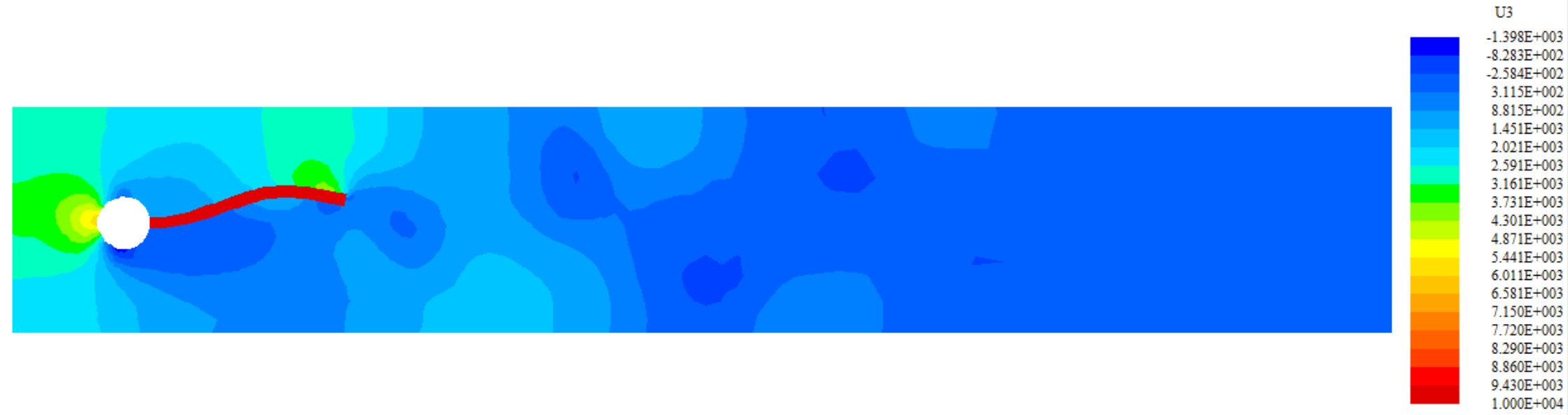}
    }
    \subfigure[\scriptsize $(\hat{U}, \mu_s) = (1.50, 0.6 \times 10^6)$: ROM $p$]{
        \includegraphics[width = 0.4\textwidth]{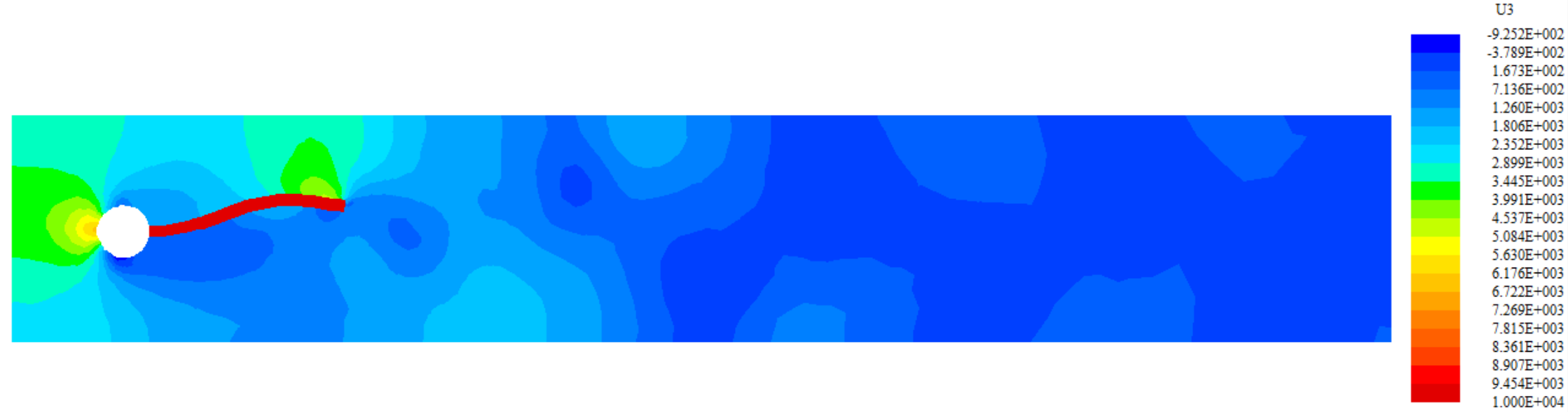}
    }
    \subfigure[\scriptsize $(\hat{U}, \mu_s) = (1.54, 0.46 \times 10^6)$: FOM $p$]{
        \includegraphics[width = 0.4\textwidth]{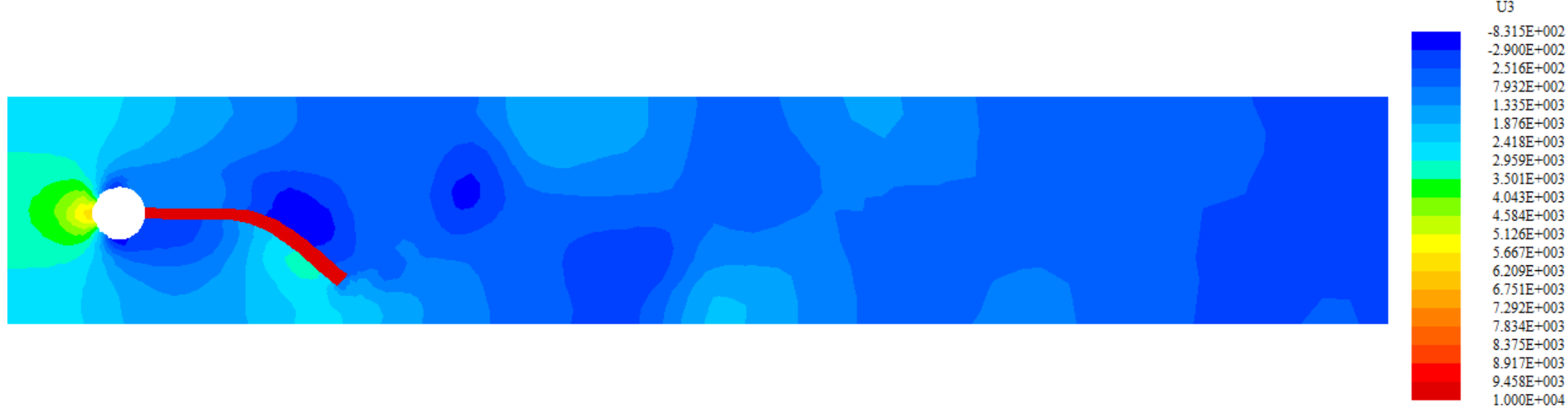}
    }
    \subfigure[\scriptsize $(\hat{U}, \mu_s) = (1.54, 0.46 \times 10^6)$: ROM $p$]{
        \includegraphics[width = 0.4\textwidth]{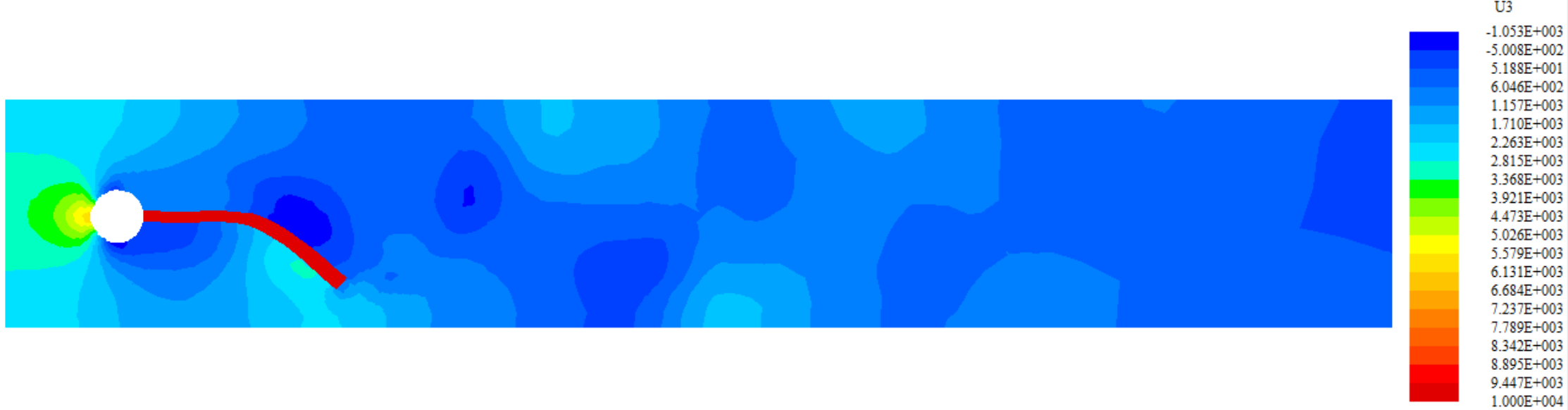}
    }
    \subfigure[\scriptsize $(\hat{U}, \mu_s) = (1.6, 0.54 \times 10^6)$: FOM $p$]{
        \includegraphics[width = 0.4\textwidth]{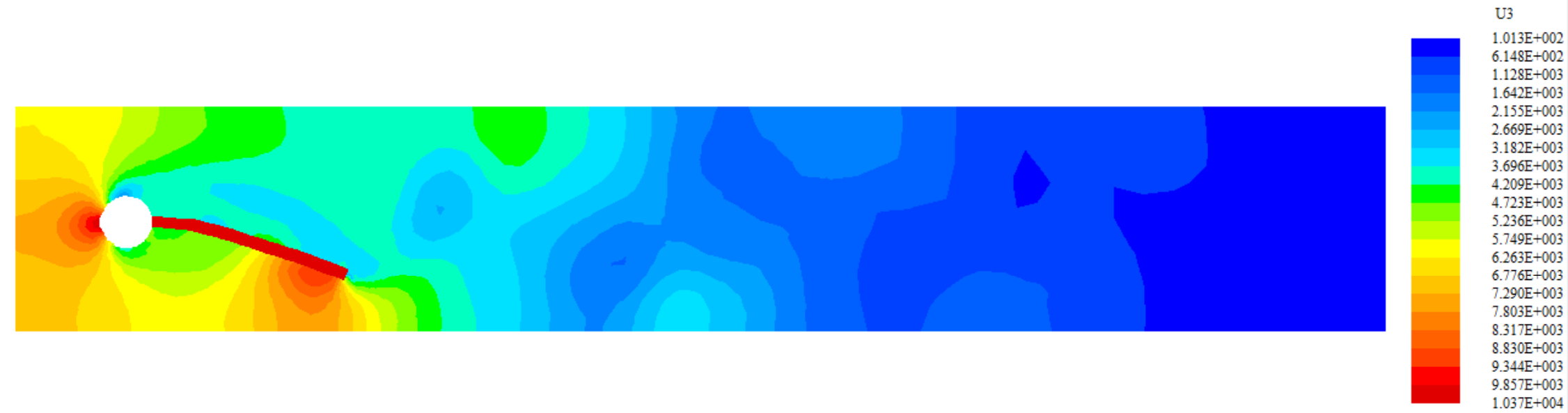}
    }
    \subfigure[\scriptsize $(\hat{U}, \mu_s) = (1.6, 0.54 \times 10^6)$: ROM $p$]{
        \includegraphics[width = 0.4\textwidth]{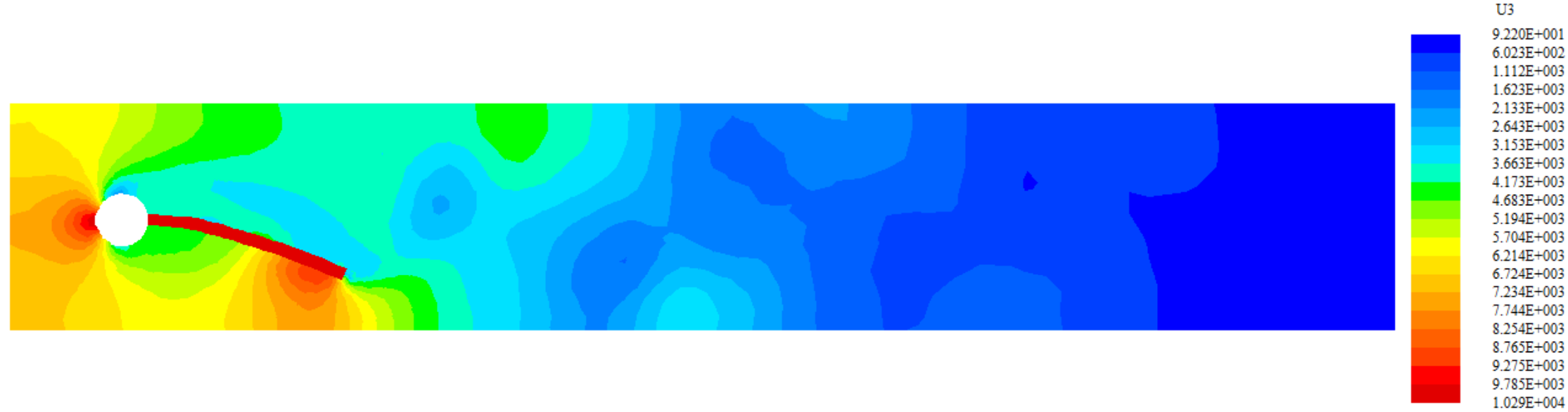}
    }
    \caption{Spatial distribution of pressure $p$ at $T=15$s for training parameter combinations.}
    \label{spatial_p_train}
\end{figure}

\begin{figure}[H]
    \centering
    \subfigure[\scriptsize $(\hat{U}, \mu_s) = (1.41, 0.45 \times 10^6)$: FOM $u_1$]{
        \includegraphics[width = 0.4\textwidth]{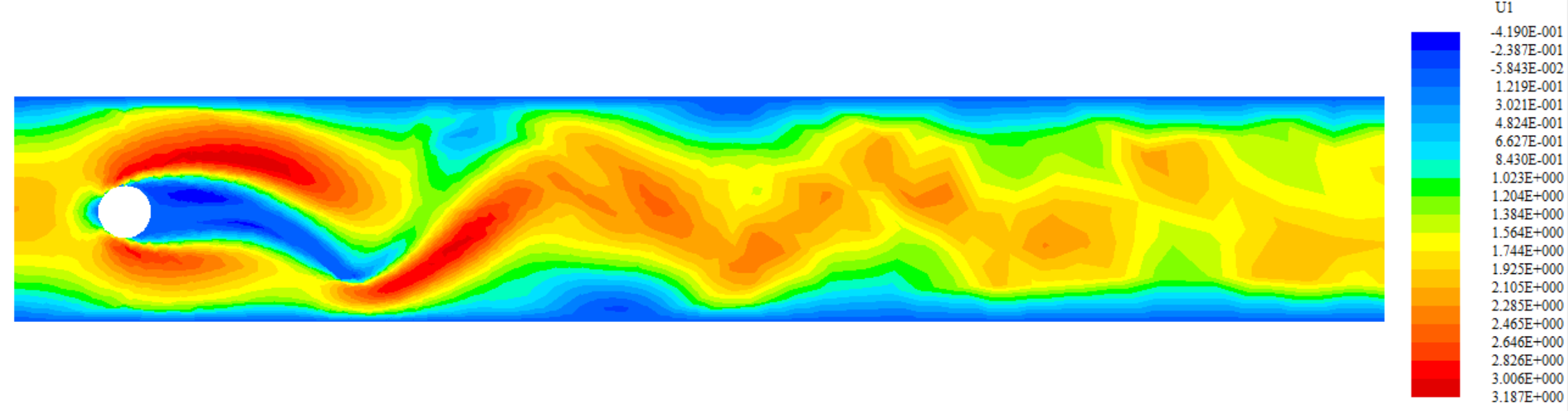}
    }
    \subfigure[\scriptsize $(\hat{U}, \mu_s) = (1.41, 0.45 \times 10^6)$: ROM $u_1$]{
        \includegraphics[width = 0.4\textwidth]{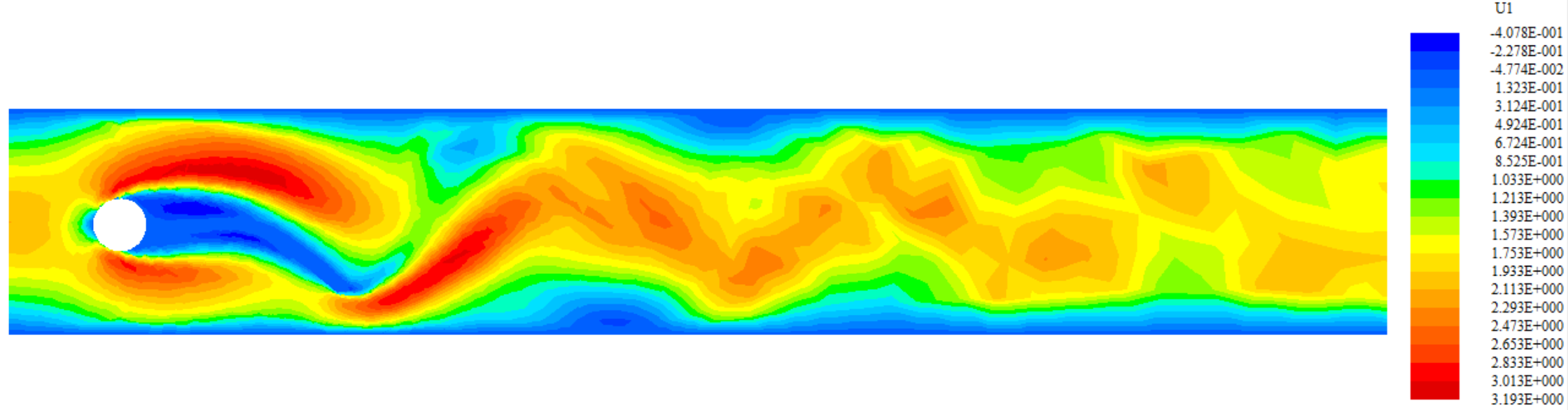}
    }
    \subfigure[\scriptsize $(\hat{U}, \mu_s) = (1.47, 0.55 \times 10^6)$: FOM $u_1$]{
        \includegraphics[width = 0.4\textwidth]{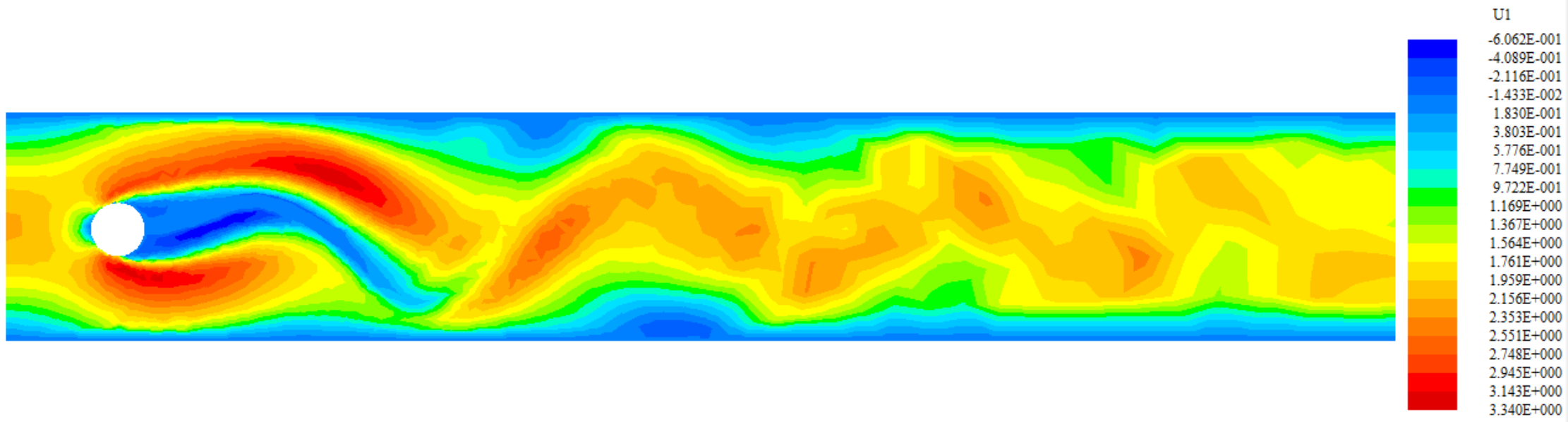}
    }
    \subfigure[\scriptsize $(\hat{U}, \mu_s) = (1.47, 0.55 \times 10^6)$: ROM $u_1$]{
        \includegraphics[width = 0.4\textwidth]{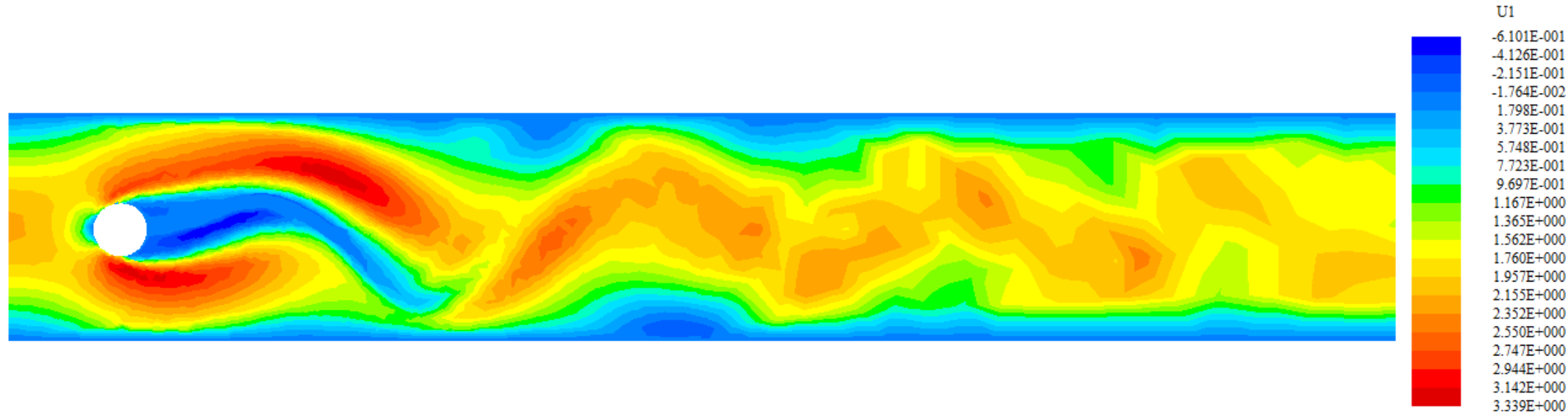}
    }
    \subfigure[\scriptsize $(\hat{U}, \mu_s) = (1.53, 0.51 \times 10^6)$: FOM $u_1$]{
        \includegraphics[width = 0.4\textwidth]{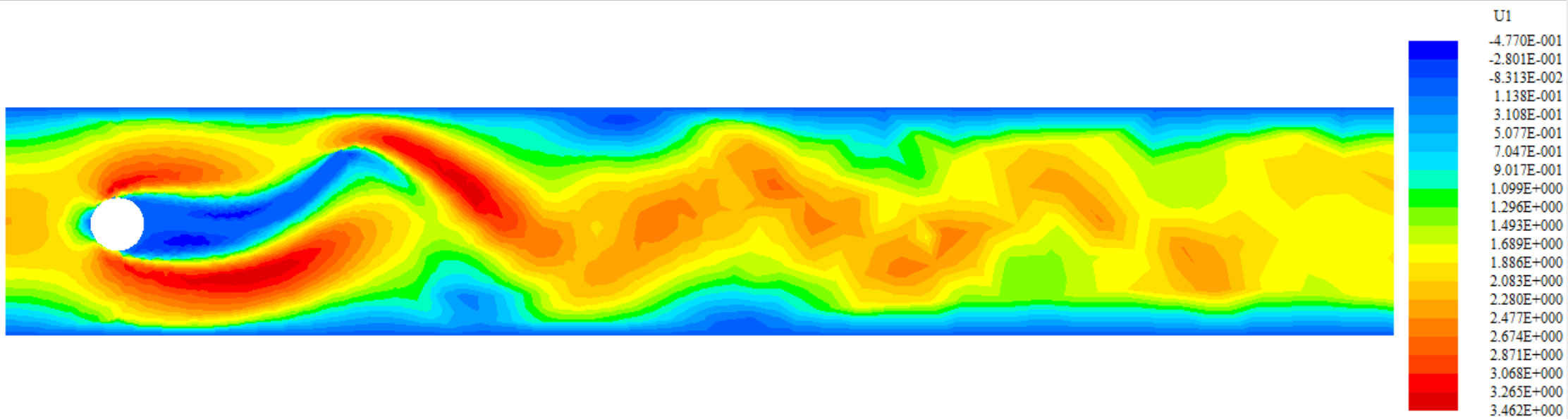}
    }
    \subfigure[\scriptsize $(\hat{U}, \mu_s) = (1.53, 0.51 \times 10^6)$: ROM $u_1$]{
        \includegraphics[width = 0.4\textwidth]{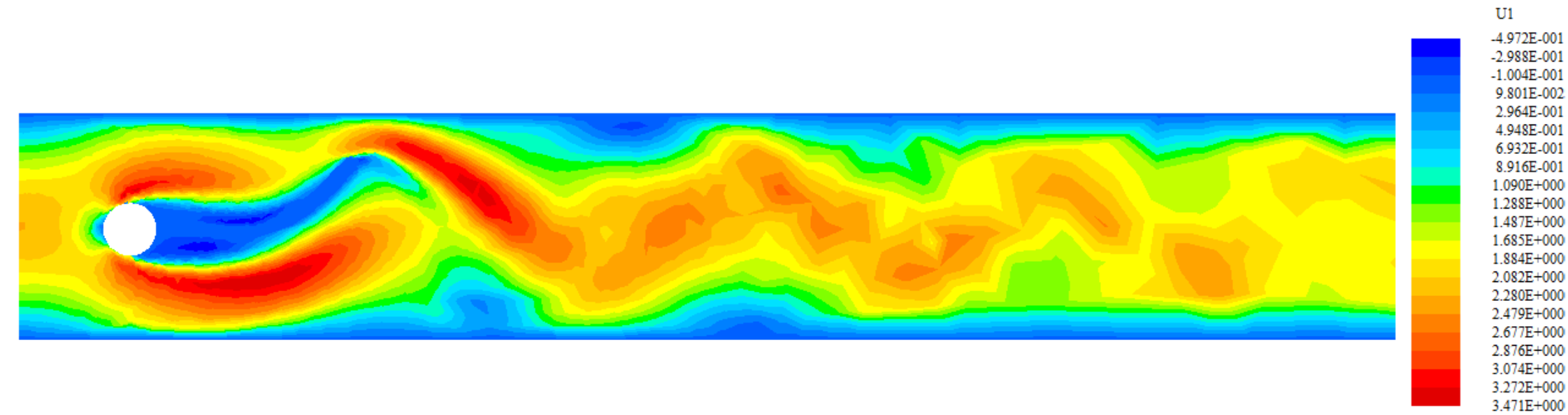}
    }
    \subfigure[\scriptsize $(\hat{U}, \mu_s) = (1.55, 0.49 \times 10^6)$: FOM $u_1$]{
        \includegraphics[width = 0.4\textwidth]{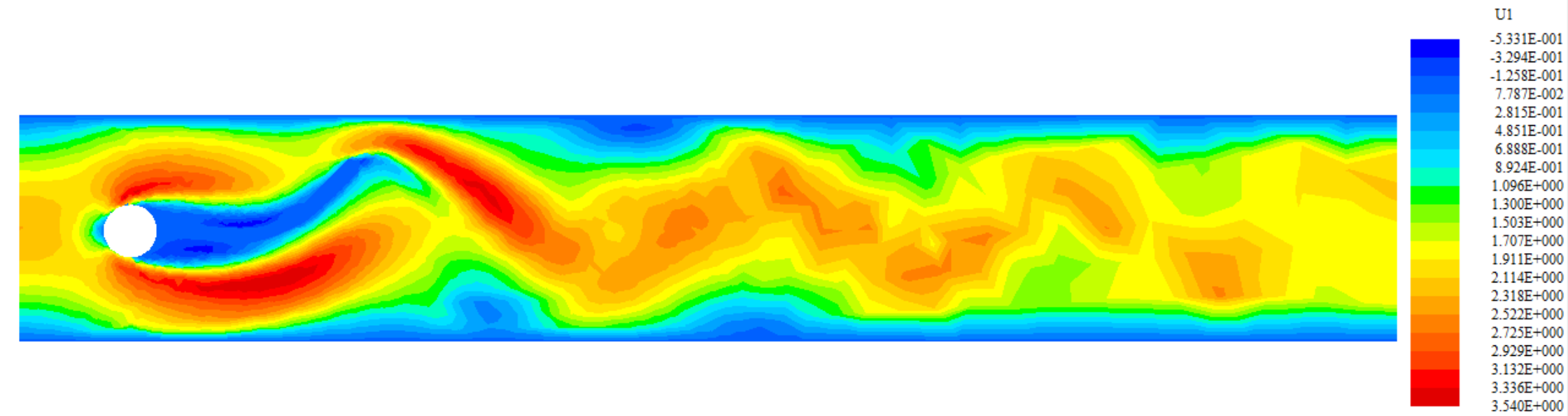}
    }
    \subfigure[\scriptsize $(\hat{U}, \mu_s) = (1.55, 0.49 \times 10^6)$: ROM $u_1$]{
        \includegraphics[width = 0.4\textwidth]{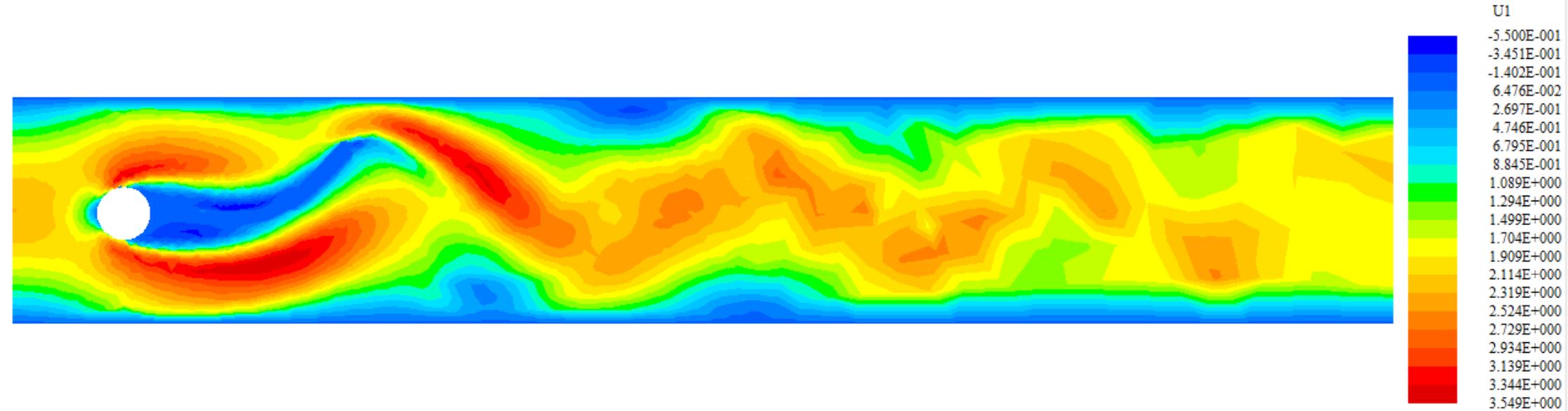}
    }
    \subfigure[\scriptsize $(\hat{U}, \mu_s) = (1.59, 0.59 \times 10^6)$: FOM $u_1$]{
        \includegraphics[width = 0.4\textwidth]{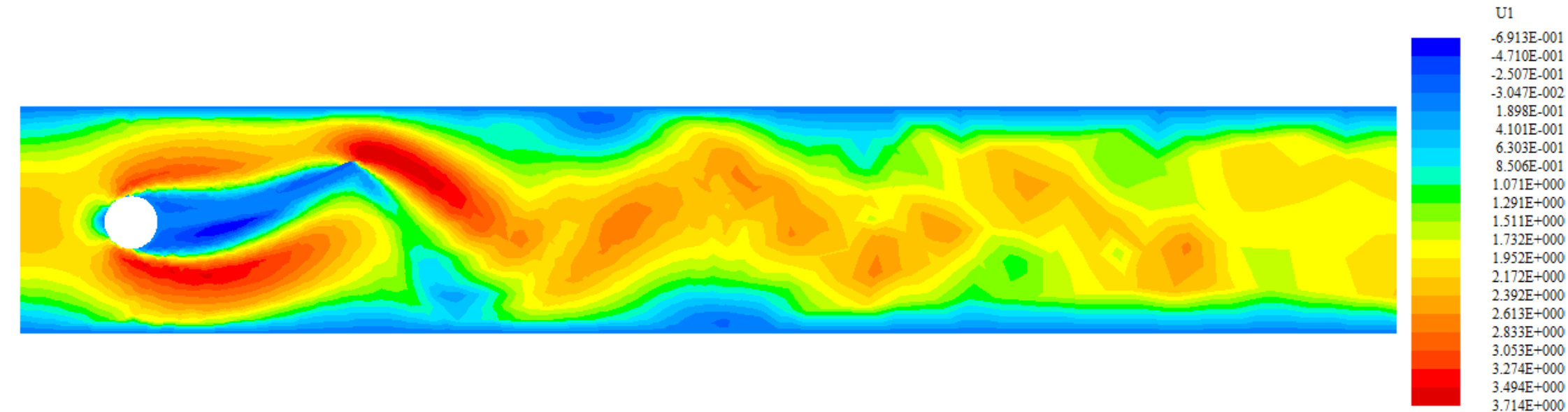}
    }
    \subfigure[\scriptsize $(\hat{U}, \mu_s) = (1.59, 0.59 \times 10^6)$: ROM $u_1$]{
        \includegraphics[width = 0.4\textwidth]{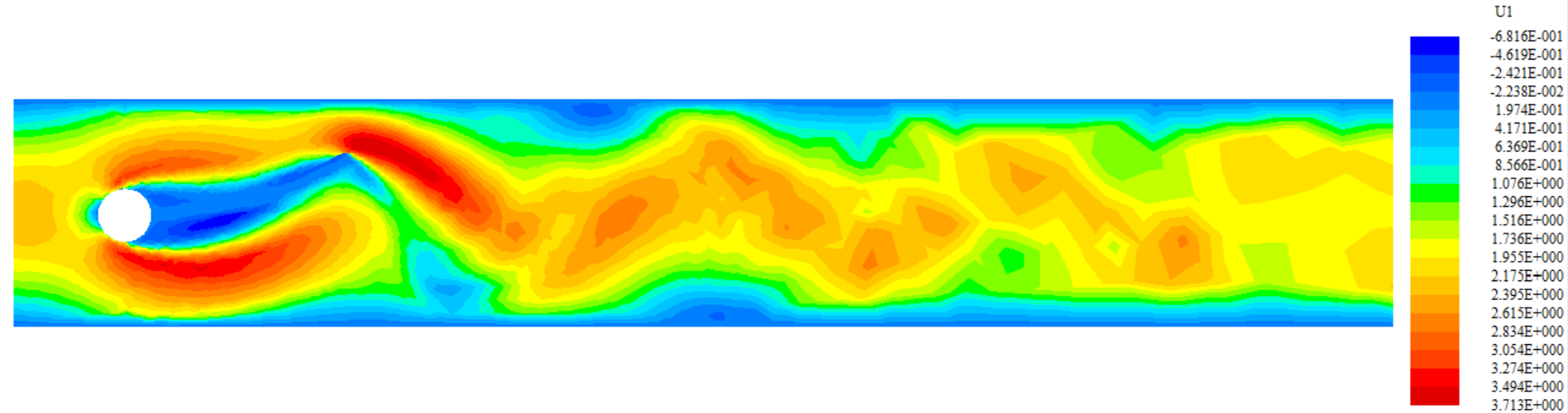}
    }
    \caption{Spatial distribution of horizontal velocity component $u_1$ at $T=15$s for test parameter combinations.}
    \label{spatial_u1_test}
\end{figure}

\begin{figure}[H]
    \centering
    \subfigure[\scriptsize $(\hat{U}, \mu_s) = (1.41, 0.45 \times 10^6)$: FOM $u_2$]{
        \includegraphics[width = 0.4\textwidth]{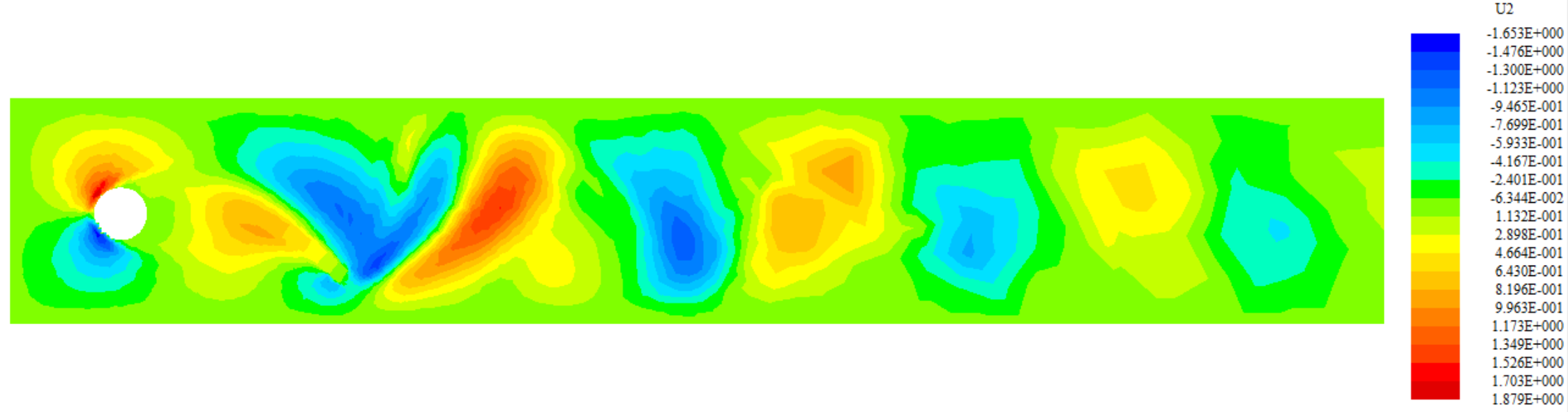}
    }
    \subfigure[\scriptsize $(\hat{U}, \mu_s) = (1.41, 0.45 \times 10^6)$: ROM $u_2$]{
        \includegraphics[width = 0.4\textwidth]{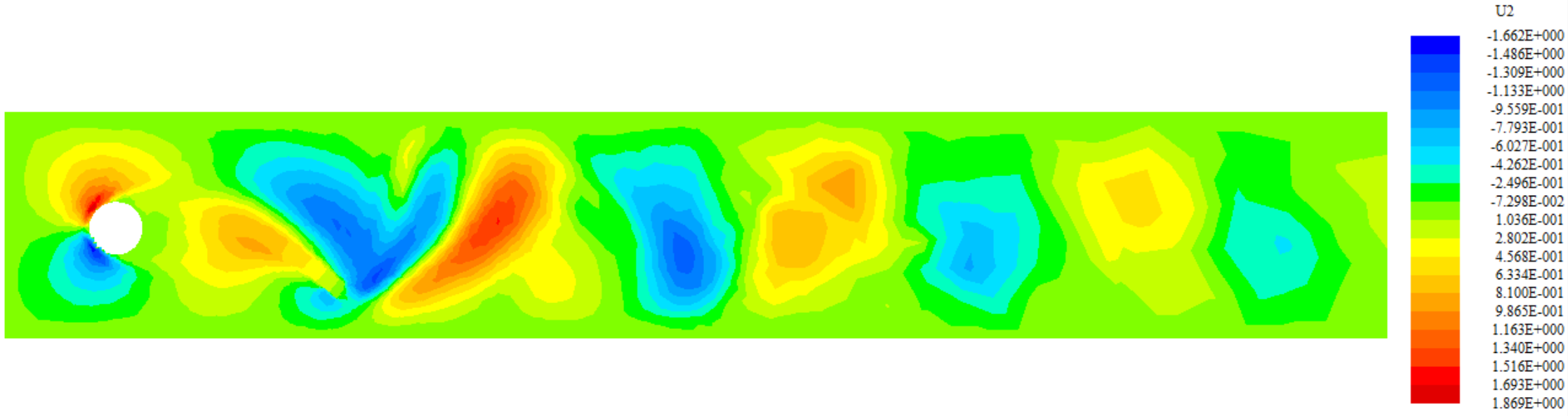}
    }
    \subfigure[\scriptsize $(\hat{U}, \mu_s) = (1.47, 0.55 \times 10^6)$: FOM $u_2$]{
        \includegraphics[width = 0.4\textwidth]{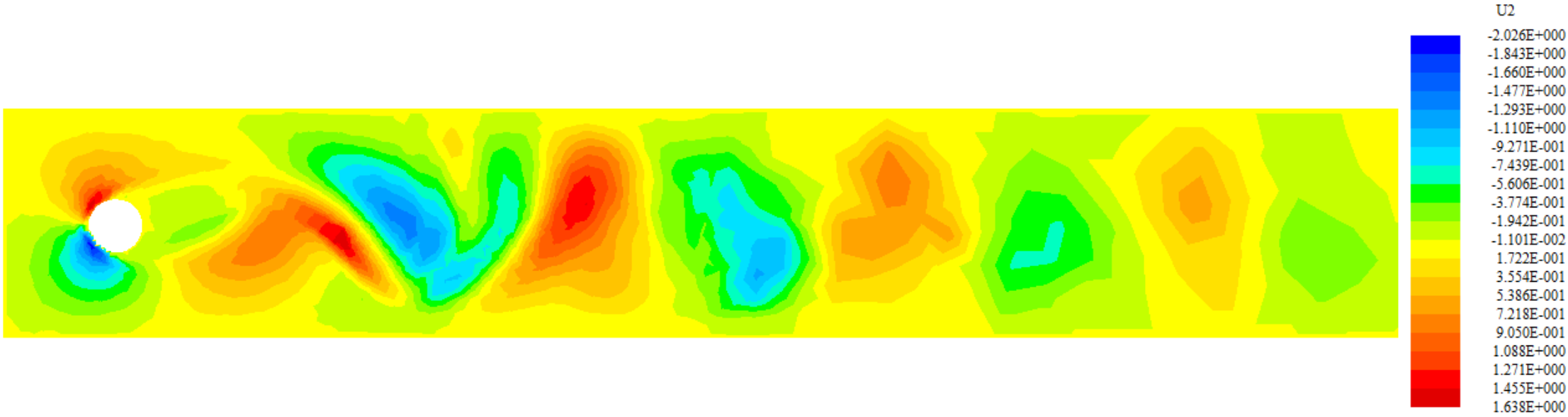}
    }
    \subfigure[\scriptsize $(\hat{U}, \mu_s) = (1.47, 0.55 \times 10^6)$: ROM $u_2$]{
        \includegraphics[width = 0.4\textwidth]{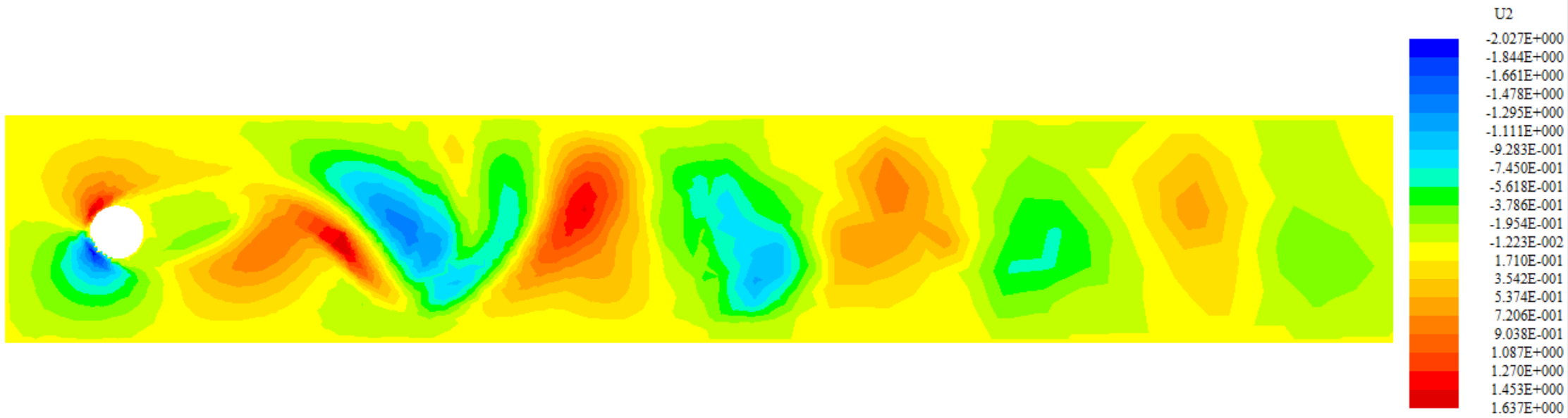}
    }
    \subfigure[\scriptsize $(\hat{U}, \mu_s) = (1.53, 0.51 \times 10^6)$: FOM $u_2$]{
        \includegraphics[width = 0.4\textwidth]{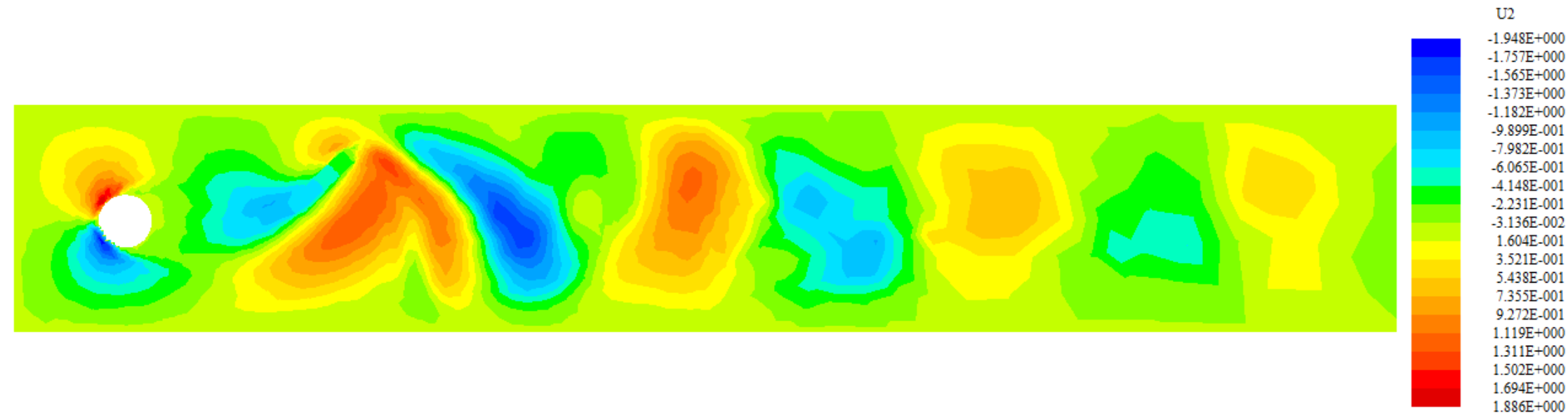}
    }
    \subfigure[\scriptsize $(\hat{U}, \mu_s) = (1.53, 0.51 \times 10^6)$: ROM $u_2$]{
        \includegraphics[width = 0.4\textwidth]{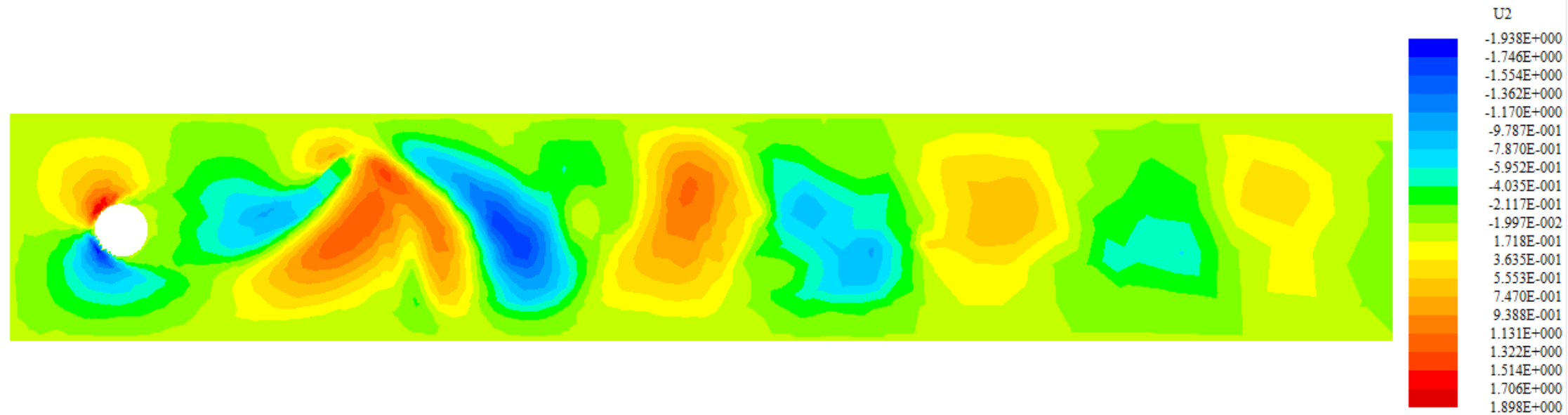}
    }
    \subfigure[\scriptsize $(\hat{U}, \mu_s) = (1.55, 0.49 \times 10^6)$: FOM $u_2$]{
        \includegraphics[width = 0.4\textwidth]{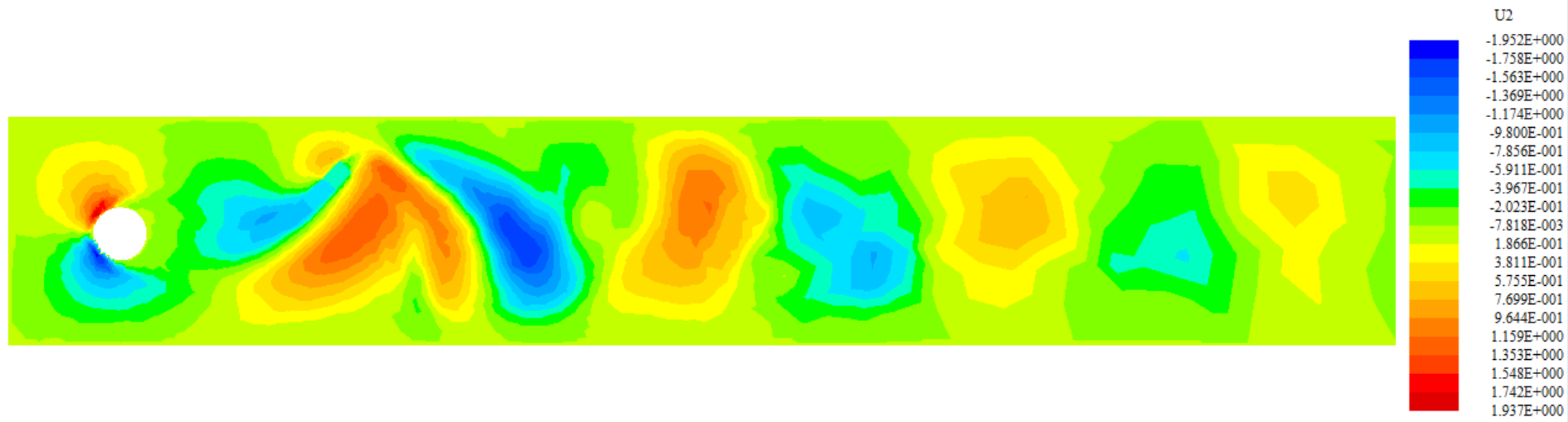}
    }
    \subfigure[\scriptsize $(\hat{U}, \mu_s) = (1.55, 0.49 \times 10^6)$: ROM $u_2$]{
        \includegraphics[width = 0.4\textwidth]{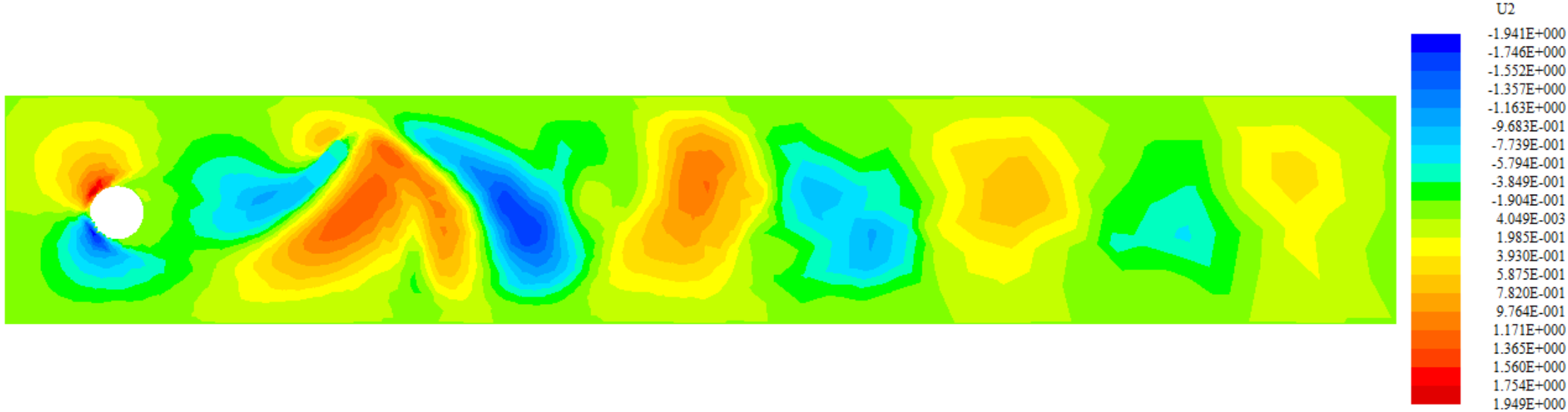}
    }
    \subfigure[\scriptsize $(\hat{U}, \mu_s) = (1.59, 0.59 \times 10^6)$: FOM $u_2$]{
        \includegraphics[width = 0.4\textwidth]{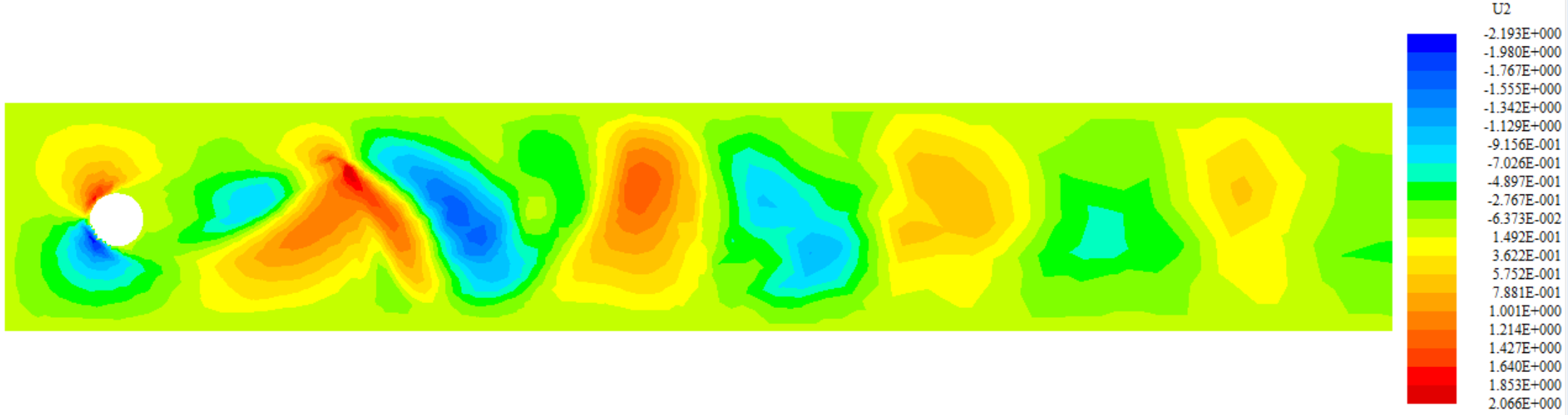}
    }
    \subfigure[\scriptsize $(\hat{U}, \mu_s) = (1.59, 0.59 \times 10^6)$: ROM $u_2$]{
        \includegraphics[width = 0.4\textwidth]{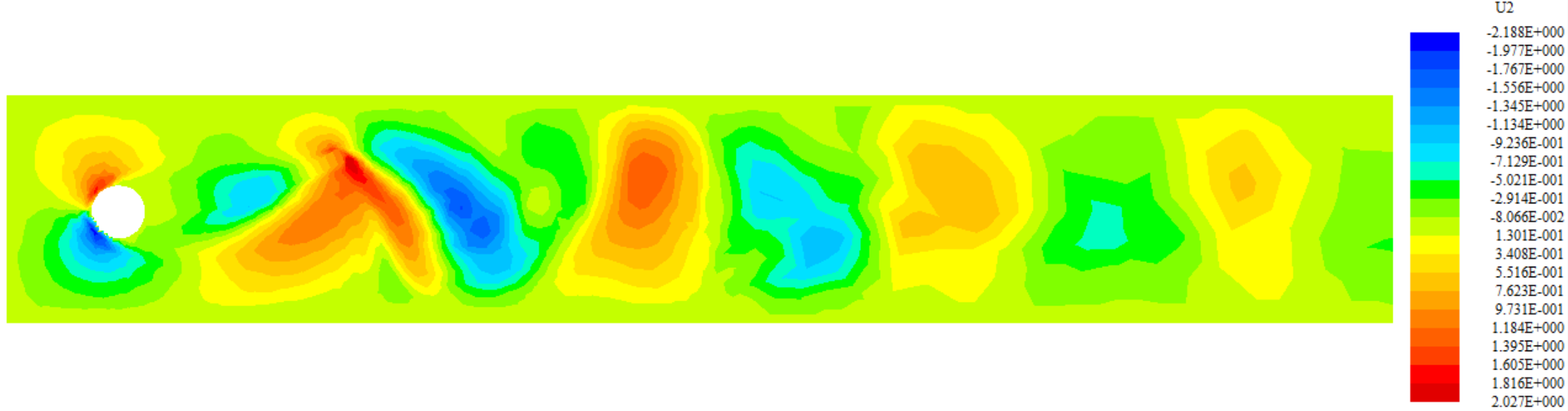}
    }
    \caption{Spatial distribution of vertical velocity component $u_2$ at $T=15$s for test parameter combinations.}
    \label{spatial_u2_test}
\end{figure}

\begin{figure}[H]
    \centering
    \subfigure[\scriptsize $(\hat{U}, \mu_s) = (1.41, 0.45 \times 10^6)$: FOM $p$]{
        \includegraphics[width = 0.4\textwidth]{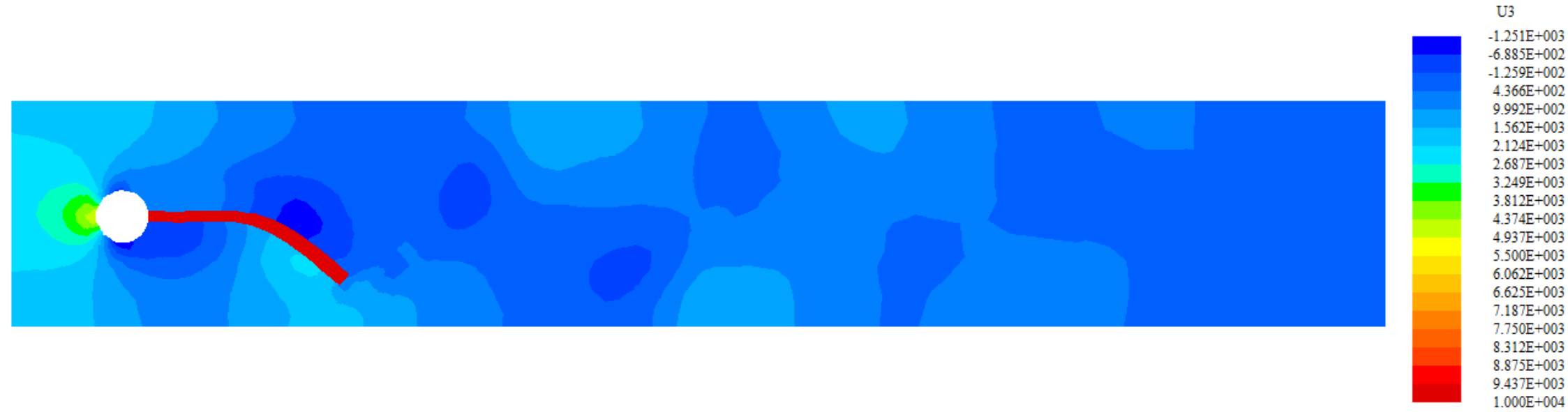}
    }
    \subfigure[\scriptsize $(\hat{U}, \mu_s) = (1.41, 0.45 \times 10^6)$: ROM $p$]{
        \includegraphics[width = 0.4\textwidth]{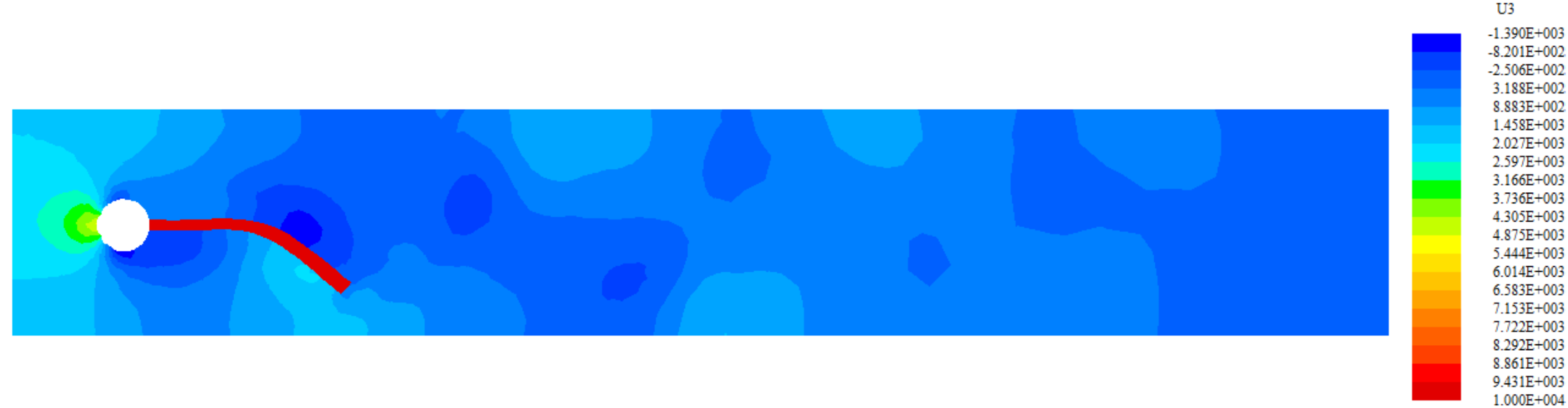}
    }
    \subfigure[\scriptsize $(\hat{U}, \mu_s) = (1.47, 0.55 \times 10^6)$: FOM $p$]{
        \includegraphics[width = 0.4\textwidth]{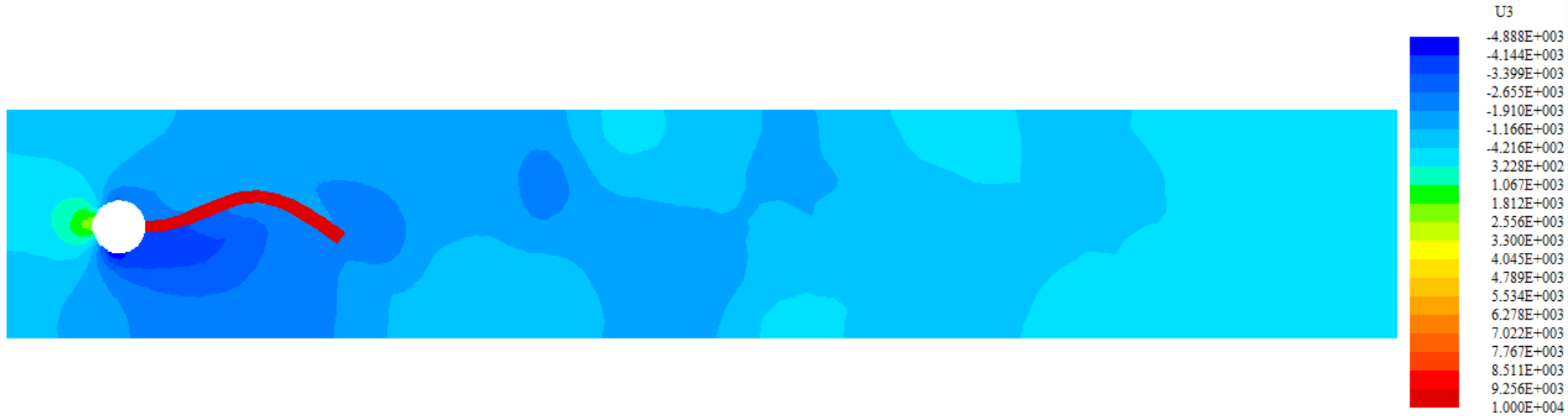}
    }
    \subfigure[\scriptsize $(\hat{U}, \mu_s) = (1.47, 0.55 \times 10^6)$: ROM $p$]{
        \includegraphics[width = 0.4\textwidth]{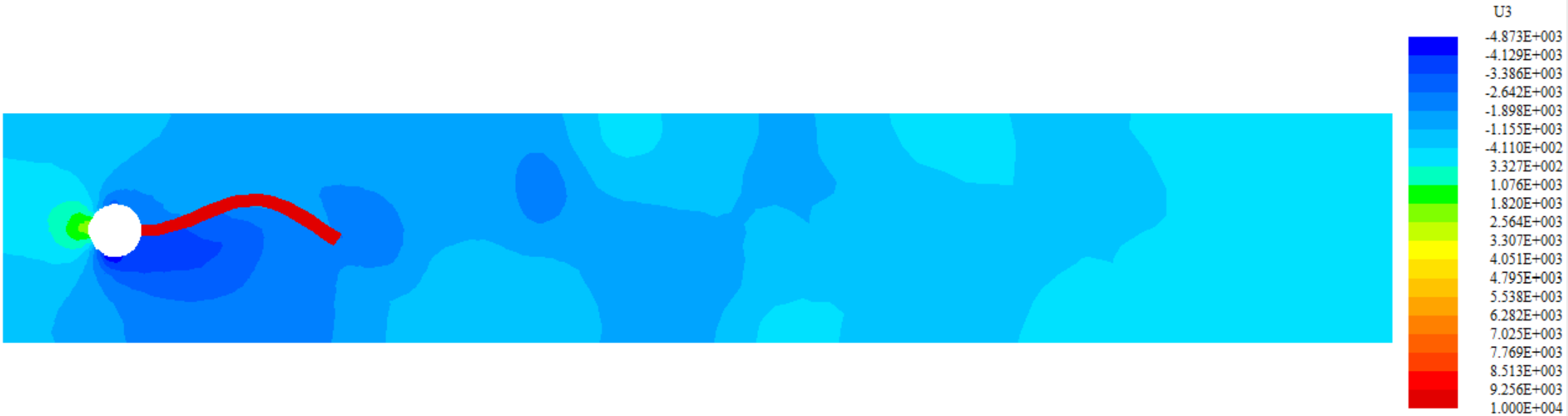}
    }
    \subfigure[\scriptsize $(\hat{U}, \mu_s) = (1.53, 0.51 \times 10^6)$: FOM $p$]{
        \includegraphics[width = 0.4\textwidth]{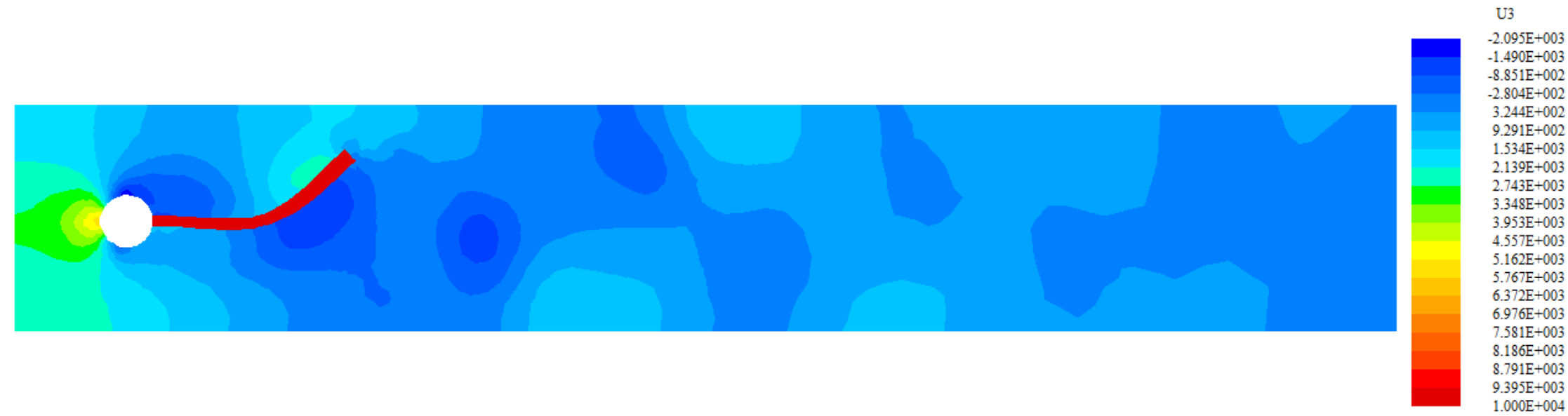}
    }
    \subfigure[\scriptsize $(\hat{U}, \mu_s) = (1.53, 0.51 \times 10^6)$: ROM $p$]{
        \includegraphics[width = 0.4\textwidth]{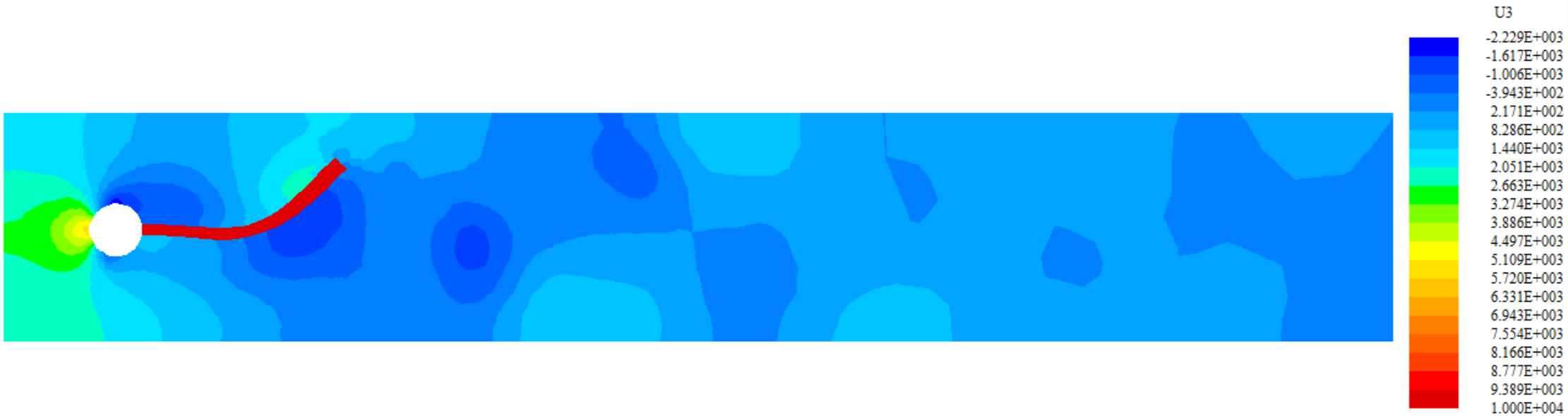}
    }
    \subfigure[\scriptsize $(\hat{U}, \mu_s) = (1.55, 0.49 \times 10^6)$: FOM $p$]{
        \includegraphics[width = 0.4\textwidth]{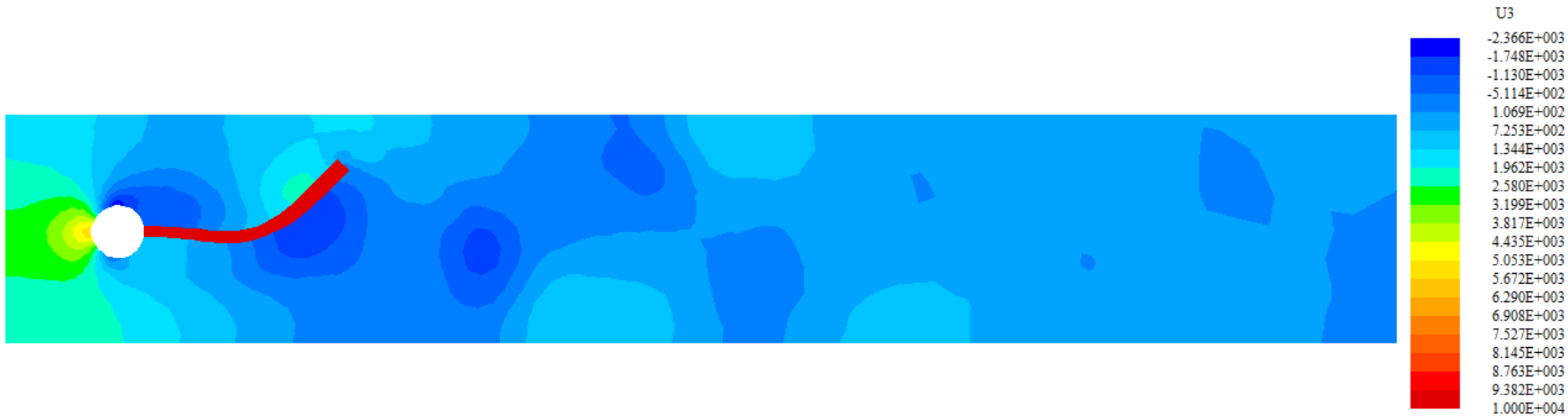}
    }
    \subfigure[\scriptsize $(\hat{U}, \mu_s) = (1.55, 0.49 \times 10^6)$: ROM $p$]{
        \includegraphics[width = 0.4\textwidth]{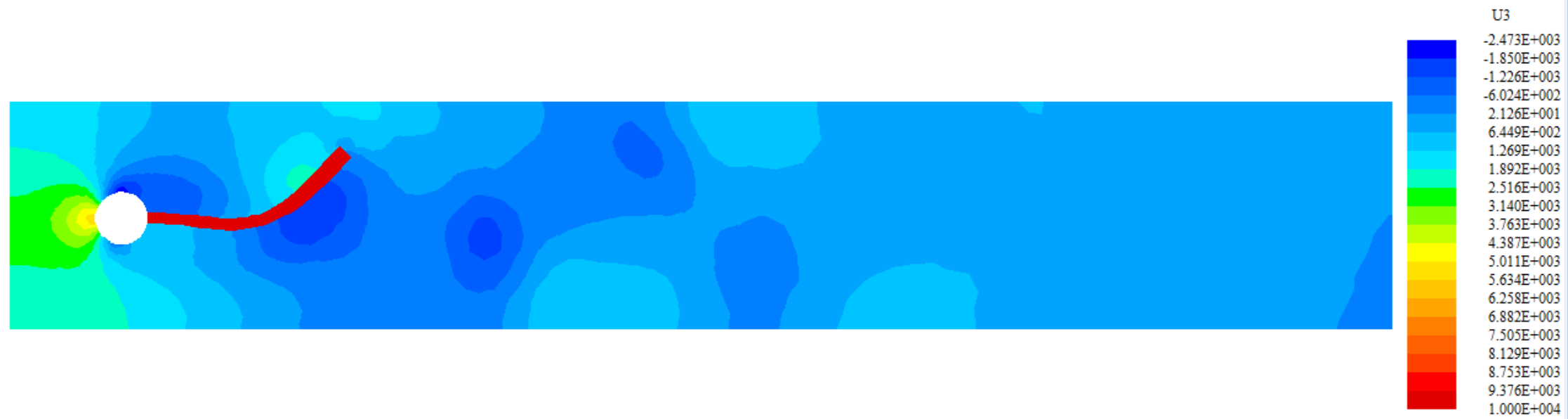}
    }
    \subfigure[\scriptsize $(\hat{U}, \mu_s) = (1.59, 0.59 \times 10^6)$: FOM $p$]{
        \includegraphics[width = 0.4\textwidth]{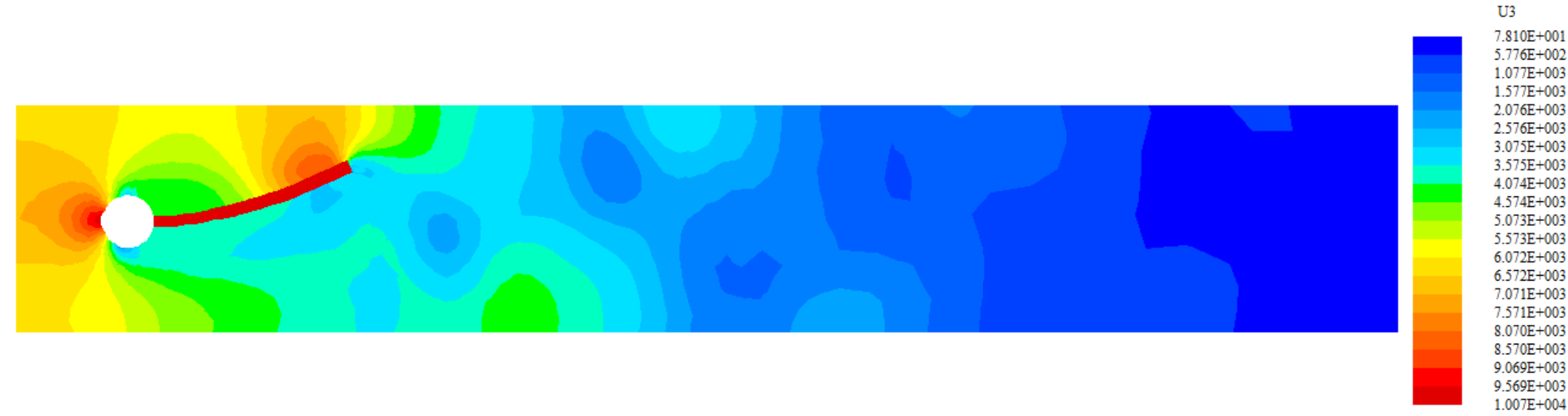}
    }
    \subfigure[\scriptsize $(\hat{U}, \mu_s) = (1.59, 0.59 \times 10^6)$: ROM $p$]{
        \includegraphics[width = 0.4\textwidth]{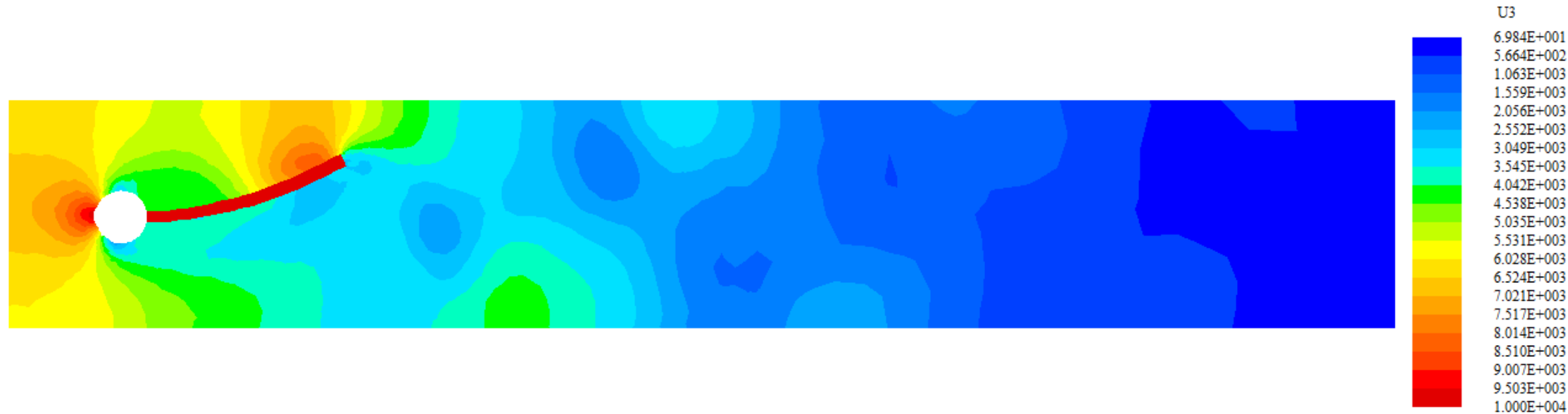}
    }
    \caption{Spatial distribution of pressure $p$ at $T=15$s for test parameter combinations.}
    \label{spatial_p_test}
\end{figure}

To assess the ROM's accuracy in capturing spatial field distributions, we compare FOM and ROM solutions at the final simulation time $T=15$s. Figures \ref{spatial_u1_train}--\ref{spatial_p_test} present comprehensive spatial distribution comparisons for all parameter combinations, including both training and test sets.

The spatial distribution comparisons reveal that the ROM accurately captures the complex flow features, including:
\begin{itemize}
    \item \textbf{Vortex structures}: The ROM preserves the vortex patterns in the fluid domain, which are crucial for understanding the FSI dynamics.
    \item \textbf{Velocity gradients}: Both horizontal ($u_1$) and vertical ($u_2$) velocity components are accurately reproduced, maintaining the correct boundary layer behavior and flow separation characteristics.
    \item \textbf{Pressure distributions}: The pressure field, which plays a critical role in the fluid-structure coupling, is captured with high fidelity, preserving the pressure gradients and their spatial variations.
    \item \textbf{Interface dynamics}: The ROM maintains accurate representation of the fluid-structure interface, where the coupling effects are most pronounced.
\end{itemize}

To provide a comprehensive view of the ROM's spatial accuracy, we present spatial distribution comparisons for all parameter combinations at $T=15$s. Figures \ref{spatial_u1_train}--\ref{spatial_p_train} show the FOM and ROM solutions for all training parameter combinations, while Figures \ref{spatial_u1_test}--\ref{spatial_p_test} show the results for all test parameter combinations.

The comprehensive spatial distribution comparisons demonstrate that the ROM accurately reproduces the FOM solutions across all parameter combinations, both in the training set and test set. The visual comparisons show excellent agreement in flow patterns, vortex structures, and pressure distributions, validating the ROM's capability to capture the complex spatial features of FSI problems.

\subsubsection{Computational Efficiency}

The primary advantage of the ROM approach lies in its computational efficiency. While each FOM simulation requires approximately 40 hours of CPU time, the ROM online phase can provide predictions in 0.97 seconds per parameter combination. This represents a speedup factor of approximately 100,000 times, making the ROM approach highly suitable for real-time applications, parameter studies, and uncertainty quantification. An important feature of the online phase is its inherent parallelizability. Once the PDE parameters $\Theta^*$ are given, all computations for different variables in different time steps can be performed in parallel, since the neural network evaluations and basis reconstructions at each time step are independent of each other. This parallel execution capability further enhances the computational efficiency of the ROM approach, allowing for even faster predictions when parallel computing resources are available.

The offline cost consists of four distinct parts:
\begin{itemize}
    \item generating the 121 FOM trajectories, at approximately 40 hours per trajectory: $121\times40=4{,}840$ CPU-hours;
    \item greedy basis search and construction, averaging 0.5 hours for each variable--segment job;
    \item preparation of the coefficient-training data, averaging 0.1 hours for each variable--segment job; and
    \item OGA training, averaging 1.5 hours for each variable--segment job.
\end{itemize}
There are $5\times55=275$ independent variable--segment jobs.  Hence the aggregate computational work is
\begin{equation}
4{,}840+275(0.5+0.1+1.5)=5{,}417.5\ \text{CPU-hours}.
\end{equation}
The last three stages are parallel over both variables and time segments.  With enough workers, their wall-clock contributions are 0.5, 0.1, and 1.5 hours, respectively.  Under the stated parallelization and a sequential generation of the 121 FOM trajectories, the offline critical-path time is
\begin{equation}
4{,}840+0.5+0.1+1.5=4{,}842.1\ \text{hours}.
\end{equation}
 If the 121 mutually independent FOM simulations are also executed concurrently, the idealized offline wall-clock time is approximately $40+0.5+0.1+1.5=42.1$ hours.  We report both aggregate CPU work and wall-clock time because the latter depends on the available parallel resources.

Table \ref{rom_fom_comparison} provides a comprehensive comparison between FOM and ROM approaches, highlighting the key advantages of the proposed method.

\begin{table}
    \centering
    \begingroup
    \small
    \setlength{\tabcolsep}{4pt}
    \renewcommand{\arraystretch}{0.92}
    \begin{tabular}{|c|>{\raggedright\arraybackslash}p{0.72\textwidth}|}
    \hline
    \textbf{Model} & \textbf{Specification} \\
    \hline
    \multirow{12}{*}{\textbf{FOM}} & \textbf{Total DOFs (all time steps):} \\
    & \quad Horizontal fluid velocity ($u_f$): $5,848 \times 15,000 = 87,720,000$ \\
    & \quad Vertical fluid velocity ($v_f$): $5,848 \times 15,000 = 87,720,000$ \\
    & \quad Fluid pressure ($p_f$): $6,170 \times 15,000 = 92,550,000$ \\
    & \quad Horizontal structure displacement ($u_s$): $242 \times 15,000 = 3,630,000$ \\
    & \quad Vertical structure displacement ($v_s$): $242 \times 15,000 = 3,630,000$ \\
    & \quad \textbf{Total DOFs:} $18,350 \times 15,000 = 275,250,000$ \\
    & \textbf{Time steps:} 15,000 \\
    & \textbf{Time segments:} 1 \\
    & \textbf{Computational cost (Online phase):} \\
    & \quad Single simulation: approximately 40 hours \\
    & \textbf{Memory requirements:} \\
    & \quad Per time step: 1.4 MB \\
    & \quad Total storage: 21 GB \\
    \hline
    \multirow{17}{*}{\textbf{ROM}} & \textbf{Total DOFs (all time steps):} \\
    & \quad Horizontal fluid velocity ($u_f$): $500 + 600 \times 4 + 600 \times 50 = 32,900$ \\
    & \quad Vertical fluid velocity ($v_f$): $500 + 600 \times 4 + 600 \times 50 = 32,900$ \\
    & \quad Fluid pressure ($p_f$): $1,500 + 1,800 \times 4 + 2,000 \times 50 = 108,700$ \\
    & \quad Horizontal structure displacement ($u_s$): $300 + 400 \times 4 + 400 \times 50 = 21,900$ \\
    & \quad Vertical structure displacement ($v_s$): $300 + 400 \times 4 + 400 \times 50 = 21,900$ \\
    & \quad \textbf{Total DOFs:} $218,300$ \\
    & \textbf{Time steps:} 15,000 \\
    & \textbf{Time segments:} 55 \\
    & \textbf{Computational cost:} \\
    & Offline aggregate work: $5{,}417.5$ CPU-hours\\
    & Offline wall-clock (parallel postprocessing): $4{,}842.1$ hours\\
    & Idealized wall-clock (121 parallel FOM runs): $42.1$ hours\\
    & \quad Online phase: \\
    & \quad \quad Single simulation: 0.97 seconds \\
    & \textbf{Speedup factor (online):} approximately 100000 times \\
    & \textbf{Accuracy:} \\
    & \quad Spatial-temporal Relative $L^2$ error (total): 0.0209 \\
    & \quad Spatial-temporal Relative $L^2$ error (fluid): 0.0087 \\
    & \quad Spatial-temporal Relative $L^2$ error (structure): 0.0016 \\
    & \quad Spatial-temporal Relative $L^2$ error (pressure): 0.0190 \\
    & \textbf{Memory requirements:} \\
    & \quad Per time step: 0.03 MB \\
    & \quad Total storage: 0.5 GB \\
    \hline
    \end{tabular}
    \endgroup
    \caption{Comprehensive comparison between FOM and ROM approaches. The ROM achieves several orders of magnitude speedup while maintaining high accuracy, making it suitable for real-time applications and extensive parametric studies.}
    \label{rom_fom_comparison}
\end{table}

The comparison table demonstrates several key advantages of the ROM approach:

\begin{itemize}
    \item \textbf{Dimensionality reduction}: The ROM reduces the total system dimension from 275,250,000 DOFs (FOM) to 218,300 DOFs (ROM) across all time steps, representing a reduction factor of approximately 1,260$\times$. The basis dimensions vary across time segments and variables, with the number of basis functions per variable per segment ranging from 300 to 2,000 depending on the solution complexity.
    \item \textbf{Computational speedup}: The online phase achieves a speedup factor of approximately 100,000 times, reducing simulation time from approximately 40 hours (FOM) to 0.97 seconds (ROM) per parameter combination.
    \item \textbf{Memory efficiency}: The ROM reduces memory requirements by a factor of approximately 42$\times$, from 21 GB (FOM) to 0.5 GB (ROM) for the complete simulation dataset. The per time step memory requirement is reduced from 1.4 MB (FOM) to 0.03 MB (ROM).
    \item \textbf{Accuracy preservation}: Despite the massive computational savings, the ROM maintains high accuracy with spatial-temporal relative $L^2$ errors of 0.0209 (total), 0.0087 (fluid), 0.0016 (structure), and 0.0190 (pressure), demonstrating excellent accuracy preservation across all variables and parameter combinations.
    \item \textbf{Enhanced capabilities}: The ROM enables real-time predictions and efficient uncertainty quantification, which are impractical with the FOM approach. The parallelizability of the online phase further enhances computational efficiency when parallel computing resources are available.
\end{itemize}

These advantages make the ROM approach particularly valuable for applications requiring rapid parameter exploration, real-time decision making, and comprehensive uncertainty quantification in FSI problems.

\section{Conclusion}\label{s5}
This paper presents a novel non-intrusive reduced order modeling
approach for fluid-structure interaction problems that combines
greedy parameter-time selection with orthogonal greedy algorithm
training. The method addresses the computational challenges of
high-fidelity FSI simulations by constructing efficient surrogate
models without requiring access to or modification of the underlying
source code. Our approach employs time segmentation to handle
long-term simulations, where separate ROM bases are constructed for
each time segment. The greedy parameter-time selection strategy
iteratively identifies the most informative parameter-time pairs by
maximizing projection errors onto current reduced bases, while
Gram-Schmidt orthogonalization ensures stable and efficient basis
construction. The integration of shallow neural networks trained
using the orthogonal greedy algorithm provides theoretical
convergence guarantees and avoids common issues of overfitting and
local minima encountered in traditional neural network training.

Comprehensive numerical experiments on a two-dimensional FSI
benchmark problem demonstrate that our method achieves comparable
accuracy to the full-order model (with relative $L^2$ errors
consistently around 1-2\%) while providing significant computational
speedup (approximately 10$\times$ in the online phase). The method
effectively handles parameter perturbations in the online phase and
maintains dimension-independence for extension to high-dimensional
FSI problems. The integration of greedy parameter selection with
convergent machine learning methods represents a significant
advancement in non-intrusive surrogate modeling for complex
multiphysics systems, making it highly suitable for real-time
applications, extensive parametric studies, and uncertainty
quantification in FSI problems.

\section{Acknowledgment}\label{s6}
X. Xie was partially supported by National Natural Science
Foundation of China (12571434) and Open Research Project
of National Key Laboratory of Fundamental Algorithms and Models for
Engineering Simulation. P. Sun was partially supported by a grant
from the Simons Foundation (MPS-TSM-00706640).

\bibliographystyle{siam}
\bibliography{ref,FSI, publicationlist4Sun}
\end{document}